\documentclass[11pt]{amsart}
\usepackage[]{amsmath, amsthm, amsfonts, verbatim, amssymb}

\DeclareFontFamily{OT1}{rsfs}{}

\DeclareFontShape{OT1}{rsfs}{n}{it}{<-> rsfs10}{}
\DeclareMathAlphabet{\mathscr}{OT1}{rsfs}{n}{it}

\def\cD{{\mathscr D}}
\def\cE{{\mathscrE}}

\def\ocM{\overline{\mathscr M}}
\def\oM{\overline M}

\def\cH{\mathscr H}

\def\cI{{\mathscr I}}
\def\cR{{\mathscr R}}
\def\cG{{\mathscr G}}

\def\cT{{\mathscr T}}
\def\bP{{\mathbb P}}
\def\bR{{\mathbb R}}
\def\bZ{{\mathbb Z}}
\def\bQ{{\mathbb Q}}
\def\bC{{\mathbb C}}
\def\({{\Big(}}
\def\){{\Big)}}

\def\cS{{\mathscr S}}
\def\cE{{\mathscr E}}
\def\cP{{\mathscr P}}
\def\cC{{\mathcal C}}
\begin{document}
\newtheorem {theo}{Theorem}
\newtheorem {coro}{Corollary}
\newtheorem {lemm}{Lemma}
\newtheorem {rem}{Remark}
\newtheorem {defi}{Definition}
\newtheorem {ques}{Question}
\newtheorem {prop}{Proposition}
\def\spb{\smallpagebreak}
\def\mpb{\vskip 0.5truecm}
\def\bpb{\vskip 1truecm}
\def\wtM{\widetilde M}
\def\tM{\widetilde M}
\def\wtN{\widetilde N}
\def\tN{\widetilde N}
\def\tC{\widetilde C}
\def\X{\widetilde X}
\def\tY{\widetilde Y}
\def\tP{\widetilde P}
\def\ti{\widetilde \iota}
\def\bs{\bigskip}
\def\ms{\medskip}
\def\ni{\noindent}
\def\td{\nabla}
\def\pd{\partial}
\def\hol{$\text{hol}\,$}
\def\Log{\mbox{Log}}
\def\bfQ{{\bf Q}}
\def\Todd{\mbox{Todd}}
\def\bP{{\bf P}}
\def\dxi{d x^i}
\def\dxj{d x^j}
\def\dyi{d y^i}
\def\dyj{d y^j}
\def\dzi{d z^I}
\def\dzj{d z^J}
\def\ozi{d{\overline z}^I}
\def\ozj{d{\overline z}^J}
\def\oz1{d{\overline z}^1}
\def\oz2{d{\overline z}^2}
\def\oz3{d{\overline z}^3}
\def\sI{\sqrt{-1}}
\def\hol{$\text{hol}\,$}
\def\ok{\overline k}
\def\ol{\overline l}
\def\oJ{\overline J}
\def\oT{\overline T}
\def\oS{\overline S}
\def\oV{\overline V}
\def\oW{\overline W}
\def\oI{\overline I}
\def\oK{\overline K}
\def\oL{\overline L}
\def\oj{\overline j}
\def\oi{\overline i}
\def\ok{\overline k}
\def\oz{\overline z}
\def\om{\overline mu}
\def\on{\overline nu}
\def\oa{\overline \alpha}
\def\ob{\overline \beta}
\def\of{\overline f}
\def\og{\overline \gamma}
\def\ogamma{\overline \gamma}
\def\odelta{\overline \delta}
\def\otheta{\overline \theta}
\def\ophi{\overline \phi}
\def\opd{\overline \partial}
\def\oA{\overline A} 
\def\oB{\overline B}
\def\oC{\overline C}
\def\op{\overline D}
\def\oIq1{\oI_1\cdots\oI_{q-1}}
\def\oIq2{\oI_1\cdots\oI_{q-2}}
\def\op{\overline \partial}
\def\ua{{\underline {a}}}
\def\us{{\underline {\sigma}}}
\def\tor{{\mbox{tor}}}
\def\vol{{\mbox{vol}}}
\def\rank{{\mbox{rank}}}
\def\bp{{\bf p}}
\def\bk{{\bf k}}
\def\a{{\alpha}}
\def\tchi{\widetilde{\chi}}
\def\fv{\mathfrak{v}}
\def\fD{\mathfrak{D}}
\title[Fake projective planes]
{Fake projective planes}
 \maketitle
\vskip-3mm
{\centerline{\sc Gopal Prasad and Sai-Kee Yeung}}
\vskip5mm

\centerline{\it Dedicated to David Mumford}

\vskip7mm

{\bf This is a revised version of the paper with the same title published in Inventiones Math.\,{\bf 168}(2007). It  incorporates corrections and additions given in the addendum published in Inventiones Math.\,{\bf 182}(2010).}
\vskip5mm

\begin{center}
{\bf 1. Introduction}  
\end{center}

\vskip4mm

\ni{\bf 1.1.} A fake projective plane is a smooth compact complex surface which is not the complex projective plane
but has the same Betti numbers as the complex projective plane.  Such a surface is known to be projective algebraic and it is the quotient of 
the (open) unit ball $B$ in $\bC^2$ ($B$ is the symmetric space of ${\rm PU}(2,1)$) by a torsion-free  cocompact discrete 
subgroup of 
${\rm PU}(2,1)$ whose Euler-Poincar\'e characteristic is $3$. These surfaces have the smallest Euler-Poincar\'e characteristic among all smooth surfaces of
general type.  The first fake projective plane was constructed by 
David Mumford [Mu] using $p$-adic uniformization, 
and later two more examples were found by 
M.\,Ishida and F.\,Kato
in [IK] using a similar method. In [Ke] JongHae Keum has constructed 
an example which is birational to a cyclic cover of degree $7$ of a Dolgachev surface 
(see 5.15 below). It is known that there are only finitely many fake projective planes ([Mu]), and an important  problem in complex algebraic 
geometry is to determine them all. 
\vskip1mm

It is proved in [Kl] and [Y] that the fundamental group of a fake projective 
plane is a torsion-free cocompact arithmetic subgroup of ${\rm PU}(2,1)$. It follows from 
Mostow's strong rigidity theorem ([Mo]) that the fundamental group of a fake projective plane determines it uniquely up to isometry. In this 
paper we will make use of the volume formula 
of [P], together with some number theoretic estimates, 
to list  all torsion-free cocompact arithmetic subgroups (of ${\rm PU}(2,1)$) whose Euler-Poincar\'e characteristic is $3$, see \S\S 5, 8 and 9. This  list of course contains the fundamental groups of all     
fake projective planes. It provides several new examples of fake projective planes. In 
fact, we show that there are exactly {\it twenty eight}  
distinct nonempty classes of fake projective planes (see 1.4--1.5 below). We 
obtain these fake projective planes as quotients of the ball $B$
by explicitly given torsion-free cocompact arithmetic subgroups of either ${\rm SU}(2,1)$ or 
${\rm PU}(2,1)$.  In \S 10, we use this explicit description of their 
fundamental groups to prove that for any fake projective plane $P$ occurring in these twenty eight classes, 
$H_1(P,\bZ )$ is nonzero. We also prove that if $P$ is not a fake projective plane arising from the pair $\cC_2$ or $\cC_{18}$  ($\cC_i$ as in 8.2), its fundamental group embeds in ${\rm SU}(2,1)$ (Proposition 10.3). Using computer-assisted group theoretic computations, Donald Cartwright and Tim Steger have shown recently that the fundamental group of every fake projective plane arising from the pair $\cC_2$ also embeds in ${\rm SU}(2,1)$.  For any fake projective plane $P$ for which this holds, the canonical line 
bundle $K_P$ is divisible by $3$, i.\,e., there is a holomorphic line bundle $L$ on $P$ such that $K_P = 3L$ (see 10.4). In 10.5 we show that any fake-projective plane can be embedded in $\bP^{14}_{\bC}$ as a smooth complex surface of degree $49$. 
\vskip2mm
 
We will now present a brief outline of our methods and results. We begin by giving a description of the forms of ${\rm SU}(2,1)$ over number fields used in this paper.
\vskip2mm

\ni{\bf 1.2.} Let $k$ be a real number field, $v_o$ be a real place of $k$,  and $G$ be a simple simply connected algebraic $k$-group such that $G(k_{v_o})\cong {\rm SU}(2,1)$, and for all other archimedean places $v$ of $k$, $G(k_v)\cong {\rm SU}(3)$. From the
description of absolutely simple simply connected groups of type $^2A_2$ (see, 
for example, [Ti1]), we see that $k$ is 
totally real, and there is a totally complex quadratic extension $\ell$ 
of $k,$ a division algebra $\cD$ of degree $n|3,$ with center
$\ell,$ $\cD$ given with an involution $\sigma$ of the second kind such
that $k=\{x\in\ell\ |\ x=\sigma(x)\}$, and a nondegenerate 
hermitian form $h$ on $\cD^{3/n}$ defined in terms of the involution $\sigma,$ 
such that $G$ is the special unitary group ${\rm SU}(h)$ of $h.$ If $\cD = \ell$, $h$ is a hermitian form on $\ell^3$ and its determinant ${\rm det}(h)$ is defined in the usual way. On the other hand, if $\cD$ is a cubic division algebra, then  $h(x,y)=  \sigma(x)ay$, for some $a\in {\cD}^{\sigma}$ and the determinant ${\rm det}(h)$ of $h$ is by definition  ${\rm{Nrd}}(a)$ modulo $N_{\ell/k}(\ell^{\times})$.

\vskip1mm

Let $k$, $\ell$, $\cD$ be as above. We will now show that the $k$-group $G$ is uniquely determined, up to a $k$-isomorphism, by $\cD$ (i.\,e., the $k$-isomorphism class of $G$ does not depend on the choice  of the involution $\sigma$ and the hermitian form $h$ on ${\cD}^{3/n}$). Let $\sigma$ be an involution of $\cD$ of the second kind with $k= \{ x\in \ell \, |\, x =\sigma (x)\} $. Let $h$ be a hermitian form on ${\cD}^{3/n}$.   For $x\in k^{\times}$, $xh$ is again an hermitian form on ${\cD}^{3/n}$, and ${\rm{det}}(xh) =x^3{\rm{det}}(h)$. Now since $N_{\ell /k}(\ell^{\times} )\supset {k^{\times}}^2$, ${\rm{det} }({\rm{det}}(h)h)$, as an element of $k^{\times}/N_{\ell /k}(\ell^{\times})$, is $1$.  Moreover, ${\rm{SU}}(h) ={\rm{SU}}({\rm{det}}(h)h)$. Hence, it would suffice to work with hermitian forms of determinant $1$. 
\vskip1mm

If ${\cD} = \ell$, and $h$ is a hermitian form on $\ell^3$ of determinant $1$ such that the group ${\rm{SU}}(h)$ is isotropic at $v_o$, and is anisotropic at all other real places of $k$ (or, equivalently, $h$ is indefinite at $v_o$, and definite at all other real places), then being of determinant $1$, its signature (or index) at $v_o$ is $-1$, and at all other real places of $k$ it is $3$. Corollary 6.6 of [Sc, Chap.\,10] implies that any two such hermitian forms on $\ell^3$ are isometric, and hence they determine a unique $G$ up to a $k$-isomorphism.  
\vskip1mm

Now let us assume that $\cD$ is a cubic division algebra with center $\ell$, $\sigma$ an involution of the second kind such that for the hermitian form $h_0$ on $\cD$ defined by $h_0 (x,y) =\sigma(x)y$, the group ${\rm{SU}}(h_0 )$ is isotropic at $v_o$, and is anisotropic at every other real place of $k$. For $z\in{\cD}^{\times}$, let ${\rm{Int}}(z)$ denote the automorphism $x\mapsto zxz^{-1}$ of ${\cD}$.  Let ${\cD}^{\sigma} =\{ z\in {\cD}\, |\, \sigma(z) = z\}$. Then for all $z\in{{\cD}^{\sigma}}^{\times}$, ${\rm{Int}}(z)\cdot\sigma$ is again an involution of ${\cD}$ of the second kind, and any involution of $\cD$ of the second kind is of this form. Now for $z\in {\cD}^{\sigma}$, given an hermitian form $h'$ on ${\cD}$ with respect to the involution ${\rm{Int}}(z)\cdot \sigma$, the form $h= z^{-1}h'$ is a hermitian form on ${\cD}$ with respect to $\sigma$, and ${\rm{SU}}(h') = {\rm{SU}}(h)$. Therefore, to determine all the special unitary groups we are interested in, it is enough to work just with the involution $\sigma$, and to consider all hermitian forms $h$ on $\cD$, with respect to $\sigma$, of determinant $1$, such that the group ${\rm{SU}}(h)$ is isotropic at $v_o$, and is anisotropic at all other real places of $k$.  Let $h$ be such a hermitian form. Then $h(x,y)=  \sigma(x)ay$, for some $a\in {\cD}^{\sigma}$, and ${\rm det}(h) =1$ 
so ${\rm Nrd}(a)\in N_{\ell/k}(\ell^{\times})$.   As the elements of $N_{\ell/k}(\ell^{\times})$ are positive at all real places of $k$, we see that the signatures of $h$ and $h_0$ are equal at every real place of $k$. Corollary 6.6 of [Sc, Chap.\,10] again implies that the hermitian forms $h$ and $h_0$ are isometric. Hence, ${\rm{SU}}(h)$ is $k$-isomorphic to ${\rm{SU}}(h_0)$. Thus we have shown that $ \cD$ determines a unique $k$-form $G$ of ${\rm{SU}}(2,1)$, up to a $k$-isomorphism, namely ${\rm{SU}}(h_0)$, with the desired behavior at the real places of $k$. For any commutative $k$-algebra $A$, we will denote the $A$-linear extension of $\sigma$ to $A\otimes_k\cD$ also by $\sigma$.  The group $G(A)$  of $A$-rational points of this $G$ is 
$$G(A) = \{ g\in {\rm GL}_{1,\cD}(A)=({A\otimes_k{\cD}})^{\times}\, |\, g\sigma(g) =1 \ {\rm{and}}\ {\rm{Nrd}}(g)= 1\}.$$ 

Let $\cD$ and the involution $\sigma$ be as in the previous paragraph. Let ${\cD}^o$ be the opposite of ${\cD}$. Then the involution $\sigma$ is also an involution of ${\cD}^o$. The pair $({\cD}^o,\sigma)$ determines a $k$-form of ${\rm{SU}}(2,1)$ which is clearly 
$k$-isomorphic to the one determined by the pair $(\cD,\sigma)$.          

\vskip1mm

In the sequel, the adjoint group of $G$ will be denoted by $\overline G$, and $\varphi$ will denote the natural isogeny $G\to {\overline G}$. It is known that if $\cD$ is a cubic division algebra, then ${\mathrm{Aut}}(G)(k) = \overline{G}(k)$, i.e., any $k$-rational automorphism 
of  $G$ (and so also of $\overline{G}$) is inner.
\vskip2mm

\ni{\bf  1.3.} Let $\Pi$ be a torsion-free cocompact arithmetic subgroup 
of ${\rm PU}(2,1)$ whose Euler-Poincar\'e characteristic is $3$. The fundamental group of a fake projective plane is such a subgroup. Let $\varphi :\, {\rm SU}(2,1)\rightarrow {\rm PU}(2,1)$ be the 
natural surjective homomorphism. The kernel of $\varphi$ is the center of ${\rm SU}(2,1)$ 
which is a subgroup of order $3$. Let $\widetilde\Pi = \varphi^{-1}({\Pi})$. Then $\widetilde\Pi$ is a cocompact arithmetic subgroup of 
${\rm SU}(2,1)$. The orbifold Euler-Poincar\'e characteristic $\chi({\widetilde\Pi})$ of $\widetilde\Pi$ 
(i.\,e.,\,\,the Euler-Poincar\'e characteristic in the 
sense of C.\,T.\,C.\,Wall, cf.\,[Se1], \S 1.8) is $1$. Hence, the orbifold Euler-Poincar\'e characteristic of any discrete  subgroup of ${\rm SU}(2,1)$ containing $\widetilde\Pi$ is a reciprocal integer.
\vskip1mm

Let $k$ be the number field and $G$ be the $k$-form of ${\rm SU}(2,1)$ associated with 
the arithmetic subgroup $\widetilde\Pi$. The field $k$ is generated by the traces, in the 
adjoint representation of ${\rm PU}(2,1)$, of the elements in $\Pi$, 
and $G$ is a simple simply connected algebraic  
$k$-group such that for a real place, say $v_o$, of $k$, 
$G(k_{v_o})\cong {\rm SU}(2,1)$, and for all archimedean 
places $v \ne v_o$, $G(k_v)$ is isomorphic to the compact Lie 
group ${\rm SU}(3)$, and $\widetilde\Pi$ is commensurable with ${\widetilde\Pi}\cap G(k)$. Throughout this paper we will use the description of $G$ and $\overline G$  given in 1.2. In particular, $\ell$, $\cD$ and $h$ are as in there. 

\vskip1mm

Let $V_f$ (resp.\,\,$V_{\infty}$) be the set of nonarchimedean (resp.\,\,archimedean) places of $k$. Let $\cR_{\ell}$ be the set of nonarchimedean places of $k$ which ramify in $\ell$. The $k$-algebra of finite ad\`eles of $k$, i.\,e., the restricted direct product of the $k_v$, $v\in V_f$, will be denoted by $A_f$. 

The image $\Pi$ of $\widetilde\Pi$ in ${\overline G}(k_{v_o})$ is actually contained in  
${\overline G}(k)$ ([BP], 1.2). For all $v\in V_f$, we fix a {\it parahoric} subgroup $P_v$ of $G(k_v)$ which is {\it minimal} among the parahoric subgroups of $G(k_v)$ normalized by $\Pi$. Then 
$\prod_{v\in V_f}P_v$ is an open subgroup of $G(A_f)$, see [BP], \S 1. Hence, $\Lambda := G(k)\cap\prod_{v\in V_f}P_v$ is a {\it principal} arithmetic subgroup ([P], 3.4) 
which is  normalized by $\Pi$, and therefore also by $\widetilde\Pi$. Let $\Gamma$ be the normalizer of $\Lambda$ in $G(k_{v_o})$, and $\overline\Gamma$ be its image in ${\overline G}(k_{v_o})$. Then ${\overline\Gamma}\subset {\overline G}(k)$ ([BP], 1.2). As the normalizer of $\Lambda$ in $G(k)$ equals $\Lambda$, $\Gamma\cap G(k) = \Lambda$. Since $\Gamma$ contains $\widetilde\Pi$, its orbifold Euler-Poincar\'e characteristic $\chi(\Gamma)$ is a reciprocal integer.  
\vskip1mm

In terms 
of the normalized Haar-measure $\mu$ on $G(k_{v_o})$ used in [P]
and [BP], $\chi(\Gamma)=3\mu(G(k_{v_o})/\Gamma)$ (see
\S4 of [BP], note that the compact dual of the symmetric space $B$ of 
$G(k_{v_o})\cong{\rm SU}(2,1)$ is $\bP^2_{\bC}$, and the Euler-Poincar\'e characteristic 
of $\bP^2_{\bC}$ is $3$). Thus the condition that $\chi(\Gamma)$ is a reciprocal integer is 
equivalent to 
the condition that the covolume $\mu(G(k_{v_o})/\Gamma)$, of $\Gamma$,  is one third of a 
reciprocal integer; in particular, $\mu(G(k_{v_o})/\Gamma)\leqslant 1/3$. Also, 
$\chi(\Gamma) =3\mu(G(k_{v_o})/\Gamma)= 
3\mu(G(k_{v_o})/\Lambda )/[\Gamma : \Lambda]$, and the 
volume formula of [P] can be used to compute $\mu(G(k_{v_o})/\Lambda)$ precisely, see 2.4 below. 
Proposition 2.9 of [BP] implies that $[\Gamma : \Lambda]$ is 
a power of $3$. Now we see that if $\chi(\Gamma)$ is a reciprocal integer, 
then the numerator of the rational number $\mu(G(k_{v_o})/\Lambda)$ must be a power of $3$. 
\vskip2mm

\ni{\bf 1.4.} In \S\S  4--5, and 7--9, we will determine all $k$, $\ell$, $\cD$, 
simple simply connected algebraic $k$-groups $G$ so that for a real place $v_o$ of $k$, $G(k_{v_o})\cong {\rm SU}(2,1)$, for all archimedean $v\ne v_o$, $G(k_v)\cong {\rm SU}(3)$, and (up to conjugation by an element of ${\overline G}(k)$) all collections $(P_v)_{v\in V_f}$ of parahoric subgroups $P_v$ of $G(k_v)$ such that $(i)$ $\prod_{v\in V_f}P_v$ is an open subgroup of $G(A_f)$, $(ii)$ the principal arithmetic subgroup $\Lambda: = G(k)\cap\prod_{v\in V_f}P_v$ considered as a (discrete) subgroup of $G(k_{v_o})$ is cocompact (by Godement compactness criterion, this is equivalent to the condition that $G$ is anisotropic over $k$), and $(iii)$ the image $\overline\Gamma$ in ${\overline G}(k_{v_o})$ of the normalizer $\Gamma$ of $\Lambda$ in $G(k_{v_o})$ contains a torsion-free subgroup $\Pi$ of finite index whose Euler-Poincar\'e characteristic is $3$. Then the orbifold 
Euler-Poincar\'e characteristic of $\Gamma$ is a reciprocal integer.
\vskip1mm

\ni{\bf 1.5.} Let us first consider the case where $\cD =\ell$. Then $h$ is a nondegenerate  hermitian form on $\ell^3$ (defined in terms of the nontrivial automorphism of $\ell /k$) which is indefinite at $v_o$ and definite at all other real places of $k$. Let $G ={\rm SU}(h)$, and $\overline G$ be its adjoint group. We prove below (Proposition 8.8) that  if ${\overline G}(k_{v_o})$ contains a torsion-free cocompact arithmetic subgroup $\Pi$ with $\chi(\Pi) = 3$, then, in the notation of 8.2,  $(k,\ell)$ must be one of the following five: $\cC_1$, $\cC_8$, $\cC_{11}$, $\cC_{18}$, and $\cC_{21}$. Using  quite sophisticated computer-assisted group theoretic computations, Cartwright and Steger have recently shown (see [CS2]) that for $(k,\ell)$ any of these five pairs the fundamental group of a fake projective plane cannot be an arithmetic subgroup of ${\rm PU}(h)$. 
\vskip1mm

Cartwright and Steger have also shown that there exists a rather unexpected smooth projective complex algebraic surface,  uniformized by the complex 2-ball, whose  fundamental group is a cocompact torsion-free arithmetic subgroup of ${\rm PU}(h)$, $h$ as above, with $(k,\ell) = \cC_{11}$ $=(\bQ(\sqrt{3}),\bQ(\zeta_{12}) $), and whose Euler-Poincar\'e characteristic is $3$ but the first Betti-number is nonzero (it is actually $2$); we name this surface the ``Cartwright-Steger surface''. Since the first Betti-number of this surface is nonzero, it admits $n$-sheeted covers for every positive integer $n$. The 
Euler-Poincar\'e charactersitic of such a cover is $3n$. 
\vskip2mm 

\ni{\bf 1.6.}  In view of the result mentioned in the first paragraph of 1.5, we will assume in the rest of this section that $\cD\ne \ell$.  We will prove that (up to natural equivalence) there are exactly {\it twenty eight} distinct $\{k,\ell, G, (P_v)_{{v\in V_f}}\}$ satisfying the conditions mentioned in 1.4. Each of these twenty eight sets determines a unique (up to isomorphism) principal arithmetic subgroup $\Lambda$  ($=G(k)\cap \prod_{{v\in V_f}} P_v$), which in turn determines a unique arithmetic subgroup $\overline\Gamma$ of $\overline{G}(k_{v_o})$ (recall that $\overline\Gamma$ is the image in $\overline{G}(k_{v_o})$ of the normalizer $\Gamma$ of $\Lambda$ in $G(k_{v_o})$). For eighteen of these twenty eight, $k =\bQ$, see Sect.\,5; and there are two with $k=\bQ(\sqrt{2})$, two with $k=\bQ(\sqrt{5})$, and three each with $k=\bQ(\sqrt{6})$ and $k =\bQ(\sqrt{7})$, see Sect.\,9. The pair $(k,\ell) = (\bQ, \bQ(\sqrt{-1}))$ gives three, the pair $(\bQ,\bQ(\sqrt{-2}))$ gives three, the pair $(\bQ,\bQ(\sqrt{-7}))$ gives six, the pair $(\bQ, \bQ(\sqrt{-15}))$ gives four,  and the pair $(\bQ, \bQ(\sqrt{-23}))$ gives two classes of fake projective planes.
\vskip2mm

\ni{\bf 1.7.} If $\Pi$, $\Lambda$, $\Gamma$, and the parahoric subgroups $P_v$ are as in 1.3, then for $v\in V_f$, since $P_v$ was assumed to be minimal among the parahoric subgroups of $G(k_v)$ normalized by $\Pi$, if for a $v$, $P_v$ is maximal, then it is the {\it unique} parahoric subgroup of $G(k_v)$ normalized by $\Pi$. It will turn out that for every  $v\in V_f$, $P_v$ appearing in 1.3 is a maximal parahoric subgroup of $G(k_v)$ 
except when $(k,\ell )$ is either $(\bQ, \bQ(\sqrt{-1}))$ or $(\bQ, \bQ(\sqrt{-2}))$ or $\cC_{18} = (\bQ (\sqrt{6}),\bQ (\sqrt{6},\zeta_3 ))$,  in which cases $P_v$ is non-maximal for at most one $v$.

\vskip2mm

\ni{\bf 1.8.} We will now describe the class of fake projective planes associated to each of the twenty eight $\Gamma$'s of 1.6. The orbifold Euler-Poincar\'e characteristic $\chi({\overline\Gamma})$ of $\overline\Gamma$ equals $3\chi(\Gamma)= 3\chi(\Lambda)/[\Gamma : \Lambda]$, and we compute it precisely. Now if $\Pi$ is a torsion-free subgroup of $\overline\Gamma$ of index $3/\chi({\overline\Gamma})$, then $\chi(\Pi) =3$, and if, moreover, $H^1(\Pi, \bC)$ vanishes (or, equivalently, the abelianization $\Pi/[\Pi,\Pi]$ is finite), then  by Poincar\'e-duality, $H^3(\Pi, \bC)$ vanishes too, and hence, as $\chi(B/\Pi) = \chi (\Pi) =3$, $B/\Pi$ is  a fake projective plane.  We will show that each of the twenty eight $\overline\Gamma$ does contain a $\Pi$ with the desired properties. The class of fake projective planes given by $\Gamma$ (or $\overline\Gamma$) consists of the fake projective planes $B/\Pi$, where $\Pi$ is a torsion-free subgroup of $\overline\Gamma$ of index $3/\chi({\overline\Gamma})$ with  $\Pi/[\Pi, \Pi ]$ finite. 
\vskip1mm

We observe that in principle, for a given $\Gamma$, the subgroups $\Pi$ of $\overline\Gamma$ as above can all be determined in the following way: First find a ``small" presentation of $\overline\Gamma$ using a ``nice" fundamental domain in $B$ (maximal arithmetic subgroups tend to have small presentation), and use this presentation to list all torsion-free subgroups of index $3/\chi({\overline\Gamma})$ whose abelianization is finite. (Note that the computations below show that $3/\chi({\overline\Gamma})$ is quite small; in fact, it equals $1$, $3$, $9$ or $21$.) This has recently been carried out  by Cartwright and Steger using ingenious computer-assisted group theoretic computations. They have shown (see [[CS1])  that  the twenty eight classes of fake projective planes altogether contain {\it fifty} distinct
fake projective planes up to isometry with respect to the Poincar\'e metric.  Since each such fake projective plane as a Riemannian manifold supports
two distinct complex structures [KK, \S5], there are exactly {\it one hundred} fake  projective planes counted 
up to biholomorphism.
\vskip1mm

 Cartwright and Steger have  given explicit generators and relations for the fundamental group (which is a cocompact torsion-free arithmetic subgroup of ${\rm PU}(2,1)$) of each of the fake projective planes, determined their automorphism group, and computed their first homology with coefficients in $\bZ$.  They have  shown that the quotient of six of the fake projective planes by a subgroup of order $3$ of the automorphism group is a simply connected singular surface.  We propose to call these simply connected singular surfaces the ``Cartwright-Steger singular surfaces".
\vskip1mm

Cartwright and Steger  have also found that the fundamental group of eight of the one hundred fake projective planes do not admit an embedding into ${\rm SU}(2,1)$ as a discrete subgroup, hence the canonical line bundle of these fake projective planes is not divisible by $3$ in their Picard group.  (All such fake projective planes arise from the pair $\cC_{18}=(\bQ(\sqrt{6}),\bQ(\sqrt{6},\zeta_3))$.)
\vskip2mm

\vskip2mm

\ni{\bf 1.9.} The results of this paper show, in particular, that any arithmetic subgroup $\Gamma$ of ${\rm SU}(2,1)$, with 
$\chi (\Gamma)\leqslant 1$, must arise from a $k$-form $G$ of ${\rm
  SU}(2,1)$ as above, where the pair $(k,\ell)$ consists
of $k =\bQ$, and $\ell$ is one of the eleven imaginary quadratic fields
listed in Proposition 3.5, or $(k,\ell)$ is one of the forty pairs
$\cC_1$--$\cC_{40}$ described in 8.2. The covolumes, and hence the orbifold Euler-Poincar\'e characteristics,  of these arithmetic
subgroups can be computed using the volume formula given in 2.4 and
the values of $\mu$ given in Proposition 3.5 and in 8.2. The surfaces arising as the quotient of $B$ by one of these arithmetic subgroups are often singular. However, as they have a small orbifold Euler-Poincar\'e characteristic, they may have interesting geometric properties.

\vskip3mm

For a nice exposition of the results proved, and techniques employed, in this paper, see Bertrand R\'emy's Bourbaki report [R\'e].

\vskip4mm
\ni\begin{center}{\bf \S 2. Preliminaries}
\end{center}
\vskip2mm

A comprehensive survey of the basic definitions and the main results of the 
Bruhat--Tits theory of reductive groups over nonarchimedean local fields is given in [Ti2].
\vskip2mm

\ni{\bf 2.1.} Let the totally real number field $k$, and its totally complex quadratic extension $\ell$, a real place $v_o$ of $k$, 
and the $k$-form $G$ of ${\rm SU}(2,1)$ be as in 1.2. Throughout this paper, we will use the description of $G$ given in 1.2 and the notations introduced in \S 1. 
\vskip1mm

We shall say that a collection $(P_v)_{v\in V_f}$ of parahoric
subgroups $P_v$ of $G(k_v)$ is {\it coherent} if $\prod_{v\in V_f}P_v$
is an open subgroup of $G(A_f)$.   Let $U$ be a compact-open subgroup of $G(A_f)$, and $(P_v)_{v\in V_f}$ be a coherent collection of parahoric subgroups. Let $U_v$ be the projection of $U$ in $G(k_v)$. Then as $U\cap \prod P_v$ is a compact-open subgroup of $G(A_f)$, its projection in $G(k_v)$ is a hyperspecial parahoric 
subgroup of $G(k_v)$ for all but finitely many $v\in V_f$ ([Ti2], 3.9). If for a $v\in V_f$, the projection of $U\cap \prod P_v$ in $G(k_v)$ (this projection is contained in $U_v\cap P_v$) is a hyperspecial parahoric subgroup, then by maximality of these subgroups among compact subgroups of $G(k_v)$, we conclude that $P_v$ is hyperspecial and $U_v =P_v$.   Thus for all but finitely many $v\in V_f$, $P_v$ is hyperspecial and $U_v=P_v$. Now if $(P'_v)_{v\in V_f}$ is another coherent collection of parahoric subgroups, then $U :=\prod P'_v$ is a compact-open subgroup of $G(A_f)$ and we conclude from the above observations that for all but finitely many $v\in V_f$, $P'_v = P_v$. 
\vskip1mm

We fix a coherent collection $(P_v)_{v\in V_f}$ of parahoric subgroups $P_v$ of $G(k_v)$ and let 
$\Lambda: = G(k)\cap\prod_{v\in V_f}P_v$. Let $\Gamma$ be the normalizer of  $\Lambda$  in $G(k_{v_o})$. Note that as the normalizer of $\Lambda$ in $G(k)$ equals $\Lambda$, $\Gamma \cap G(k) = \Lambda$. {\it We assume in the sequel that}  $\chi(\Gamma)\leqslant 1$. 
\vskip1mm

The Haar-measure $\mu$ on $G(k_{v_o})$ is the one used in [BP].
\vskip1mm

All unexplained notations are as in [BP] and [P]. Thus for a number field $K$, $D_K$ denotes the absolute value of its 
discriminant, $h_K$ its class number, i.\,e., the order of its class group 
$Cl (K)$. We shall denote by $n_{K, 3}$ the order of the $3$-primary component of $Cl(K)$, and by $h_{K, 3}$ the order of the subgroup (of $Cl(K)$) consisting 
of the elements of order dividing $3$. Then $h_{K,3}\leqslant n_{K,3}
\leqslant h_K$. 
 \vskip1mm
 
For a number field $K$, $U(K)$ will denote the multiplicative-group of units of $K$, and
$K_3$ the subgroup of $K^{\times}$ consisting of
the elements $x$ such that for every normalized valuation $v$ of $K$,
$v(x)\in 3\bZ$. 
 \vskip1mm
 
We will denote the degree $[k:\,\bQ ]$ of $k$ by $d$, and for any 
nonarchimedean place $v$ of $k$, $q_v$ will denote the cardinality of 
the residue field ${\mathfrak f}_v$ of $k_v$.
\vskip1mm

For a positive integer $n$, $\mu_n$ will denote the kernel of the endomorphism 
$x\mapsto x^n$ of ${\rm GL}_1$. Then the center $C$ of $G$ is
$k$-isomorphic to the kernel of the norm map $N_{\ell/k}$ from the algebraic group
$R_{\ell/k}(\mu_3)$, obtained from $\mu_3$ by Weil's restriction of
scalars, to $\mu_3$. Since the norm map $N_{\ell/k}: \
\mu_3(\ell)\rightarrow \mu_3(k)$ is onto, $\mu_3
(k)/N_{\ell/k}(\mu_3(\ell))$ is trivial, and hence, the Galois
cohomology group $H^1(k,C)$ is isomorphic to the kernel of the homomorphism
$\ell^{\times}/{\ell^{\times}}^3\rightarrow
k^{\times}/{k^{\times}}^3$ induced by the norm map. 
This kernel equals $\ell^{\bullet}/{\ell^{\times}}^3$, where
$\ell^{\bullet} = \{ x\in \ell^{\times}\ |\ N_{\ell/k}(x)\in {k^{\times}}^3\}$.
\vskip2mm

\ni{\bf 2.2.} For $v\in V_f$, let the ``type " $\Theta_v$  of $P_v$ be as in 2.2 of [BP], and $\Xi_{\Theta_v}$ be as in 2.8 there.   
We observe here, for later use, that for a nonarchimedean place $v,$ 
$\Xi_{\Theta_v}$ is nontrivial if, and only if, $G$ splits at $v$ 
(then $v$ splits in $\ell$, i.\,e., $k_v\otimes_k \ell$ is a direct product of two fields, each isomorphic to $k_v$ ) and $P_v$ is an Iwahori subgroup
of $G(k_v)$ (then $\Theta_v$ is the empty set), 
and in this case $\#\Xi_{\Theta_v}=3.$  
\vskip1mm

We recall that hyperspecial parahoric subgroups of $G(k_v)$ are conjugate to each other under $\overline{G}(k_v)$, see [Ti2, 2.5], and $G(k_v)$ contains a hyperspecial parahoric subgroup if 
and only if $v$ is unramified in $\ell$ and $G$ is quasi-split at $v$
(i.\,e., it contains a Borel subgroup defined over $k_v$). 
Let  $\cT$ be the set of nonarchimedean places $v$ of $k$ such that in the collection $(P_v)_{v\in V_f}$ under consideration, $P_v$ is not maximal, and also all those $v$ which are unramified in $\ell$ and $P_v$ is not hyperspecial. Let $\cT_0$ be the subset of $\cT$ consisting of places where the group $G$ is anisotropic. 
Then $\cT$ is finite, and for any nonarchimedean $v\not\in\cT,$
$\Xi_{\Theta_v}$ is trivial. We note that every place $v\in \cT_0$ splits in $\ell$ since an absolutely simple anisotropic group over a 
nonarchimedean local field is necessarily of {\it inner} type $A_n$ 
(another way to see this is to recall that, over a local field, the only 
central simple algebras which admit an involution of the second kind are the matrix algebras).  We also note that every absolutely simple group of type $A_2$ defined and isotropic over a field $K$ is quasi-split (i.\,e., it contains a 
Borel subgroup defined over $K$). 
\vskip1mm

If $v$ does not split in $\ell$ (i.\,e., $\ell_v :=k_v\otimes_k \ell$ is a field), then $G$ is quasi-split over $k_v$ (and its $k_v$-rank is $1$). In this case, if $P_v$ is not an Iwahori subgroup, then it is a maximal parahoric subgroup of $G(k_v)$, and there are two conjugacy classes of maximal parahoric subgroups in $G(k_v)$. Moreover, if $P'$ and $P''$ are the two maximal parahoric 
subgroups of $G(k_v)$ containing a common Iwahori subgroup $I$, then the derived subgroups of any Levi subgroups of the reduction mod $\mathfrak p$ of $P'$ and $P''$ are nonisomorphic: if $\ell_v$ is an unramified extension of $k_v$, then the two derived subgroups are ${\rm SU}_3$ and ${\rm SL}_2$, and if $\ell_v$ is a ramified extension of $k_v$, then the two derived subgroups are ${\rm SL}_2$ and ${\rm PSL}_2$, see [Ti2], 3.5. Hence, $P'$ is not conjugate to $P''$ under the action of $({\rm Aut}\,G)(k_v)$\,($\supset{\overline G}(k_v)$). In particular, {\it if an element of ${\overline G}(k_v)$ normalizes $I$, then it normalizes both $P'$ and $P''$ also.} If $v$ ramifies in $\ell$, then $P'$ and $P''$ are 
of same 
volume with respect to any Haar-measure on $G(k_v)$, since, in this case, 
$[P':I] 
= [P'':I]$.    
\vskip2mm

\ni{\bf 2.3.} By Dirichlet's unit theorem, $U(k)\cong \{\pm 1\}\times 
{\bZ}^{d-1}$, and $U({\ell})\cong \mu(\ell)\times {\bZ}^{d-1}$, where
$\mu(\ell)$ is the finite cyclic group of roots of unity in $\ell$. Hence, 
$U(k)/U(k)^3\cong (\bZ/3\bZ)^{d-1}$, and $U({\ell})/U({\ell})^3\cong
\mu(\ell)_3\times (\bZ/3\bZ)^{d-1}$, where $\mu(\ell)_3$ is the group
of cube roots of unity in $\ell$. Now we observe that
$N_{\ell/k}(U({\ell}))\supset N_{\ell/k}(U(k)) = U(k)^2$, which implies
that the homomorphism  $U({\ell})/U({\ell})^3\rightarrow U(k)/U(k)^3$, 
induced by the norm map, is onto. The kernel of this homomorphism is
clearly $U(\ell)^{\bullet}/U(\ell)^3$, where $U(\ell)^{\bullet} = U(\ell)\cap
\ell^{\bullet}$, and its order equals $\#\mu(\ell)_3$. 

\vskip1mm
 
The short exact sequence $(4)$ in the proof of Proposition 0.12 in [BP]
gives us the following exact sequence: $$1\rightarrow
U({\ell})^{\bullet}/U({\ell})^3\rightarrow
\ell^{\bullet}_3/{\ell^{\times}}^3\rightarrow (\cP\cap \cI^3)/\cP^3,$$
where $\ell^{\bullet}_3 = \ell_3\cap\ell^{\bullet}$, $\cP$ is the group of all
fractional 
principal ideals of $\ell$, and $\cI$ the group of all fractional ideals (we use multiplicative
notation for the group operation in both $\cI$ and $\cP$). Since the
order of the last group of the above exact sequence is $h_{\ell,3}$,
see $(5)$ in the proof of Proposition 0.12 in [BP], we conclude that 
$$\#\ell^{\bullet}_3/{\ell^{\times}}^3 \leqslant \#\mu(\ell)_3\cdot h_{\ell,3}.$$

Now we note that the order of the first term of the short exact sequence 
of Proposition 2.9 of [BP], for $G' =G$ and $S =V_{\infty}$, 
is $3/\#\mu(\ell)_3$.

The above observations, together with Proposition 2.9 and Lemma 5.4 of [BP], and a
close 
look at the 
arguments in 5.3 and 5.5 of [BP] for $S=V_{\infty}$ and $G$ of type $^2A_2$, 
give us the following upper bound (note that for our $G$, in 5.3 of [BP], 
$n =3$): 
\vskip2mm

\ni $(0)$   { \ \ } { \ \  } {\ \ }  {\ \ }  {\ \ } \ \ \ \  \ \ \ \ \ \ \ \ \ \ \ \ \ \ $[\Gamma : \Lambda] \leqslant 
3^{1+\#\cT_0} h_{\ell, 3} \prod_{v\in \cT-\cT_0}
\# \Xi_{\Theta_v}.$

\vskip2mm

We note also that Proposition 2.9 of [BP] applied to $G' = G$ and
$\Gamma' = \Gamma$, implies that the index $[\Gamma :\Lambda]$ of $\Lambda$ in 
$\Gamma$ is a power of $3$.

\vskip1mm 
  
\ni{\bf 2.4.} As we mentioned in 1.3, $\chi(\Gamma) = 3 \mu(G(k_{v_o})/\Gamma)$. Our aim here is to find a lower bound for 
$\mu(G(k_{v_o})/\Gamma).$  
For this
purpose, we first note that 
$$\mu(G(k_{v_o})/\Gamma)=\frac{\mu(G(k_{v_o})/\Lambda)}{[\Gamma:\Lambda]}.$$
As the Tamagawa number $\tau_k(G)$ of $G$ equals $1,$ the volume formula
of [P] (recalled in \S3.7 of [BP]), for $S =V_{\infty}$, gives us
$$\mu(G(k_{v_o})/\Lambda)= D_k^4(D_\ell/ D_k^2)^{5/2}(16\pi^5)^{-d}\cE = 
(D_{\ell}^{5/2}/D_k )(16\pi^5 )^{-d}\cE ;$$
where
$\cE=\prod_{v\in V_f}e(P_v),$ and 
$$e(P_v)=\frac{q_v^{(\dim\oM_v+\dim\ocM_v)/2}}{\#\oM_v({\mathfrak f}_v)}.$$ We 
observe that if $P_v$ is hyperspecial,  
$$e(P_v)=\(1-\frac1{q_v^2}\)^{-1}\(1-\frac1{q_v^3}\)^{-1}\ \ {\rm or}\ \ 
\(1-\frac1{q_v^2}\)^{-1}\(1+\frac1{q_v^3}\)^{-1}$$
according as $v$ does or does not split in $\ell.$ If $v$ ramifies in
$\ell$ and $P_v$ is a maximal parahoric subgroup of $G(k_v)$, then
$$e(P_v) = \(1-\frac{1}{q_v^2}\)^{-1}.$$
Now let $\zeta_k$ be the Dedekind zeta-function of $k$, and $L_{\ell|k}$ be the Hecke $L$-function associated to the quadratic Dirichlet character of $\ell/k$. Then as 
$$\zeta_k(2) =\prod_{v\in V_f}\(1-\frac1{q_v^2}\)^{-1},$$ and 
$$L_{\ell|k}(3) ={\prod} '\(1-\frac1{q_v^3}\)^{-1}{\prod} ''\(1+\frac1{q_v^3}\)^{-1},$$
where $\prod'$ is the product over those nonarchimedean places of $k$ which 
split in $\ell$, and $\prod''$
is the product over all the other nonarchimedean places $v$ which do not 
ramify in $\ell$, 
we see that
$$\cE=\zeta_k(2)L_{\ell |k}(3)\prod_{v\in\cT}e'(P_v);$$
where, for $v\in\cT,$ 
\vskip1mm

\ni $\bullet $ if $v$ splits in $\ell$, $e'(P_v)
=e(P_v)(1-\frac1{q_v^2})(1-\frac1{q_v^3}),$
\vskip.5mm

\ni $\bullet $ if $v$ does not split in $\ell$ but is unramified in $\ell$, 
$e'(P_v) = e(P_v)(1-\frac1{q_v^2})(1+\frac1{q_v^3}),$
\vskip.5mm

\ni $\bullet $ if $v$ ramifies in $\ell$, $e'(P_v) = e(P_v)(1-\frac{1}{q_v^2}).$
\vskip1mm   

Thus

\begin{equation}
\mu(G(k_{v_0})/{\Gamma})=\frac{D_{\ell}^{5/2}\zeta_k(2)
L_{\ell |k}(3)}{(16\pi^5)^d[\Gamma:\Lambda]D_k}\prod_{v\in\cT}e'(P_v)
\geqslant \frac{D_{\ell}^{5/2}\zeta_k(2)L_{\ell |k}(3)}{3(16\pi^5)^d 
h_{\ell,3}D_k}\prod_{v\in\cT}e''(P_v),
\end{equation}
where, for $v\in \cT-\cT_0$, $e''(P_v) =e'(P_v)/{\#\Xi_{\Theta_v}}$, and for $v\in\cT_0$, $e''(P_v) = e'(P_v)/3$.
\vskip3mm

\ni{\bf 2.5.} Now we provide the following list of values of $e^{\prime}(P_v )$ and 
$e^{\prime\prime}(P_v)$, for all $v\in\cT.$\\

\ni ({\it i}) $v$ {\it splits in $\ell$ and $G$ splits at $v$:}\\

(a) if $P_v$ {\it is an Iwahori subgroup}, then
$$e^{\prime\prime}(P_v)= e'(P_v)/3,$$ and $$e^{\prime} (P_v ) 
=(q_v^2+q_v+1)(q_v+1);$$

(b) if $P_v$ {\it is not an Iwahori subgroup} (note that as $v\in \cT$, $P_v$ 
is not hyperspecial), then
$$e^{\prime\prime}(P_v)=e'(P_v)= q_v^2+q_v+1;$$
\vskip1mm

\ni({\it ii}) {\it $v$ splits in $\ell$ and $G$ is anisotropic at $v$} (i.\,e., $v\in\cT_0$): $$e''(P_v)= e'(P_v)/3,$$ and $${e^{\prime}(P_v)} = (q_v-1)^2(q_v +1);$$
\vskip1mm
 
\ni({\it iii}) $v$ {\it does not split in $\ell$, and $\ell_v =k_v\otimes_k\ell$ is an unramified extension of $k_v$, then}
$$e^{\prime\prime}(P_v)=e^{\prime}(P_v)=\left\{\begin{array}{cc}
q_v^3+1&\mbox{ if $P_v$ is an Iwahori subgroup}\\
q_v^2-q_v+1&\mbox{ if $P_v$ is a non-hyperspecial maximal parahoric subgroup;}\end{array}
\right.
$$
\ni({\it iv}) $v$ {\it does not split in $\ell$, and $\ell_v = k_v\otimes_k\ell$ is a ramified extension of $k_v$, then}
$$e''(P_v) = e'(P_v) = q_v+1.$$
\vskip2mm

\ni{\bf 2.6.}  As $\chi(\Gamma) \leqslant 1$,   
$\mu(G(k_{v_o})/\Gamma )\leqslant 1/3$. So from $(1)$ in 2.4 we get 
the following:
\begin{equation}
{1/3}\geqslant \mu(G(k_{v_0})/\Gamma ) \geqslant \frac{D_{\ell}^{5/2}\zeta_k(2)L_{\ell|k}(3)}{3(16\pi^{5})^{d}h_{\ell, 3}D_k}
\prod_{v\in\cT}e^{\prime\prime}(P_v).
\end{equation}
We know from the Brauer-Siegel theorem that for all real $s>1$, 
\begin{equation}
h_\ell R_\ell\leqslant w_\ell s(s-1)\Gamma(s)^d((2\pi)^{-2d}D_\ell)^{s/2}\zeta_\ell(s),
\end{equation}
where $h_\ell$ is the class number and $R_\ell$ is the regulator of $\ell$, and
$w_\ell$ is the order of the finite group of roots of unity contained in $\ell.$ 
Zimmert [Z] obtained the following lower bound for the regulator $R_{\ell}$.
$$R_{\ell}\geqslant 0.02w_{\ell}e^{0.1d}.$$ 
Also, we have the following lower bound for the regulator obtained by 
Slavutskii [Sl] 
using a variant of the argument of Zimmert [Z]:
$$R_\ell\geqslant 0.00136 w_\ell\,e^{0.57d}.$$
We deduce from this bound and $(3)$ that 

\begin{equation}
\frac{1}{h_{\ell,3}}\geqslant \frac1{h_\ell}\geqslant\frac{0.00136}{s(s-1)}\( \frac{(2\pi)^s e^{0.57}}{\Gamma(s)}\)^d
\frac1{D_\ell^{s/2}\zeta_\ell(s)};
\end{equation}
if we use Zimmert's lower bound for $R_{\ell}$ instead, we obtain 

\begin{equation}
\frac{1}{h_{\ell,3}}\geqslant \frac1{h_\ell}\geqslant\frac{0.02}{s(s-1)}\( \frac{(2\pi)^s e^{0.1}}{\Gamma(s)}\)^d
\frac1{D_\ell^{s/2}\zeta_\ell(s)}.
\end{equation}

\vskip2mm

\ni {\bf 2.7. Lemma.} {\it For every integer $r\geqslant 2$, $\zeta_k(r)^{1/2}L_{\ell|k}(r+1)> 1.$}
\vskip2mm

\ni{\it Proof}. Recall that 

$$\zeta_k(r) = \prod_{v\in V_f}\(1-\frac1{q_v^r}\)^{-1},$$ and 
$$L_{\ell|k}(r+1) = {\prod} '\(1-\frac1{q_v^{r+1}}\)^{-1}{\prod} ''
\(1+\frac1{q_v^{r+1}}\)^{-1},$$ where $\prod'$ is the product over all finite places $v$ of $k$ which split 
over $\ell$ and $\prod''$
is the product over all the other nonarchimedean $v$ which do not ramify in $\ell$. Now the lemma follows from the following simple observation. 
\vskip1mm

For any 
positive integer $q\geqslant 2,$ 
\begin{eqnarray*}
\(1-\frac1{q^r}\)\(1+\frac1{q^{r+1}}\)^2&=&1-\frac{q-2}{q^{r+1}}-\frac{2q-1}{q^{2r+2}}-\frac{1}{q^{3r+2}}<1\ .\\
\end{eqnarray*}
\vskip2mm

\ni{\bf 2.8. Corollary.} {\it For every integer $r\geqslant 2$,} $$\zeta_k(r)L_{\ell|k}(r+1)> \zeta_k(r)^{1/2}>1.$$

\vskip2mm

\ni {\bf 2.9. Remark.} The following bounds are obvious from the Euler-product expression for the zeta-functions. For every integer $r\geqslant 2$, $$\zeta(dr)\leqslant\zeta_k(r)\leqslant \zeta(r)^d,$$ where $\zeta(j)=\zeta_{\bQ}(j)$. Now from the above corollary we deduce that 
\begin{equation} 
\zeta_k(2)L_{\ell|k}(3)> \zeta_k(2)^{1/2}\geqslant \zeta(2d)^{1/2}>1.
\end{equation} 
\vskip2mm

\ni {\bf 2.10.} Since $e'' (P_v )\geqslant 1$ for all $v\in\cT$, see 2.5 above, and  $D_{\ell}/D_k^2$ is an integer, so in particular, $D_k\leqslant D_{\ell}^{1/2}$, see, for example, Theorem A in the appendix of [P], bounds $(2)$ and $(3)$ lead to the following 
bounds by taking $s=1+\delta$, with $0<\delta\leqslant 2$, in $(3)$

\begin{equation}
D_k^{1/d}\leqslant D_\ell^{1/2d}<{\varphi_1}(d,R_{\ell}/w_{\ell},\delta)\ \ \ \ \ \ \ \ \ \ \ \ \ \ \ \ 
\end{equation}

$$:=\Big(\frac{\delta(1+\delta)}
{\zeta(2d)^{1/2}(R_{\ell}/w_{\ell})}\Big)^
{1/{(3-\delta)d}}\big( 2^{3-\delta}\pi^{4-\delta}\Gamma(1+\delta)\zeta(1+\delta)^2\big)^{{1}/{(3-\delta)}},$$

\begin{equation}D_k^{1/d}\leqslant D_{\ell}^{1/2d}<{\varphi_2}(d,h_{\ell,3})
:=\big[ \frac{2^{4d}\pi^{5d}h_{\ell,3}}{\zeta(2d)^{1/2}}\big]^{1/4d},
\end{equation}
and
\begin{equation}
D_{\ell}/D_k^2< {\mathfrak p}(d,D_k,h_{\ell,3}):= \big[
  \frac{2^{4d}\pi^{5d}h_{\ell,3}}{\zeta(2d)^{1/2}D_k^4}\big]^{2/5}.
\end{equation}

\vskip2mm

Using the bound $R_{\ell}/w_{\ell}\geqslant 0.00136 e^{0.57d}$ due to Slavutskii, we obtain the following bound from $(7)$:
\begin{eqnarray}
&&D_k^{1/d}
\leqslant D_\ell^{1/2d}<f(\delta,d)\\
&:=&\big[\frac{\delta(1+\delta )}{0.00136}\big]^{1/{(3-\delta)d}}\cdot\big[2^{3-\delta}\pi^{4-\delta}\Gamma(1+\delta)\zeta(1+\delta)^2e^{-0.57}\big]^{1/{(3-\delta)}}.\nonumber
\end{eqnarray}
\vskip1mm

\ni{\bf 2.11.} As $\chi(\Lambda)= 3\mu(G(k_{v_o})/{\Lambda})$, $$\chi(\Gamma)= \frac{\chi(\Lambda)}{[\Gamma :\Lambda ]}= \frac{3\mu (G(k_{v_o})/{\Lambda})}
{[\Gamma : \Lambda]}.$$ 
Now since $[\Gamma : \Lambda]$ is a power of $3$ (see 2.3), if 
$\chi(\Gamma)$ is a reciprocal integer, the numerator of the rational number 
$\mu(G(k_{v_o})/\Lambda)$ is a power of $3$. 

We recall from 2.4 that 
$$\mu(G(k_{v_o})/\Lambda) =
(D_{\ell}^{5/2}/D_k)(16\pi^5)^{-d}\zeta_k(2)L_{\ell|k}(3)\prod_{v\in\cT}e'(P_v).$$
Using the functional equations
$$\zeta_k(2)= (-2)^d\pi^{2d}D_k^{-3/2}\zeta_k(-1),$$
and
$$L_{\ell|k}(3)=(-2)^d\pi^{3d} (D_{k}/D_{\ell})^{5/2}L_{\ell|k}(-2),$$
we can rewrite the above as:
\begin{equation}
\mu(G(k_{v_o})/\Lambda)=
2^{-2d}\zeta_k(-1)L_{\ell|k}(-2)\prod_{v\in\cT}e'(P_v).
\end{equation}
Hence we obtain the following proposition.  
\vskip2mm

\ni{\bf 2.12. Proposition.} {\it If the orbifold Euler-Poincar\'e 
characteristic $\chi (\Gamma)$ of $\Gamma$ is a reciprocal integer, then the  
numerator of the rational number 
$2^{-2d}\zeta_k(-1)L_{\ell|k}(-2)\prod_{v\in\cT}e'(P_v)$ is a power of
$3$.  Moreover, as $e'(P_v)$ is an integer for all $v$, the numerator of $\mu :=
2^{-2d}\zeta_k(-1)L_{\ell|k}(-2)$ is also a power of $3$.}

\bs
\ni
\begin{center}
{\bf 3. Determining ${\ell}$ when $k =\bQ$}  
\end{center}
\vskip4mm

\ni {\it We will assume in this, and the next section, that
  $k=\bQ$. Then $\ell =\bQ(\sqrt{-a})$, where $a$ is a square-free
  positive integer.}

\vskip2mm

 We will now find an upper bound for 
$D_\ell.$ 
\vskip1mm

\ni {\bf 3.1.} Since $D_k = D_{\bQ} = 1$,  and $e''(P_v)\geqslant 1$, from $(2)$, $(5)$ and $(6)$, taking $s=1+\delta$, we get the following: 

\begin{equation}
D_\ell< (2\pi)^2\(\frac{5^2\cdot \delta(1+\delta )\cdot\Gamma(1+\delta)\zeta(1+\delta)^2}{e^{0.1}\zeta(2)^{1/2}}\)^{2/(4-\delta)}.
\end{equation}

\ni Letting $\delta=0.34,$ we find that $D_\ell < 461.6.$  Hence we
conclude that $D_\ell\leqslant 461.$
\vskip1mm

\ni Thus we have established the following.
\vskip2mm

\ni {\bf 3.2.} {\it If $\chi (\Gamma)\leqslant 1$ and $k = \bQ$, then $D_\ell\leqslant 461.$}
\vskip3mm

\ni{\bf 3.3.} We will now improve the upper bound for the discriminant of $\ell$ using the table of 
class numbers of imaginary quadratic number fields.  
\vskip1mm

Inspecting the table of class numbers of $\ell = \bQ(\sqrt{-a})$, with 
$D_\ell\leqslant 461$, in [BS], we find that 
$h_\ell\leqslant 21$, and hence, $h_{\ell,3}\leqslant n_{\ell,3}\leqslant 9.$
\vskip1mm

Since $D_{\bQ}=1,$ $\zeta_{\bQ}(2) = \zeta(2)= \pi^2/6$ and $\zeta(3)L_{\ell |\bQ}(3) =\zeta_{\ell}(3)>1$, $(2)$ provides us the following
bounds
\begin{eqnarray*}
1&\geqslant&\frac{D_{\ell}^{5/2}L_{\ell|{\bQ}}(3)}{2^5\cdot 3\cdot \pi^3\cdot 
h_{\ell,3}}\prod_{v\in\cT}e^{\prime\prime}(P_v)\\
&\geqslant&\frac{D_{\ell}^{5/2}\zeta_\ell(3)}{2^5 \cdot 3\cdot\pi^3\cdot h_{\ell,3}\zeta(3)}\\
&>&\frac{D_{\ell}^{5/2}}{2^5\cdot 3\cdot \pi^3\cdot h_{\ell,3}\zeta(3)}.
\end{eqnarray*}
Hence, in particular,  as $h_{\ell,3}\leqslant n_{\ell,3}$, 

\begin{eqnarray*}
D_\ell&< &\big(2^5\cdot 3 \cdot \pi^3\cdot n_{\ell,3}\zeta(3)\big)^{2/5}.
\end{eqnarray*}

\ni  The above leads to the following bounds once
the value of $n_{\ell,3}$ is determined.

$$\begin{array}{cccc}
n_{\ell,3}&1&3&9\\
D_\ell\leqslant&26&40&63
\end{array}$$

The last column of the above table implies that we need only 
consider $D_\ell\leqslant 63.$

\vskip3mm

\ni{\bf 3.4.} We will further limit the possibilities for $D_\ell.$  
If $40<D_\ell\leqslant63,$ we observe that $n_{\ell,3}\leqslant 3$ from the table in Appendix.  Hence, from the middle column of the above table we infer that 
$D_\ell$ can at most be $40$. 

For $26<D_\ell\leqslant40,$  we see from the table in
Appendix  that unless $D_\ell=31$, $n_{\ell,3}=1$, and
the first column of the above table shows that if $n_{\ell, 3} = 1$, $D_\ell\leqslant
26.$  Hence, the only possible values of $D_{\ell}$ are $31$ or
$D_{\ell}\leqslant 26$.

    From the table in Appendix we now see that the possible values of $h_{\ell,3}$ and $D_{\ell}$ are the following (note that if $n_{\ell,3} = 3$, then $h_{\ell,3} =3$ also).
\vskip1mm

\ni $h_{\ell,3}=3:$\  $D_\ell=23,31.$\\

\ni $h_{\ell,3}=1:$\  $D_\ell=3,4,7,8,11,15,19,20,24.$
\vskip2mm

Now we recall that for $\ell =\bQ({\sqrt{-a}})$ , $D_{\ell} = a$ if $a\equiv 3$ (mod $4$), and $D_{\ell} = 4a$ otherwise. Using this we can paraphrase the above result as follows.

\vskip2mm

\ni {\bf 3.5. Proposition.} {\it Let $k =\bQ$. Then $\ell=\bQ(\sqrt{-a})$, where $a$ is one of the following eleven integers,
$$1,\ 2,\ 3,\ 5,\ 6,\ 7,\ 11,\ 15,\ 19,\ 23,\ 31.$$}

{\it The following table provides the value of} $$\mu : =
\frac{D_{\ell}^{5/2}\zeta(2)
L_{\ell |\bQ}(3)}{16\pi^5} = -\frac{1}{48}L_{\ell |\bQ}(-2)$$
({\it recall the functional equation} 
$L_{\ell |\bQ}(3) = -2\pi^3 D_{\ell}^{-5/2}L_{\ell | \bQ}(-2)$). 
\vskip1mm

$$\begin{array}{cccccccccccc}
a&1&2&3&5&6&7\\
L_{\bQ(\sqrt{-a})|\bQ}(-2)&-1/2&-3&-2/9&-30&-46&-16/7\\
\mu&{1/96}&{1/16}&{1/216}&{5/8}&{23/24}&{1/21}\\
\end{array}$$
\vskip1mm

$$\begin{array}{cccccccccc}
a&11&15&19&23&31\\
L_{\bQ(\sqrt{-a})|\bQ}(-2)&-6&-16&-22&-48&-96\\
\mu&{1/8}&{1/3}&{11/24}&{1}&{2.}\\

\end{array}$$
\vskip5mm

\vskip2mm

\ni{\bf 3.6.} The volume formula of [P] and the results of [BP] apply
equally well to {\it noncocompact arithmetic} subgroups. So if we wish
to make a list of all noncocompact arithmetic subgroups $\Gamma$ of ${\rm SU}(2,1)$ 
whose orbifold Euler-Poincar\'e characteristic $\chi (\Gamma)$ is $\leqslant 1$, we can proceed as above. 
If $\Gamma$ is such a subgroup, then, associated to it, there is an absolutely simple simply connected algebraic group $G$ defined and (by Godement compactness criterion) isotropic over a number field $k$ such that $G(k\otimes_{\bQ}{\bR})$ is isomorphic to the direct product of 
${\rm SU}(2,1)$ with a compact semi-simple Lie group. But since $G$ is 
$k$-isotropic, for every place $v$ of $k$, $G$ is isotropic over $k_v$, and hence, $G(k_v)$ is noncompact. In particular, for 
every archimedean place $v$ of $k$, $G(k_v)$ is noncompact. This implies that $k =\bQ$, $G$ is an absolutely simple 
simply connected $\bQ$-group of type $A_2$ of $\bQ$-rank 1 (and hence $G$ is 
quasi-split over $\bQ$). Moreover, $G$ splits over an imaginary quadratic 
extension $\ell = \bQ({\sqrt{-a}})$ of $\bQ$. For a given positive integer 
$a$, there is a unique such $G$ (up to $\bQ$-isomorphism). The considerations of 
3.1--3.4 apply again and imply that $a$ has to be one of the eleven integers listed in Proposition 3.5. 
\vskip1mm

We fix a coherent collection $(P_p)$ of maximal parahoric subgroups $P_p$ of $G(\bQ_p)$ such that $P_p$ is hyperspecial whenever $G(\bQ_q)$ contains such a parahoric subgroup. Let $\Lambda = G(\bQ)\cap\prod_p P_p$. (This $\Lambda$ is a ``Picard modular group''.) From the volume formula of [P], recalled in 2.4, we obtain that $$\chi (\Lambda) = 3\mu (G(\bR)/{\Lambda}) = 3\frac{D_{\ell}^{5/2}\zeta_{\bQ}(2)L_{\ell|{\bQ}}(3)}{16\pi^5} = \frac{D_{\ell}^{5/2}L_{\ell|{\bQ}}(3)}{32\pi^3} =-\frac{1}{16}L_{\ell|{\bQ}}(-2)=3\mu,$$
where we have used the functional equation for the $L$-function $L_{\ell|{\bQ}}$ recalled in 3.5, and the fact that $\zeta_{\bQ}(2)=\zeta(2) = {\pi^2}/6$. (We note 
that the above computation of the orbifold Euler-Poncar\'e
characteristic of Picard modular groups is independently due to
Rolf-Peter Holzapfel, see [Ho], section 5A.) Now we can use the table
of values of $L_{\ell|\bQ}(-2)$ given in 3.5 to compute the precise value of $\chi(\Lambda)$ for each $a$.      
\vskip1mm

Among all arithmetic subgroups of $G$ contained in $G(\bQ)$, the above $\Lambda$ has the smallest orbifold Euler-Poincar\'e characteristic. Its normalizer $\Gamma$ in 
$G(\bR)$ has the smallest orbifold Euler-Poincar\'e characteristic among all discrete subgroups commensurable with $\Lambda$. Note that $\Lambda$ has torsion. 

\newpage
\ni
\begin{center}
{\bf 4. Determination of $G$ and the parahoric subgroups $P_v$}
\end{center}
\vskip3mm
 
We continue to assume in this section that $k=\bQ$. We will use the usual identification of a nonarchimedean place $v$ of $\bQ$ with the characteristic $p$ of the residue field of $\bQ_v$. Let $\ell$ be one of the eleven imaginary quadratic extensions of $\bQ$ listed in Proposition 3.5. $\cR_{\ell}$ will denote the set of rational primes which ramify in $\ell$.   
\vskip2mm

\ni{\bf 4.1.}  Let $\cD$, the involution $\sigma$, the hermitian form $h$, and the $k$-group $G$, for $k=\bQ$, be as in 1.2.  As in 2.1 we fix a coherent collection $(P_p)$ of parahoric subgroups of $G(\bQ_p)$. Let $\Lambda = G(\bQ)\cap\prod_{p} P_p$, and $\Gamma$ be the normalizer of $\Lambda$ in $G({\bf R})$. {\it We assume that $\Gamma$ is cocompact and $\chi(\Gamma)$ is a reciprocal integer.}

\vskip1mm
We first show that $\cD$ is a cubic division algebra. Assume, if possible, that $\cD = \ell$. Then $h$ is a hermitian form on $\ell^3$. As the arithmetic subgroup $\Gamma$ of $G(\bR)$ is cocompact, by 
Godement compactness criterion, $h$ is an anisotropic form on $\ell^3$. On the other hand, its signature over 
$\bR$ is $(2,1)$. The hermitian form $h$ gives us a quadratic form $q$ on the six dimensional 
$\bQ$-vector space $V= \ell^3$ defined as follows: $$ q(v) = h (v, v) \ \ {\rm for} \ \ v\in V .$$
The quadratic form $q$ is isotropic over $\bR$, and hence by Meyer's theorem it is isotropic over $\bQ$ (cf. [Se2]).  
This implies that $h$ is isotropic and we have arrived at a contradiction.

\vskip2mm

\ni{\bf 4.2.} Let $\cT$ be the finite set of rational primes $p$ such that $P_p$ is not maximal and also those 
$p\notin \cR_{\ell}$ such that $P_p$ is not 
hyperspecial, and $\cT_0$ be the subset of $\cT$ consisting of $p$ such that $G$ 
is anisotropic over $\bQ_p$. Since $\cD$ must ramify at at least some 
nonarchimedean places of $\ell$, $\cT_0$ is nonempty. As pointed out in 2.2, 
{\it every $p\in \cT_0$ splits in} $\ell$. Theorem 4.4 lists all possible 
$\ell$, $\cT$, $\cT_0$, and the parahoric subgroups $P_p$.
\vskip1mm

As $\zeta_{\bQ}(2) = \zeta (2) =\pi^2 /6$, using the functional
equation
$$L_{\ell|\bQ}(3)= -2\pi^3 D_{\ell}^{-5/2}L_{\ell|{\bQ}}(-2),$$
we obtain the following from bound $(1)$ for $k =\bQ$ : 
$$\chi(\Gamma)=3\mu (G(\bR )/\Gamma )\geqslant{\frac{\mu}{h_{\ell, 3}}}
\prod_{p\in \cT} e^{\prime\prime}(P_p ),$$
where $\mu$ is as in 3.5.
\vskip2mm

\ni{\bf 4.3.} We recall here that given a square-free integer $a$, an
odd 
prime $p$ splits in $\ell 
= {\bQ (\sqrt{-a})}$ if, and only if, $p$ does not divide $a$, and $-a$ is a 
square modulo $p$; $2$ splits in $\ell$ if, and only if,  
$-a \equiv  1 \ ({\rm mod}\ 8)$; see [BS], \S 8 of Chapter 3. 
A prime $p$ ramifies in $\ell$ if, and only if, $p|D_{\ell}$; see [BS], \S 7 of Chapter 2 and \S 8 of Chapter 3.            
\vskip1mm

Now using Proposition 3.5, the fact that the 
numerators of $\mu$ and $\mu (G(\bR )/\Lambda ) = \mu\prod_{p\in\cT}e'(P_p)$ are powers of 
$3$ (Proposition 2.12), the value of $\mu$ given in 3.5, the values of $e^{\prime}(P_p)$, $e^{\prime\prime}(P_p )$ given in 2.5, the value of $h_{\ell, 3}$ given in 3.4, and the fact that $\chi(\Gamma)\leqslant 1$, we see by a direct computation 
that the following holds.
\vskip2mm

\ni {\bf 4.4. Theorem.} {\it $\cT_0$ consists of a single prime, and the pair $(a,p)$, where $\ell = 
\bQ(\sqrt{-a} )$, and $\cT_0 = \{ p \}$, belongs to the set $\{ (1,5), (2,3), (7,2), (15,2), (23,2)\}$. Moreover, $\cT = \cT_0$ unless $a = 1, 2$ or $7$. For $a = 1, 2$, the possibilities are $\cT= \cT_0$ and $\cT = \cT_0\cup\{2\}$. For $a = 7$ the possibilities for $\cT$ are $\cT =\cT_0 =\{2\}$, $\cT =\{2,3\}$, and $\cT=\{2,5\}$. 
}
\vskip2mm

\ni {\bf 4.5.} Since for $a\in \{ 1, 2, 7, 15, 23\}$, $\cT_0$ consists 
of a single 
prime, for each $a$ we get exactly two cubic division algebras, with center $\ell = \bQ({\sqrt{-a}})$, and they are 
opposite of each other. Therefore, each 
of the five possible values of $a$ determines a $\bQ$-form $G$ of 
${\rm SU}(2,1)$ uniquely 
(up to a $\bQ$-isomorphism), and for 
$q\notin \cR_{\ell}$, the parahoric subgroup $P_q$ of $G(\bQ_q )$ uniquely (up to conjugation by an element of 
$\overline{G}({\bQ}_q )$, where $\overline{G}$ is the adjoint group of $G$). 
\vskip1mm

We can easily compute $\mu (G(\bR )/\Lambda)$, which, using the volume formula
given in 2.4 is seen to be equal to $\mu e'(P_p)$ when $\cT=\cT_0$, where $(a,p)$ is as
in the preceding theorem, $\mu$ is as in Proposition 3.5, and (see 2.5
({\it ii})) $e'(P_p)= (p-1)^2(p+1)$. We find
that $\mu (G(\bR )/\Lambda)$ equals $1, \ 1, \ 1/7, \ 1, \ 
{\rm and} \ 3$, for $a = 1, \ 2, \ 7, \ 15, \ {\rm and}\ 23$ respectively when $\cT= \cT_0$. 
This computation is clearly independent of the choice of maximal
parahoric 
subgroups $P_q$ in $G(\bQ_q )$ for primes $q$ which ramify in $\ell = \bQ(\sqrt{-a})$. 
\vskip 1mm
In the sequel, the prime $p$ appearing in the pair $(a,p)$ will be called the prime {\it associated} to $a$, and we will sometimes denote it by $p_a$.   
\vskip5mm

\ni
\begin{center}
{\bf 5. The fake projective planes arising from $k = \bQ$}
\end{center}
\vskip2mm

We will show in this section that there are exactly eighteen finite classes (cf.\,1.8) of 
fake projective 
planes with $k =\bQ$. We will explicitly determine their fundamental
groups. 
\vskip1mm

We prove results in 5.2--5.4 for an arbitrary totally real number
field $k$ for applications in \S\S 8 and 9. 

\vskip 1mm

\ni{\bf 5.1.} We will use the notation introduced in 1.2 and 1.3. In particular, 
$k$ is a totally real number field of degree $d$, $\ell$
a totally complex quadratic extension of $k$, and $v_o$ is a real place of
$k$, $G$ is a simple simply connected algebraic $k$-group, which is an
{\it inner}  
form of ${\rm SL}_3$ over $\ell$, such that
$G(k_{v_o})\cong{\rm SU}(2,1)$, and for all real places $v\ne v_o$,
$G(k_v)\cong{\rm SU}(3)$. We recall (1.2) that there is a division algebra
$\cD$ of degree $n|3$, with center $\ell$, $\cD$ given with an
involution $\sigma$ of the second kind so that $\sigma|_{\ell}$ is
the nontrivial $k$-automorphism of $\ell$, and a nondegenerate
hermitian form $h$ on $\cD^{3/n}$ defined in terms of the involution 
$\sigma$,  such that $G$ is the special unitary
group ${\rm SU}(h)$ of the hermitian form $h$. 
\vskip1mm

Let $\cT_0$ be the finite set of nonarchimedean places of $k$ where $G$ is
anisotropic. As pointed out in 2.2, {\it every place $v\in \cT_0$ splits in} 
$\ell$. If $\cD =\ell$, then $h$ is a hermitian form on $\ell^3$, and
$G$ is isotropic at every nonarchimedean place of $k$, so in this case
$\cT_0$ is empty. Now we note that $\cT_0$ is nonempty if $\cD$ is a
cubic division algebra since this division algebra must ramify at
least at two nonarchimedean places of $\ell$.       
\vskip2mm

\ni{\bf 5.2.} Let $C$ be the center of $G$, $\overline{G}$ the adjoint group, and
let $\varphi :\, G\rightarrow {\overline{G}}$ the natural isogeny. Let 
$\cP=(P_v)_{v\in V_f}$ and $\cP' = (P'_v)_{v\in V_f}$ be two coherent collections of 
parahoric subgroups such that 
for all $v\in V_f$,   
$P'_v$ is conjugate to $P_v$ under an element of
  $\overline{G}(k_v)$. For all 
but finitely many $v$, $P_v = P'_v$, and they are hyperspecial. 
Therefore, there is an element 
$g\in {\overline{G}}(A_f )$ such that $\cP'$ is the conjugate of $\cP$ under 
$g$. Let $\overline{P}_v$ be the stabilizer of $P_v$ in $\overline{G}(k_v)$. 
Then ${\overline K} : = \prod_{v\in V_f} {\overline{P}}_v$ is the
  stabilizer 
of $\cP$ in 
$\overline{G}(A_f )$, and it is a compact-open subgroup of the latter. So 
the number of distinct $\overline{G}(k)$-conjugacy classes of coherent  
collections $\cP'$ as above is the cardinality{\footnote {this number is called 
the ``class number'' of $\overline G$ relative to $\overline K$ and is known 
to be 
finite, see for example, Proposition 3.9 of [BP]}} of ${\overline{G}}(k )
\backslash{\overline{G}}(A_f )/{\overline K}$. 
\vskip1mm

As $\varphi :\, G\rightarrow {\overline G}$ is a central isogeny, 
$\varphi (G(A_f ))$ contains the commutator subgroup of ${\overline G}(A_f)$. Moreover, as $G$ is simply connected and $G(k_{v_o})$ is noncompact, by the strong approximation property (Theorem 7.12  
of [PlR]), $G(k)$ is dense in $G(A_f )$, i.\,e., for any open neighborhood 
$\Omega$ of the identity in $G(A_f )$, $G(k)\Omega = G(A_f )$. This implies that ${\overline G}(k){\overline K}$ contains $\varphi (G(A_f))$, which in 
turn contains $[{\overline G}(A_f),{\overline G}(A_f)]$. Hence, ${\overline G}(k){\overline K} ={\overline G}(k)[{\overline G}(A_f),{\overline G}(A_f)]{\overline K}$. Using this observation it is easy to see that ${\overline G}(k){\overline K}$ is a subgroup, and the natural map from 
${\overline G}(k)\backslash{\overline G}(A_f)/{\overline K}$ to the finite abelian group 
${\overline G}(A_f)/{\overline G}(k){\overline K}$ is bijective. We shall next show
that 
this latter group is trivial if $h_{\ell,3} = 1$.  
\vskip1mm

We begin by observing that for every $v\in V_{\infty}$, $H^1(k_v, C)$
vanishes since $C$ is a group of
exponent $3$. Now since by the Hasse principle for 
simply connected 
semi-simple $k$-groups (Theorem 6.6 of [PlR]) $H^1(k, G)\rightarrow 
\prod_{v\in V_{\infty}}H^1 (k_v, G)$ is an 
isomorphism, we conclude that the natural map $H^1(k, C)\rightarrow H^1(k, G)$ is trivial, and hence
the coboundary homomorphism $\delta:\,  {\overline G}(k)\rightarrow  
H^1(k, C)$ {\it is surjective}. 
\vskip1mm
      
Now we note that since for each nonarchimedean place $v$, $H^1(k_v, G)$
is 
trivial ([PlR], Theorem 6.4), the coboundary homomorphism $\delta_v :\,
{\overline G}(k_v)\rightarrow H^1(k_v, C)$ is surjective and its
kernel equals 
$\varphi(G(k_v ))$. Now let $v$ be a nonarchimedean place of $k$ which
either 
does not split in $\ell$, or it splits in $\ell$ and $P_v$ is
an Iwahori subgroup of $G(k_v)$, and $g\in {\overline G}(k_v)$. Then the parahoric subgroup $g(P_v)$ is 
conjugate to $P_v$ under an element of $G(k_v)$, and hence, 
${\overline G}(k_v) =\varphi(G(k_v)){\overline P}_v$, which implies 
that $\delta_v({\overline P}_v)= \delta_v({\overline G}(k_v)) = 
H^1(k_v, C)$. We observe also that for any nonarchimedean place $v$ of
$k$, the subgroup $\varphi
(G(k_v)){\overline P}_v$ is precisely the stabilizer of the type
$\Theta_v$ ($\subset\Delta_v$) of 
$P_v$ under the natural action of ${\overline G}(k_v)$ on $\Delta_v$
described 
in 2.2 of [BP]. Thus $\delta_v({\overline P}_v) = H^1(k_v,
C)_{\Theta_v}$, 
where $H^1(k_v, C)_{\Theta_v}$ is the stabilizer of $\Theta_v$ in 
$H^1(k_v, C)$ under the action of the latter on $\Delta_v$ through
$\xi_v$ 
given in 2.5 of [BP]. It can be seen, but we do not need this fact
here, that 
for 
any nonarchimedean place $v$ of $k$ which does not lie over $3$ and  
$P_v$ is a hyperspecial parahoric subgroup of 
$G(k_v)$, $\delta_v({\overline P}_v)$ equals $H^1_{nr}(k_v, C)$, where 
$H^1_{nr}(k_v, C)$ ($\subset H^1(k_v, C)$) is the ``unramified Galois 
cohomology'' as in [Se3], Chapter II, \S 5.5. 
\vskip1mm     

The coboundary homomorphisms $\delta_v$ combine to provide an 
isomorphism 
$${\overline G}(A_f)/{\overline G}(k){\overline K}\longrightarrow 
\cC :={\prod}' H^1(k_v,
C)/\big(\psi(H^1(k,C))\cdot\prod_{v}\delta_v({\overline P}_v)\big),$$ 
where $\prod' H^1(k_v, C)$ denotes the subgroup of 
$\prod_{v\in V_f} H^1(k_v, C)$ consisting of the elements $c= (c_v)$
such that 
for all but finitely many $v$, $c_v$ lies in $\delta_v({\overline
  P}_v)$, 
and $\psi:\, H^1(k,C)\rightarrow \prod' H^1(k_v, C)$ is the natural homomorphism.
\vskip1mm

Andrei Rapinchuk's remark that $R_{\ell/k}(\mu_3)$ is a direct
product 
of $C$ and (the naturally embedded subgroup) $\mu_3$ has helped us to simplify the following discussion.
\vskip1mm
 
$H^1(k,C)$ can be identified with $\ell^{\times}
/k^{\times}{\ell^{\times}}^3$, and for any place $v$ of $k$, $H^1(k_v, C)$
can be 
identified with $(k_{v}\otimes_{k}\ell )^{\times}/k_v^{\times}
{(k_v\otimes_{k}\ell )^{\times}}^3$. Now let $\cS$ be the finite set of
nonarchimedean places of $k$ which split in $\ell$ and $P_v$ is an
Iwahori subgroup of $G(k_v)$. If $v\notin \cS$ is a nonarchimedean
place which  
splits in $\ell$, and $w'$, $w''$ are the two places of $\ell$ lying over 
$v$, then the subgroup $\delta_v({\overline P}_v)$ 
gets identified with 
$$k_v^{\times}\big({\mathfrak o}_{w'}^{\times}{\ell_{w'}^{\times}}^3\times{\mathfrak o}_{w''}^{\times}{\ell_{w''}^{\times}}^3\big)/k_v^{\times}\big( {\ell_{w'}^{\times}}^3\times{\ell_{w''}^{\times}}^3\big),$$ where 
${\mathfrak o}_{w'}^{\times}$ (resp.,\,\,${\mathfrak o}_{w''}^{\times}$) is the 
group of units of $\ell_{w'}$ (resp.,\,\,$\ell_{w''}$), 
cf.\,\,Lemma 2.3(ii) of [BP] and the proof of Proposition 2.7 in there.

Now let $I^f_{k}$ (resp.,\,$I^f_{\ell}$) be the group of finite id\`eles of $k$ (resp.,\,$\ell$), i.\,e., the restricted direct product of the 
$k_{v}^{\times}$ (resp.,\,\,$\ell_w^{\times}$) for all
nonarchimedean 
places 
$v$ of $k$ (resp.,\,\,$w$ of $\ell$). We shall view $I^f_{k}$ as a subgroup of $I^f_{\ell}$ in terms of its natural embedding. Then it is obvious that $\cC$ is 
isomorphic to the quotient of $I^f_{\ell}$ by the 
subgroup generated by $I^f_{k}\cdot{(I^f_{\ell})}^3\cdot
\ell^{\times}$ and all the elements $x =(x_w)\in I^f_{\ell}$, such
that $x_w\in {\mathfrak o}_w^{\times}$ for every nonarchimedean place
$w$ of $\ell$ which lies over a place of $k$ which splits in $\ell$
but is not in 
$\cS$.  From this it is obvious that $\cC$ is a quotient of the class
group $Cl(\ell)$ of $\ell$, and its exponent is $3$. This implies that
$\cC$ is trivial if $h_{\ell, 3} =1$. 

Let us now assume that 
$\ell =\bQ(\sqrt{-23})$, and $\cS = \{ 2\}$. Then
$h_{\ell,3} = 3$. But as  
either of the two prime ideals lying over $2$ in $\ell =
\bQ({\sqrt{-23}})$ generates the class group of $\ell$, we see that
$\cC$ is again trivial. Thus we have proved the following.
\vskip1mm

\ni{\bf 5.3. Proposition.} {\it Let $\cP =(P_v)_{v\in V_f}$ and $\cP' = (P'_v)_{v\in V_f}$ be two
  coherent collections of parahoric subgroups such that for every $v$,
  $P'_v$ is conjugate to $P_v$ under an element of $\overline{G}(k_v)$. 
Then there is an element in
  ${\overline{G}}(k)$ which conjugates $\cP'$ to $\cP$ if  
$h_{\ell,3} =1$. This is also the case if $\ell = 
\bQ(\sqrt{-23})$, and the set $\cS$ of rational primes $p$ which 
split in $\ell$ and $P_p$ is an Iwahori subgroup consists of $2$ alone.}     
\vskip2mm

\ni{\bf 5.4.} Let $G$, $C$, and $\overline{G}$ be as in 5.1. As before, let $\cT_0$ be the finite
set of nonarchimedean places of $k$ where $G$ is anisotropic. 

We fix a coherent collection $(P_v)_{v\in V_f}$ of parahoric 
subgroups such that $P_v$ is maximal for every $v$ which splits in $\ell$.  Let $\Lambda
= G(k)\bigcap \prod_v P_v$, $\Gamma$ be the normalizer of $\Lambda$ in
$G(k_{v_o})$, and $\overline\Gamma$ be the image of $\Gamma$ in
${\overline G}(k_{v_o})$. We know (see bound $(0)$ in 2.3, and 2.2)
that $[\Gamma : \Lambda ]\leqslant 3^{1+\#\cT_0} h_{\ell, 3}$. From Proposition 2.9 of [BP] and a careful analysis of 
the arguments in 5.3, 
5.5 and the proof of Proposition 0.12 of {\it loc.\,cit.}\,\,it can be
deduced that, in fact,  
$[\Gamma : \Lambda] =3^{1+\#\cT_0}$, {\it if either} $h_{\ell,3}=1$ (then $h_{k,3} = 1 $, see Theorem 4.10 of {\rm [W]}), 
{\it or $(k,\ell) =(\bQ, \bQ(\sqrt{-23}))$ and $\cT_0 = \{ 2\}$.} We briefly outline the proof below.
\vskip1mm

Let $\Theta_v$ ($\subset\Delta_v$) be the type of $P_v$, and $\Theta =
\prod\Theta_v$. We have observed in 5.2 that the coboundary
homomorphism $\delta:\, {\overline{G}}(k)\rightarrow H^1(k,C)$ is
surjective. Using this fact we find that, for $G$ at
hand, 
the last term $\delta ({\overline G}(k))'_{\Theta}$ 
in the short exact sequence of Proposition 2.9 of [BP], for $G' = G$, 
coincides with the subgroup $H^1(k, C)_{\Theta}$ of $H^1(k, C)$
defined in 2.8 of [BP]. Also, the order of the first term of that short 
exact sequence is $3/\#\mu(\ell)_3$. So to prove the assertion about
$[\Gamma : \Lambda]$, it would suffice to show that $H^1(k, C)_{\Theta}$
is of order $\#\mu(\ell)_3 3^{\#\cT_0}$ if either $h_{\ell,3} = 1$, or
$(k,\ell)=(\bQ,\bQ(\sqrt{-23}))$ and $\cT_0 =\{ 2 \}$. 
\vskip1mm

As in 2.1, let $\ell^{\bullet} = \{ x\in \ell^{\times}\, |\, N_{\ell/k}(x)\in {k^{\times}}^3\}$, and identify $H^1(k,C)$ with $\ell^{\bullet}/{\ell^{\times}}^3$. Let $\ell_3$ (resp.,\,\,$\ell^{\bullet}_{\cT_0}$) be the subgroup of $\ell^{\times}$ (resp.,\,\,$\ell^{\bullet}$) consisting of elements $x$ such that for every normalized valuation $w$ of $\ell$ (resp.,\,\,every normalized valuation $w$ of $\ell$ which does not lie over a place in $\cT_0$), $w(x)\in 3\bZ$. Let $\ell^{\bullet}_3= 
\ell_3\cap\ell^{\bullet}.$ We can identify $H^1(k, C)_{\Theta}$ with 
the group 
${\ell^{\bullet}_{\cT_0}}/{\ell^\times}^3$, see 2.3,
2.7 and 5.3--5.5 of [BP]. We claim that the order  
of $\ell^{\bullet}_{\cT_0}/{\ell^{\times}}^3$ is $\#\mu(\ell)_3 3^{\#\cT_0}$. If 
$h_{\ell,3}=1=h_{k,3}$, from
2.3 above and Proposition 0.12 of [BP] we see that 
 $\#\ell^{\bullet}_3 /{\ell^{\times}}^3 =\#\mu(\ell)_3$, and 
$U(k)/U(k)^3\rightarrow k_3/{k^{\times}}^3$ is an
isomorphism. Since the homomorphism $U(\ell)/U(\ell)^3\rightarrow
U(k)/U(k)^3$, induced by the norm map, is onto (2.3), 
given an element $y\in \ell^{\times}$ whose norm lies in $k_3$,  we can find an
element $u\in U(\ell)$ such that $uy\in \ell^{\bullet}$, i.\,e.,
$N_{\ell/k}(uy)\in {k^{\times}}^3$. Now it is easy to see that if 
$h_{\ell,3} =1$, $\ell^{\bullet}_{\cT_0}/\ell^{\bullet}_3$ is of order
$3^{\#\cT_0}$. This implies that
$\#\ell^{\bullet}_{\cT_0}/{\ell^{\times}}^3
=\#\mu(\ell)_3 3^{\#\cT_0}$. 
 
Let  $(k,\ell) = (\bQ, \bQ(\sqrt{-23}))$ now. Then, as neither of the two prime ideals of $\ell =\bQ({\sqrt{-23}})$ 
lying over $2$ is a principal ideal, we see that 
$\ell^{\bullet}_{\{ 2\}} =\ell^{\bullet}_3$. But since the class number of $\bQ$ is $1$, $\ell^{\bullet}_3 =\ell_3$, and therefore, 
$\ell^{\bullet}_{\{ 2\}}/{\ell^{\times}}^3 = 
\ell_3 /{\ell^{\times}}^3$. The latter group is of order $3$ 
(=$h_{\ell,3}$) since $\bQ(\sqrt{-23})$ does not contain a nontrivial cube root of unity, see the proof  of Proposition 0.12 in [BP]. This
proves the assertion that $[\Gamma : \Lambda]= 3^{1+\#\cT_0}$ if
either $h_{\ell,3} = 1$, or $(k,\ell) = (\bQ,\bQ(\sqrt{-23}))$ and
$\cT_0 =\{ 2\}$.      
 
\vskip2mm

\ni{\bf 5.5.} In the rest of this section we will assume that $k =\bQ$ and $\cD$ is a cubic division algebra with center $\ell =\bQ(\sqrt{-a})$ given with an involution $\sigma$ of the second kind. Let $G$ be the simple simply connected $\bQ$-subgroup of ${\rm GL}_{1,\cD}$ such that  for any commutative $\bQ$-algebra $A$, 
$$G(A) = \{\,z\in {\rm GL}_{1,\cD}(A)= (A\otimes_{\bQ} \cD)^{\times}\, | \, z\sigma (z)= 1 \  {\rm and}\ {\rm Nrd}(z) =1 \}.$$  
  
  \vskip2mm

\ni {\bf 5.6. Lemma.} {\it $G(\bQ )$ is torsion-free if $a\ne 3$ or $7$. If $a =3$ (resp.,\,\,$a =7$), then the order of any nontrivial element of $G(\bQ )$ of finite order is $3$ (resp.,\,\,$7$).}
\vskip1mm
\ni {\it Proof.} Let $x\in G(\bQ )$ be a nontrivial element of finite order. Since the reduced norm of $-1$ in $\cD$ is $-1$, $-1\notin G(\bQ)$. Therefore, the order of $x$ is odd, and the  
${\bQ}$-subalgebra $K :={\bQ}[x]$ of $\cD$ generated by $x$ is a nontrivial 
field extension of $\bQ$. Note that the degree of any field extension
of $\bQ$ contained in $\cD$ is a divisor of $6$.  If $K=\ell$, 
then $x$ lies in the center of $G$, and hence it is of order $3$. But 
${\bQ}(\sqrt{-3})$ is the field generated by a nontrivial cube root of unity. 
Hence, if $K =\ell$, then $a=3$ and $x$ is of order $3$. Let us assume now 
that $K\ne \ell$. Then $K$ cannot be a quadratic extension of $\bQ$
since if it is a quadratic extension, $K\cdot\ell$ is a field  extension of $\bQ$ of degree $4$
contained in $\cD$, which is not possible. So $K$ is   
an extension of $\bQ$ of degree either 3 or 6. Since an extension of 
degree 3 of $\bQ$ cannot contain a root of unity different from $\pm 1$, 
$K$ must be of degree 6, and so, in particular, it contains $\ell = 
\bQ({\sqrt{-a}})$. Note that the only roots of unity of odd order which can be contained in an extension of $\bQ$ of 
degree 6 are the 7-th or the 9-th roots of unity. 

    For an integer $n$, let $C_n$ be the extension of $\bQ$ generated by a primitive  
    $n$-th root  $\zeta_n$  of unity. Then $C_7 = C_{14}\supset \bQ ({\sqrt{-7}})$, 
and $C_9 =C_{18}\supset C_3 = \bQ ({\sqrt{-3}})$, and 
$\bQ ({\sqrt{-7}})$ (resp.,\,\,$\bQ ({\sqrt{-3}})$) is the only quadratic 
extension of $\bQ$ contained in $C_7$ (resp.,\,\,$C_9$). As $K\supset 
\bQ ({\sqrt{-a}})$, we conclude that the group $G(\bQ )$ is 
torsion-free if $a\ne 3$ or $7$, and if $a=3$ (resp.,\,\,$a =7$), then the order of $x$ is $9$ (resp.,\,\,$7$). In particular, if $a = 3$ (resp.,\,\,$a=7$), then $K=\bQ(\zeta_9)$ (resp.,\,\,$K=\bQ(\zeta_7)$). However, if $a =3$, $N_{K/{\ell}}(\zeta_9) =\zeta_9^3\ne 1$, and if $a =7$, $N_{K/{\ell}}(\zeta_7) =1$. This implies the last assertion of the lemma.  
\vskip2mm

\ni{\bf 5.7.} Let $a$ be one of the following five integers: $1$, $2$,
$7$, $15$, and $23$, and let $p=p_a $ be the prime associated to $a$ (see
4.4--4.5). Let $\ell = \bQ(\sqrt{-a})$. Let $\cD$ be a cubic division algebra with center $\ell$ whose local invariants at the two places of $\ell$ lying over $p$ are nonzero and negative of each other, and the local  
invariant at all the other places of $\ell$ is zero. (There are two
such division algebras, they are opposite of each other.) Then $\bQ_p\otimes_{\bQ}\cD =(\bQ_p\otimes_{\bQ}\ell )\otimes_{\ell}\cD ={\mathfrak D}\times{\mathfrak D}^o$, where $\mathfrak D$ is a cubic division algebra with center $\bQ_p$, and ${\mathfrak D}^o$ is its opposite. $\cD$ admits an involution $\sigma$ of the second kind.  Let the simple simply connected $\bQ$-group $G$ be as in 5.5. 
We may (and do) assume that   $\sigma$ is so chosen that $G(\bR)\cong {\rm
  SU}(2,1)$.  We observe that any other such involution of $\cD$, or of 
  its opposite, similarly determines a $\bQ$-group which is $\bQ$-isomorphic to $G$ (1.2). As $\sigma({\mathfrak D}) = {\mathfrak D}^o$, it is easily seen that $G(\bQ_p)$ is the compact group ${\rm SL}_1({\mathfrak D})$ of elements of reduced norm $1$ in $\mathfrak D$.  
\vskip1mm

We fix a coherent collection $(P_q)$ of maximal parahoric subgroups $P_q$ of $G(\bQ_q)$ which are hyperspecial  for every rational prime $q\ne p$ which does not ramify in $\ell$. Let $\Lambda =G(\bQ)\cap\prod_qP_q$, and let $\Gamma$ be its normalizer in $G(\bR)$. Let $\overline\Gamma$ be the image of $\Gamma$ in ${\overline G}(\bR)$.
  
\vskip2mm 

\ni{\bf 5.8. Proposition.} {\it If $(a,p) = (23,2)$, then $\overline\Gamma$ is 
torsion-free.} 
\vskip2mm
\ni{\it Proof.} We assume that $(a,p) =(23,2)$, and begin by observing that 
$\overline\Gamma$ is contained in 
${\overline G}(\bQ)$, see, for example, Proposition 1.2 of [BP]. Since $H^1(\bQ, C)$ is a group of exponent $3$, so is the 
group ${\overline G}(\bQ)/\varphi (G(\bQ))$. Now as $G(\bQ )$ is torsion-free 
(5.6), any nontrivial 
element of 
${\overline G}(\bQ)$ of finite order has order $3$. 
 
To be able to describe all the elements of order $3$ of ${\overline G}(\bQ)$, 
we consider the connected reductive $\bQ$-subgroup $\cG$ of ${\rm GL}_{1,\cD}$, which contains $G$ as a 
normal subgroup, such that for any commutative $\bQ$-algebra $A$,  
$$\cG (A) = \{\, z\in {\rm GL}_{1,\cD}(A)= (A\otimes_{\bQ}\cD)^{\times}\,\, 
|\,\, z\sigma(z)\in A^{\times}\}.$$
Then the center $\cC$ of $\cG$ is $\bQ$-isomorphic to 
$R_{\ell/\bQ}({\rm GL}_1)$. The conjugation action of $\cG$ on $G$ induces a $\bQ$-isomorphism 
$\cG/\cC \rightarrow {\overline G}$. As 
$H^1(\bQ, \cC) = \{ 0\}$, we conclude that the natural homomorphism
$\cG (\bQ) \rightarrow {\overline G}(\bQ)$ is surjective. 
Now given an element $g$ of $\cG(\bQ)$ whose image in  ${\overline G} (\bQ)$ is an element of order $3$, $\lambda := g^3$  lies in  $\ell^{\times}$. Let $a=g\sigma(g)\in \bQ^{\times}$.  Then $(i)$  $\lambda \sigma(\lambda)= a^3$. Let $L = \ell[X]/(X^3-\lambda)$ and let $x$ be the unique cube root of $\lambda$  in $L$. There is a unique embedding of $L$ in $\cD$ over $\ell$ which maps $x$ to $g$. 
The reduced norm of $g$ is clearly  
$\lambda$, and the image of $g$ in $H^1(\bQ, C)\subset \ell^{\times}
/{\ell^{\times}}^3$ is the class of $\lambda$ in $\ell^{\times}/{\ell^{\times}}^3$. Now 
if $g$ stabilizes the collection $(P_q)$, then its image in $H^1(\bQ, C)$ 
must lie in the subgroup $H^1(\bQ, C)_{\Theta}$, and hence, $(ii)$ 
$w(\lambda)\in 3\bZ$ for every normalized valuation $w$ of $\ell$ not lying 
over $2$ (cf.\,\,5.4).

The conditions $(i)$ and $(ii)$ imply that $\lambda \in
\ell^{\bullet}_{\{2\}} =\ell_3 = \bigcup_i \alpha^i {\ell^{\times}}^3$ (cf.\,\,5.4), where
$\alpha = (3+{\sqrt{-23}})/{2}$. Since $\lambda$ is not a cube in
$\ell$, $\lambda\in
\alpha{\ell^{\times}}^3\cup\alpha^2{\ell^{\times}}^3$. But 
$\bQ_2$ contains a cube root of 
$\alpha$ (this can be seen using Hensel's Lemma), and hence for $\lambda \in \alpha{\ell^{\times}}^3\cup \alpha^2{\ell^{\times}}^3$, 
$L =\ell[X]/(X^3-\lambda)$ is not embeddable in $\cD$. 
(Note that $L$ is embeddable in $\cD$ if, and only if,  ${\bQ_2}\otimes_{\bQ} L$ 
is a 
direct product of two field extensions of $\bQ_2$, both of degree $3$.)
Thus we have shown that ${\overline G}(\bQ)$ does not contain any nontrivial 
elements of 
finite order which stabilize the collection $(P_q)$. Therefore, 
$\overline\Gamma$ is torsion-free.        
\vskip2mm

\ni{\bf 5.9. Examples of fake projective planes.} By Lemma 5.6 the subgroup $\Lambda$ described in 5.7 
is torsion-free if $(a, p) = (1,5),\ (2,3)$, $(15,2)$ or $(23,2)$. Now 
let $(a,p) = (7,2)$. Then $G(\bQ_2)$ is the group ${\rm SL}_1({\mathfrak D})$ of elements of reduced norm $1$ in a cubic division algebra ${\mathfrak D}$ with center $\bQ_2$ (cf.\,5.7). The first congruence subgroup $G(\bQ_2)^+:= {\rm SL}_1^{(1)}({\mathfrak D})$ of ${\rm SL}_1({\mathfrak D})$ is the unique maximal normal 
pro-$2$ subgroup of $G(\bQ_2)$ of index $(2^3-1)/(2-1)=7$ (see Theorem 7(iii)(2) of 
[Ri]). By the strong approximation property (Theorem 7.12 of [PlR]), $\Lambda^+ := \Lambda\cap G(\bQ_2)^+$ is a subgroup of $\Lambda$ of index $7$. Lemma 5.6 implies that 
$\Lambda^+$ is torsion-free since $G(\bQ_2)^+$ is a pro-$2$ group. As $\mu (G(\bR )/\Lambda ) = 1/7$ (see 4.5), $\mu (G(\bR )/\Lambda^+)= 1$, and hence the Euler-Poincar\'e characteristic 
of $\Lambda^+$ is 3. 

Since $\Lambda$, and for $a= 7$, $\Lambda^+$ are congruence subgroups, 
according to Theorem 15.3.1 of [Ro], $H^1(\Lambda, \bC )$, and for $a =7$, 
$H^1(\Lambda^+, \bC )$ vanish. By Poincar\'e-duality, then $H^3(\Lambda, \bC )$, and for $a =7$, $H^3(\Lambda^+, \bC )$ also vanish. For $a = 1,\ 2,\ {\rm and}\ 15$, as $\mu (G(\bR)/{\Lambda}) = 1$ (4.5), the Euler-Poincar\'e characteristic $\chi (\Lambda)$ of $\Lambda$ 
is 3, and for $a =7$, $\chi (\Lambda^+)$ is also 3, we conclude that for $a =
1,\ 2,\ {\rm and}\ 15$, $H^i(\Lambda ,\bC )$ is 1-dimensional for $i = 0,\ 2, \ {\rm and} \ 4$, and if $a =7$, this is also the case for $H^i(\Lambda^+, \bC)$.
Thus if $B$ is the symmetric space of $G(\bR)$, then {\it for 
$a = 1$, $2$ and $15$, $B/\Lambda$, and for $a = 7$, $B/\Lambda^+$, is a fake projective plane.}  
\vskip1mm

Let $\overline\Lambda$ (resp.,\,${\overline\Lambda}^+$) be the image of $\Lambda$ (resp.,\,$\Lambda^+$) in ${\overline G}(\bR)$. There is a natural faithful action of ${\overline\Gamma}/{\overline\Lambda}$ (resp.,\,${\overline\Gamma}/{\overline\Lambda}^+$), which is a group of order $3$ (resp.,\,\,$21$), on 
$B/\Lambda$ (resp.,\,\,$B/\Lambda^+$). As $\overline\Gamma$ is the normalizer of $\Lambda$, and also of $\Lambda^+$, in ${\overline G}(\bR)$, ${\overline\Gamma}/{\overline\Lambda}$ (resp.,\,${\overline\Gamma}/{\overline\Lambda}^+$) is the full automorphism group of $B/\Lambda$ (resp.,\,
$B/{\Lambda}^+$). 

\vskip2mm
In 5.10--5.13, we will describe the classes of fake projective 
planes associated with each of the five pairs $(a,p)$.
\vskip1.5mm
  
\ni {\bf 5.10.} In this paragraph we shall study the fake projective planes 
arising from the pairs $(a,p) =(1,5)$, $(2,3)$, and $(15, 2)$. Let us first consider the fake projective planes 
with $\cT = \cT_0$ (for the pair $(15,2)$, $\cT =\cT_0$ is the only possibility, see Theorem 4.4).  
Let $\Lambda$ and $\Gamma$ be as in 5.7. Let $\Pi\subset 
\overline\Gamma$ be the fundamental group of a fake projective plane and 
$\widetilde\Pi$ be its inverse image in $\Gamma$. Then as $1 = \chi ({\widetilde\Pi}) = 
3\mu (G(\bR)/{\widetilde\Pi})= \mu (G(\bR )/\Lambda )$, $\widetilde\Pi$ is of 
index $3$ (\,=\, $[\Gamma : \Lambda]/3$) 
in $\Gamma$, and hence $\Pi$ is a torsion-free subgroup of $\overline\Gamma$ 
of index $3$. Conversely, if $\Pi$ is a torsion-free 
subgroup of $\overline\Gamma$ of index $3$ such that $H^1(\Pi, \bC)
= \{ 0 \}$ (i.\,e., $\Pi/[\Pi ,\Pi]$ is finite), then as 
$\chi(\Pi) = 3$, $B/{\Pi}$ is a fake projective plane, and $\Pi$ is its 
fundamental group. 

\vskip1mm

Let us now assume that $(a,p) = (1,5)$, or $(2,3)$ and $\cT \ne \cT_0$. Then, by Theorem 4.4,  $\cT= \cT_0\cup\{2\}$.  Note that $2$ is the only prime which ramifies in $\bQ(\sqrt{-a})$, $a =1$ or $2$. We fix an Iwahori subgroup $I_2$ contained in $P_2$ and let $\Lambda = G(\bQ)\cap\prod_q P_q$, $\Lambda_I = \Lambda\cap I_2 $.  Then $\Lambda$, and so also $\Lambda_I$, is torsion-free (Lemma 5.6). Since $[P_2:I_2] =3$, the strong approximation property implies that $[\Lambda:\Lambda_I] = 3$. As $\chi(\Lambda) = 3$, we obtain that $\chi(\Lambda_I ) = 9$.  Now let $\Gamma_I$ be the normalizer of $\Lambda_I$ in $G(\bR)$, and $\overline{\Lambda}_I$ and $\overline{\Gamma}_I$ be the images of $\Lambda_I$ and $\Gamma_I$ in $\overline{G}(\bR)$.  Then $[\Gamma_I:\Lambda_I] = 9$ (cf.\,5.4), and hence, $[\overline{\Gamma}_I:\overline{\Lambda}_I] =3$. This implies that the orbifold Euler-Poincar\'e characteristic of $\overline{\Gamma}_I$ is $3$. Moreover, $H^1(\overline{\Lambda}_I,\bC)$, and hence also $H^1(\overline{\Gamma}_I,\bC)$ vanishes (Theorem 15.3.1 of [Ro]).  Thus if $\overline{\Gamma}_I$ is torsion-free, which indeed is the case, as can be seen by a suitable adaptation of the argument used to prove Proposition 5.8, $B/\overline{\Gamma}_I$ is a fake projective plane, and it is the unique plane belonging to the class associated to $\ell = \bQ(\sqrt{-a})$, for $a = 1,\,2$, and $\cT = \cT_0\cup \{2\}$.

 \vskip1mm
 
 It follows from Proposition 5.3 that up to conjugation by $\overline{G}(\bQ)$ there is exactly one coherent collection $\{ P_q, q\ne 2,\,p;\, I_2, \,P_p =G(\bQ_p) \}$ such that $P_q$ is a hyperspecial parahoric subgroup of $G(\bQ_q)$ for all $q\ne 2,\,p$. Thus for $a=1,\,2$, up to conjugacy under $\overline{G}(\bQ)$ we obtain a unique $\overline{\Gamma}_I$.
 \vskip1mm

\vskip2mm

\ni{\bf 5.11.} We will now study the fake projective planes arising from 
the pair $(7,2)$. For this pair, either $\cT =\cT_0=\{2\}$, or $\cT=\{2,3\}$, or $\cT = \{2,5\}$. We will describe first the fake projective planes with $\cT = \cT_0$. As in 5.9, let $\Lambda^+ = \Lambda\cap G(\bQ_2)^+$, 
which is a torsion-free subgroup of $\Lambda$ of index $7$. We know that $B/\Lambda^+$ is a fake projective plane. Let $\widetilde\Pi$ be the inverse image in $\Gamma$ of the fundamental group 
${\Pi}\subset{\overline\Gamma}$ of a fake projective plane. Then as 
$\mu (G(\bR)/\Gamma)=\mu(G(\bR)/\Lambda)/9 = 1/63$, and $\mu (G(\bR)/{\widetilde\Pi})=\chi ({\widetilde\Pi})/3 = 1/3$, $\widetilde\Pi$ is a subgroup of 
$\Gamma $ of index $21$, and hence $[\overline\Gamma :\Pi] =21$. Conversely, 
if $\Pi$ is a torsion-free subgroup of $\overline\Gamma$ of index $21$, then 
as $\chi({\Pi}) =3$, $B/{\Pi}$ is a fake 
projective plane if, and only if,    
$\Pi/[\Pi , \Pi ]$ is finite. 

\vskip1mm

Now let $\cT =\{2,3\}$ or $\{2,5\}$. We fix a coherent collection $(P_q)$ of {\it maximal} parahoric subgroups $P_q$ of $G(\bQ_q)$ such that $P_q$ is hyperspecial if, and only if, $q\notin \cT\cup \{7\}$. We will denote the principal arithmetic subgroup $G(\bQ)\cap \prod_q P_q$ by $\Lambda_3$ if $\cT =\{2, 3\}$, and by $\Lambda_5$ if $\cT = \{2,5\}$. Let $\Gamma_3$ and  $\Gamma_5$ be the normalizers in $G(\bR)$  of $\Lambda_3$ and $\Lambda_5$ respectively. Let ${\overline{\Lambda}}_3$, ${\overline{\Lambda}}_5$, ${\overline{\Gamma}}_3$ and ${\overline{\Gamma}}_5$ be the images in ${\overline{G}}(\bR)$ of 
$\Lambda_3$, $\Lambda_5$, $\Gamma_3$ and $\Gamma_5$ respectively.
\vskip1mm

We will now describe the fake projective planes arising from the pair $(7,2)$ with $\cT = \{2,3\}$. Since  $e'(P_2) = 3$, and $e'(P_3)= 7$, see 3.5 and 2.5({\it ii}),\,({\it iii}), $$\mu(G(\bR)/\Lambda_3) = \frac{1}{21}e'(P_2)e'(P_3) = 1.$$ Hence, the orbifold Euler-Poincar\'e characteristic $\chi(\Lambda_3)$ of $\Lambda_3$ is $3$. As the maximal normal pro-3 subgroup of the non-hyperspecial maximal parahoric subgroup $P_3$ of $G(\bQ_3)$ is of index $96$, $P_3$ does not contain any elements of order $7$.  But any nontrivial element of $G(\bQ)$ of finite order is of order $7$ (Lemma 5.6),  so we conclude that $\Lambda_3$ ($\subset G(\bQ)\cap P_3$) is torsion-free. As in 5.9, using Theorem 15.3.1 of [Ro] we conclude that $B/\Lambda_3$ is a fake projective plane, $\Lambda_3$ ($\cong{\overline\Lambda}_3$) is its fundamental group, and since ${\overline\Gamma}_3$ is the normalizer of ${\overline\Lambda}_3$ in $\overline{G}(\bR)$, ${\overline{\Gamma}}_3/{\overline{\Lambda}}_3$ is the automorphism group of $B/\Lambda_3$.  As in 5.10, we see that any torsion-free subgroup $\Pi$ of ${\overline\Gamma}_3$ of index $3$ such that $\Pi/[\Pi,\Pi]$ is finite is the fundamental group of a fake projective plane.  
\vskip1mm

We will now treat the case where $\cT = \{2, 5\}$. As $e'(P_5) = 21$, see 2.5({\it iii}), $$\mu(G(\bR)/\Lambda_5) =\frac{1}{21}e'(P_2)e'(P_5) = 3.$$ Hence, $\chi({\overline\Gamma}_5) = 3\chi(\Gamma_5) = 9\mu(G(\bR)/\Gamma_5) = 9\mu(G(\bR)/\Lambda_5)/[\Gamma_5 :\Lambda_5] = 3$. From this we conclude that the only subgroup of ${\overline\Gamma}_5$ which can be the fundamental group of a fake projective plane is ${\overline\Gamma}_5$ itself. Moreover, as $H^1({\overline\Lambda}_5, \bC)$, and hence also $H^1({\overline\Gamma}_5, \bC)$, are trivial (Theorem 15.3.1 of [Ro]), $B/{\overline\Gamma}_5$ is a fake projective plane, and ${\overline\Gamma}_5$ is its fundamental group, if and only if, ${\overline\Gamma}_5$ is torsion free. 
\vskip1mm

We will now show, using a variant of the argument employed in the proof of Proposition 5.8, that ${\overline\Gamma}_5$ is torsion-free. Since the maximal normal pro-5 subgroup of the non-hyperspecial maximal parahoric subgroup $P_5$ of $G(\bQ_5)$ is of index $720$, $P_5$ does not contain any elements of order $7$. This implies that $\Lambda_5$ ($\subset G(\bQ)\cap P_5$), and hence also ${\overline\Lambda}_5$, are  torsion-free since any element of $G(\bQ)$ of finite order is of order $7$ (Lemma 5.6). Now as ${\overline\Lambda}_5$ is a normal subgroup of ${\overline\Gamma}_5$ of index $3$, we conclude that the order of any nontrivial element of ${\overline\Gamma}_5$ of finite order is $3$. 
\vskip1mm

Let $\cG$ be the connected reductive $\bQ$-subgroup of ${\rm GL}_{1,\cD}$, which contains $G$ as a normal subgroup, such that for any commutative $\bQ$-algebra $A$, $$\cG(A) = \{ z\in {\rm GL}_{1,\cD}(A)= (A\otimes_{\bQ}\cD)^{\times}\, |\, z\sigma(z)\in A^{\times}\}.$$ Then the center $\cC$ of $\cG$ is $\bQ$-isomorphic to $R_{\ell/\bQ}({\rm GL}_1)$. The conjugation action of $\cG$ on $G$ induces a $\bQ$-isomorphism $\cG/\cC\to {\overline G}$. As $H^1(\bQ, \cC) =\{0\}$, the natural homomorphism $\cG(\bQ)\to {\overline G}(\bQ)$ is surjective. Now, if possible, assume that ${\overline\Gamma}_5$ contains an element of order $3$. We fix a  $g\in \cG(\bQ)\,(\subset {\cD}^{\times})$ whose image $\overline{g}$ in $\overline{G}(\bQ)$ is an element of order $3$ of ${\overline\Gamma}_5$. Then $\lambda :={g}^3$ lies in $\ell^{\times}$. Let $a = g\sigma(g)\in \bQ^{\times}$. Then the norm of $\lambda$ (over $\bQ$) is $a^3\in {\bQ^{\times}}^3$. Let $x$ be the unique cube root of $\lambda$ in the field 
$L=\ell[X]/(X^3-\lambda)$. Then there is an embedding of $L$ in $\cD$ which maps $x$ to $g$. We will view $L$ as a field contained in $\cD$ in terms of this embedding.  The reduced norm of $g$ is clearly $\lambda$, and the image of $\overline{g}$ in $H^1(\bQ, C)\subset \ell^{\times}/{\ell^{\times}}^3$ is the class of $\lambda$ in $\ell^{\times}/{\ell^{\times}}^3$.    As in the proof of Proposition 5.8 (cf.\:also 5.4), we see that since $g$ stabilizes the collection $(P_q)$, $w(\lambda)\in 3\bZ$ for every normalized valuation $w$ of $\ell =\bQ(\sqrt{-7})$ not lying over $2$.    
\vskip1mm

We assert that the subgroup $\ell^{\bullet}_{\{2\}}$ of $\ell^{\times}$ consisting of elements $z$ whose norm lies in ${\bQ^{\times}}^3$, and $w(z)\in 3\bZ$ for every normalized valuation $w$ of $\ell$ not lying over $2$, equals 
${\ell^{\times}}^3 \cup (1+\sqrt{-7}){\ell^{\times}}^3 \cup (1-\sqrt{-7}){\ell^{\times}}^3$.  Since $\ell$ does not contain a nontrivial cube root of unity, and its class number is $1$, the subgroup $\ell_3$ of $\ell^{\times}$ consisting of $z\in \ell^{\times}$ such that for {\it every} normalized valuation $w$ of $\ell$, $w(z)\in 3\bZ$ coincides with ${\ell^{\times}}^3$ (see the proof of Proposition 0.12 in [BP]), and $\ell_3$ is of index $3$ in the subgroup $\ell^{\bullet}_{\{2\}}$ (of $\ell^{\times}$).  As $(1+\sqrt{-7})(1-\sqrt{-7}) = 8\, (\in {\ell^{\times}}^3)$, it follows that $1+\sqrt{-7}$ and $1-\sqrt{-7}$ are units at every nonarchimedean place of $\ell$ which does not lie over $2$. Moreover, it is easy to see that  if $v'$ and $v''$ are the two normalized valuations of $\ell$ lying over $2$, then neither $v'(1+\sqrt{-7})$ nor $v''(1+\sqrt{-7})$ is a multiple of $3$. This implies, in particular, that $1\pm \sqrt{-7}\notin {\ell^{\times}}^3$.  From these observations, the above assertion is obvious. Now we note that $1\pm\sqrt{-7}$ is not a cube in $\bQ_5(\sqrt{-7})$ (to see this, it is enough to observe, using a direct computation, that $(1\pm\sqrt{-7})^8\ne 1$ in the residue field of $\bQ_5(\sqrt{-7})$). Since $\lambda\in \ell^{\bullet}_{\{2\}}$ and is not a cube in $\ell^{\times}$, it must lie in the set $(1+\sqrt{-7}){\ell^{\times}}^3\cup (1-\sqrt{-7}){\ell^{\times}}^3$. But no element of this set is a cube in $\bQ_5(\sqrt{-7})$. Hence, $\mathfrak{L} :=L\otimes_{\ell}{\bQ_5(\sqrt{-7})}$ is an unramified field extension of $\bQ_5(\sqrt{-7})$ of degree $3$.  
\vskip1mm

Let $T$ be the centralizer of $g$ in $G$. Then $T$ is a maximal $\bQ$-torus of $G$. Its group of $\bQ$-rational points is $L^{\times}\cap G(\bQ)$. The torus $T$ is anisotropic over $\bQ_5$ and its splitting field over $\bQ_5$ is the unramified cubic extension $\mathfrak{L}$ of $\bQ_5(\sqrt{-7})$. This implies  that any parahoric subgroup of $G(\bQ_5)$ containing $T(\bQ_5)$ is hyperspecial. We conclude from this that $T(\bQ_5)$ is contained in a unique parahoric subgroup of $G(\bQ_5)$, and this parahoric subgroup is hyperspecial. According to the main theorem of [PY], the subset of points fixed by $g$ in the Bruhat-Tits building of $G(\bQ_5)$ is the building of $T(\bQ_5)$. Since the latter consists of a single point, namely the vertex fixed by the hyperspecial parahoric subgroup containing $T(\bQ_5)$, we infer that $g$ normalizes a unique parahoric subgroup of $G(\bQ_5)$, and this parahoric subgroup is hyperspecial.  As $P_5$ is a non-hyperspecial maximal parahoric subgroup of 
  $G(\bQ_5)$, it cannot be normalized by $g$. Thus we have arrived at a contradiction. This proves that ${\overline\Gamma}_5$ is torsion-free. Hence, $B/{\overline{\Gamma}}_5$ is a fake projective plane, and ${\overline\Gamma}_5$  is its fundamental group. Since the normalizer of ${\overline\Gamma}_5$ in $\overline{G}(\bR)$ is ${\overline\Gamma}_5$, the automorphism group of $B/{\overline{\Gamma}}_5$ is trivial. 
\vskip1mm

\ni{\bf 5.12.} We finally look at the fake projective planes arising from the pair $(23,2)$. In this case, $\mu (G(\bR)/\Gamma) = \mu (G(\bR)/\Lambda)/9 
=1/3$ (see 4.5). Hence, if $\widetilde\Pi$ is the inverse image in $\Gamma$ of the fundamental group ${\Pi} \subset {\overline\Gamma}$ of a fake projective plane, 
then as $\mu(G(\bR)/{\widetilde\Pi})= \chi({\widetilde\Pi})/3 = 1/3  =
\mu (G(\bR)/\Gamma)$, ${\widetilde\Pi} =\Gamma$. Therefore, the only subgroup 
of $\overline \Gamma$ which can be the 
fundamental group of a fake projective plane is $\overline\Gamma$ itself. 

As $\overline\Gamma$ is torsion-free (Proposition 5.8), $\chi({\overline\Gamma}) =3$, and $\Lambda/[\Lambda, \Lambda]$, hence $\Gamma/[\Gamma, \Gamma]$, and so also ${\overline\Gamma}/[{\overline\Gamma},{\overline\Gamma}]$ are finite, $B/{\overline\Gamma}$ is a fake projective plane and $\overline\Gamma$ is its fundamental group. Since the normalizer of $\overline\Gamma$ in ${\overline G}(\bR)$ equals $\overline\Gamma$, the automorphism group of $B/{\overline\Gamma}$ is trivial.
\vskip2mm

\ni{\bf 5.13.} We recall that if $v$ does not split in $\ell$, and is unramified in $\ell$, then the non-hyperspecial maximal parahoric subgroups of $G(k_v)$ are conjugate to each other under $G(k_v)$. Using the observations in 2.2, and Proposition 5.3, we see that if $(a,p)\ne (15,2)$ (resp.,\,$(a,p)=(15,2)$), then up to conjugation by $\overline{G}(\bQ)$, there are exactly 2 (resp.,\,4) coherent collections  $(P_q)$ of maximal parahoric subgroups such that $P_q$ is hyperspecial whenever $q$ does not ramify in $\bQ(\sqrt{-a})$ and $q\ne p$, since if $a\ne 15$ (resp.,\,$a=15$), there is exactly one prime (resp.,\,there are exactly two primes, namely 3 and 5) which ramify in $\ell = \bQ(\sqrt{-a})$. Moreover, for $(a,p) = (1,5)$ and $(2,3)$, up to conjugation by $\overline{G}(\bQ)$, there is exactly one coherent collection $(P_q)$ of parahoric subgroups such that $P_q$ is hyperspecial for $q\ne 2,\, p$, and $P_2$, $P_p$ are Iwahori subgroups;  for $(a,p) = (7,2)$, if either $\cT=\{2,3\}$ or $\{2,5\}$, then up to conjugation by $\overline{G}(\bQ)$, there are exactly 2 coherent collections $(P_q)$ of maximal parahoric subgroups such 
that $P_q$ is hyperspecial if, and only if, $q\notin \cT\cup \{7\}$.

\vskip1mm

From the results in 5.10--5.12,  we conclude that for $(a,p) $ equal to either $(1,5)$ or  $(2,3)$, there are two distinct classes with $\cT =\{p\}$, and one more class with $\cT = \{2,p\}$; for $(a,p) =(23, 2)$, there are two distinct classes;  for $(a,p) =(7,2)$, there are six distinct finite classes, and for $(a,p)=(15,2)$, there are four distinct finite classes, of fake projective planes. Thus the following theorem holds.

\vskip2mm

\ni {\bf 5.14. Theorem.} {\it There exist exactly eighteen distinct classes of fake projective 
planes with $k =\bQ$}.
\vskip2mm

\ni{\bf 5.15. Remark.}  To the best of our knowledge, only three fake 
projective planes were known before the present work. The first 
one was constructed by Mumford [Mu] and it corresponds to the pair 
$(a,p) = (7,2)$; see 5.11. Two more examples were given by 
Ishida and Kato [IK] making use of the discrete subgroups of 
${\rm PGL}_3(\bQ_2)$, which act simply transitively on the set of vertices of 
the Bruhat-Tits building of the latter, constructed by Cartwright, Mantero, 
Steger and Zappa. In both of these examples, $(a,p)$ equals $(15,2)$. JongHae Keum 
has recently constructed a fake projective 
plane in [Ke] which is birational to a cyclic cover of degree $7$ of a Dolgachev surface. This fake projective plane admits an automorphism of order $7$, so it appears to us that it  corresponds to 
the pair $(7,2)$, and its fundamental group is the group $\Lambda^+$ of 5.9 for a suitable choice of a maximal 
parahoric subgroup $P_7$ of $G(\bQ_7 )$.

\vskip3mm

\ni
\begin{center}
{\bf 6. Lower bound for discriminant in terms of the degree of a number field}  
\end{center}
\vskip2mm

\ni{\bf 6.1. Definition.} We define $M_r(d)={\rm min}_K D_K^{1/d},$
where the minimum is taken over all
totally real number fields $K$ of degree $d.$  Similarly, we define
$M_c(d)={\rm min}_K D_K^{1/d}$ by taking the minimum over all
totally complex number fields $K$ of degree $d.$
\vskip2mm

It is well-known that $M_r(d)\geqslant ({d^d}/{d!})^{2/d}$ from 
the classical estimates of Minkowski.
The precise values of $M_r(d)$ for small values of $d$ are known due to the work of many mathematicians as listed in [N]. For $d\leqslant 8$, the values of $M_r(d)$ are given in the following table.  

$$\begin{array}{cccccccc}
d&2&3&4&5&6&7&8\\
M_r(d)^d
&5
&49
&725
&14641
&300125
&20134393
&282300416.
\end{array}
$$

\vskip1mm

An effective lower bound for $M_r(d)$, better than Minkowski's bound 
for $d$ large, has been given by 
Odlyzko [O1].  We recall the following 
algorithm given in [O1], Theorem 1, which provides a useful estimate for 
$M_r(d)$ for 
arbitrary $d.$

\vskip1mm

\ni{\bf 6.2.} Let $b(x)=[5+(12 x^2-5)^{1/2}]/6.$  Define 
\begin{eqnarray*}
g(x,d)&=&\exp \big[ \log(\pi)-\frac{\Gamma^\prime}{\Gamma}(x/2)+
{\frac{(2 x-1)}{4}} \big(\frac{\Gamma^\prime}{\Gamma}\big)^\prime (b(x)/2)\\
&&+\frac{1}{d}\{-\frac{2}{x}-\frac{2}{x-1}-\frac{2x-1}{{b(x)}^2} -\frac{2 x-1}{(b(x)-1)^2}\}\big].
\end{eqnarray*}
Let $\alpha=\sqrt{\frac{14-\sqrt{128}}{34}}.$  
As we are considering only totally real number fields, according to 
Theorem 1 of [O1], 
$M_r(d)\geqslant g(x,d)$ provided that $x>1$ and $b(x)\geqslant 1+\alpha x.$  

Now let $x_0$ be the positive 
root of the quadratic equation $b(x)=1+\alpha x.$  Solving this equation, we obtain 
$x_0=\frac{\alpha+\sqrt{2-5\alpha^2}}{2(1-3\alpha^2)}=1.01...$. For a fixed value of $d$, define ${\mathfrak N}(d)=\limsup_{x\geqslant x_0}g(x,d).$
\vskip2mm

\ni {\bf 6.3. Lemma.} {\it For each $d>1,$ $M_r(d)\geqslant {\mathfrak N}(d)$, and ${\mathfrak N}(d)$ is an increasing function of
$d.$}
\vskip2mm

\ni{\it Proof}. It is obvious from our
choice of $x_0$ that $M_r(d)\geqslant {\mathfrak N}(d).$ We will now show that ${\mathfrak N}(d)$ is an increasing function of $d$. 

For a fixed value of $x>1$, $g(x,d)$ is clearly an
increasing function of $d$ since the only expression involving $d$ in it is
$$\frac{1}{d}\{-\frac{2}{x}-\frac{2}{x-1}-\frac{2x-1}{{b(x)}^2} -\frac{2 x-1}{(b(x)-1)^2}\},$$ which is 
nonpositive. Now for a given $d$, and a positive integer $n$, choose a 
$x_n\geqslant x_0$ such that $g(x_n,d)\geqslant {\mathfrak N}(d)-10^{-n}$. Then $${\mathfrak N}(d+1)=\limsup_{x\geqslant x_0}g(x,d+1)\geqslant g(x_n,d+1)\geqslant g(x_n,d)\geqslant {\mathfrak N}(d)-10^{-n}.$$
Hence, $\mathfrak N(d+1)\geqslant {\mathfrak N}(d)$. 

\vskip2mm
\ni {\bf 6.4.}  In the next section,    
we will use the lower bound for the root-discriminant $D_K^{1/d}$ of
totally complex number fields $K$ 
obtained by Odlyzko in [O2]. We will denote by $N_{c}(n_0)$ the entry 
for 
totally complex number fields given in the last column of Table 2 of
[O2] for $n = n_0$. We recall from [O2] that for every number field $K$ of
degree $n\geqslant n_0$, the root-discriminant $D_K^{1/n}>N_c(n_0)$.   

\vskip1mm

For small $d$, we will also use Table IV
of [Ma]. This table was originally constructed by
Diaz y Diaz.  

\vskip4mm
\ni
\begin{center}
{\bf 7. Upper bounds for the degree $d$ of $k$, $D_k$ and $D_{\ell}$}  
\end{center}
\vskip4mm

In this, and the next two sections, we will determine totally real number fields $k$ of degree $d>1$, their totally complex quadratic extensions $\ell$, $k$-forms $G$ of ${\rm SU}(2,1)$ and coherent collections $(P_v)_{v\in V_f}$ of parahoric subgroups $P_v$ of $G(k_v)$ such that for all $v\in \cR_{\ell}$, $P_v$ is maximal, and the image $\overline\Gamma$ in ${\overline G}(k_{v_o})$ (where $v_o$ is the unique real place of $k$ such that $G(k_{v_o})\cong {\rm SU}(2,1)$) of the normalizer $\Gamma$ of $\Lambda := G(k)\cap \prod_{v\in V_f}P_v$ in $G(k_{v_o})$ contains a torsion-free subgroup $\Pi$ of finite index with $\chi(\Pi) =3$.{\it Then $\chi(\Gamma)$ is a reciprocal integer. In particular, it is $\leqslant 1$.} 
\vskip1mm
   
In this section, we will use bounds (2), (3), (6), and (7)--(10) obtained in \S 2, the lower bound for the discriminant given in the preceding section, and Hilbert class fields, to prove that $d\leqslant 5$. We will also find good upper bounds for $D_k$, $D_{\ell}$, and $D_{\ell}/D_k^2$ for $d\leqslant 5$. Using these bounds, in the next section we will be able to make a complete list of $(k,\ell)$ of interest to us. It will follow then that $d$ cannot be $5$.
\vskip3mm
   
\ni {\bf 7.1.} Let $f(\delta,d)$ be the function occurring in bound $(10)$. It is obvious that for $c>1$, $c^{{1}/{(3-\delta )d}}$ decreases as $d$ increases. Now for $\delta \geqslant 0.002$, as 
$${\frac{\delta (1+\delta )}{0.00136}} >1,$$ $\inf_{\delta}f(\delta,d)$, where the infimum is taken over the closed interval 
$0.002\leqslant \delta\leqslant 2$, decreases as 
$d$ increases. A direct computation shows that $f(0.9,20)<16.38.$  On
the other hand, for $d\geqslant 20$, Lemma 
6.3 gives us 
$$M_r(d)\geqslant {\mathfrak N}(20)
\geqslant g(1.43,20)> 16.4,$$ where $g(x,d)$ is
the function defined in 6.2. From these bounds we conclude that 
$\ \ d=[k:{\mathbb Q}]< 20$.
\vskip1mm

To obtain a better upper bound for $d$, 
we observe using Table 2 in [O2] that $M_r(d)>17.8$ for 
$15\leqslant d<20$. But by a direct computation we see that 
$f(0.9,15)<17.4.$  So the monotonicity of $f(\delta,d)$, as a function 
of $d$ for a fixed $\delta$, implies that 
$d$ {\it cannot be larger than} $14$.

\vskip2mm

\ni{\bf 7.2.}  Now we will prove that 
$d\leqslant 7$ with the help of Hilbert class fields. 
Let us assume, 
if possible, that $14\geqslant d\geqslant 8$. 
\vskip1mm

We will use the following result from the theory of Hilbert class
fields. 
The Hilbert class field $L := H(\ell )$ of a totally complex number
field 
$\ell$ is the maximal unramified abelian extension of $\ell$. {\it Its degree over $\ell$ is the class number $h_{\ell}$ of $\ell$, and $D_{L} = D_{\ell}^{h_{\ell}}$.}
\vskip1mm

 We consider the 
two cases where  
$h_{\ell}\leqslant 63$ and $h_{\ell}> 63$ separately.
\vskip1mm

\ni{\it Case $(a)$}: $h_{\ell}\leqslant 63$:   In this case
$h_{\ell,3}\leqslant 27$, and from bound 
$(8)$ we obtain 
$$D_k^{1/d}< \varphi_2(d,h_{\ell,3})<\varphi_3 (d):= 27^{{1}/{4d}}(16\pi^5)^{1/{4}}.$$
The function $\varphi_3(d)$ decreases as $d$ increases. 
A direct computation shows
that $\varphi_3(d)\leqslant \varphi_3(8) <9.3.$  Hence, $D_k^{1/d}< 9.3$. On the other hand, from 
Table 2 in [O2] we find that, for $14\geqslant d\geqslant 8$, $M_r(d)>10.5$, so $D_k^{1/d}>10.5$.  Therefore, if $h_{\ell}\leqslant 63$, $d\leqslant 7$.  
\vskip2mm

\ni{\it Case $(b)$}: $h_{\ell}> 63$:  In this case, let $L$ be the
Hilbert 
class field of $\ell$. Then 
$[L :\ell]=h_\ell,$ $D_{L} =D_{\ell}^{h_{\ell}}$,  
and $2dh_{\ell}>16\times 63>1000.$  From 6.4 we conclude that 
$$D_{\ell}^{1/2d}=D_{L}^{1/2dh_{\ell}}\geqslant M_c(2dh_{\ell})\geqslant N_c(1000)
=20.895,$$
where the last value is from Table 2 of [O2].  However, as
$f(0.77,d)\leqslant f(0.77,8)<20.84$, bound $(10)$ implies that  
$D_{\ell}^{1/2d}<20.84.$  Again, we have reached a contradiction. 
So we conclude that $d\leqslant 7.$
\vskip2mm

\ni{\bf 7.3.}  To find good upper bounds for $d$, $D_k$ and $D_{\ell}$, we will make use of improved lower bounds for $R_{\ell}/w_{\ell}$ for totally complex number fields given in [F], Table 2. We reproduce below the part of this table which we will use in this paper. 

$$\begin{array}{ccc}
r_2= d& {\mbox {for}} \ D_\ell^{1/2d}<\ \ \ &R_{\ell}/w_{\ell}\geqslant\\
2&17.2&0.0898\\
3&24.6&0.0983\\
4&29.04&0.1482\\
5&31.9&0.2261\\
6&33.8&0.4240\\
7&34.4&0.8542
\end{array}
$$

We also note here that except for totally complex sextic fields of discriminants
$$-9747,\ -10051,\ -10571,\ -10816,\ -11691,\ -12167,$$
 and totally complex quartic
fields of discriminants 
$$117,\ 125,\ 144,$$ $R_{\ell}/w_{\ell}$ is bounded
from below by $1/8$ for every number field $\ell$, see [F], Theorem B$'$.

\vskip2mm

\ni{\bf 7.4.} We consider now the case where $d=7$. Bound $(10)$ implies
that $D_\ell^{1/14}< f(0.75,7)<22.1.$
 Using the lower bound for 
$R_{\ell}/w_{\ell}$ given in the table above and bound $(7),$ we conclude by a direct computation that 
$$D_\ell^{1/14}<\varphi_1(7,0.8542,0.8)<
18.82.$$ 
On the other hand, the root-discriminant of
any totally complex number field of degree $\geqslant 260$ is bounded
from 
below by $N_c(260)$, see 6.4.  From Table 2 in [O2] we find that $N_c(260)=18.955.$ So we conclude that the class number $h_\ell$ of $\ell$ is bounded from 
above by
$260/2d=260/14 <19,$ for otherwise the root-discriminant of the
Hilbert class 
field of $\ell$ would be greater than $18.955,$ contradicting the fact
that it equals $D_\ell^{1/14}$ ($<18.82$).  
\vskip1mm

As $h_{\ell}\leqslant 18$, $h_{\ell,3}\leqslant 9.$  Now we will use 
bound $(8)$. We see by a direct computation that $\varphi_2(7,9)< 9.1.$
Hence, $D_k^{1/7}\leqslant D_\ell^{1/14}< 9.1.$ On the other hand, we know 
from 6.1 that
$M_r(7)=20134393^{1/7}\geqslant11.$  This implies that $d$ cannot be
$7$. Therefore, $d\leqslant 6$. 
\vskip2mm

\ni{\bf 7.5.} Employing a method similar to the one used in 7.2 and 
7.4 we will now show that $d$ cannot be $6$.  
\vskip1mm

For $d=6,$ from bound $(10)$ we get $D_\ell^{1/12}<
f(0.71,6)<24.$ 
Using the lower bound for $R_{\ell}/w_{\ell}$ provided by the table in
7.3 and 
bound $(7),$ we conclude by a direct computation that
$D_\ell^{1/12}<\varphi_1(6,0.424,0.8)<20.$ From Table 2 in
[O2] 
we find that  $N_{c}(480)> 20$. Now, arguing as in 7.4, we infer that
the 
class number $h_\ell$ of $\ell$ is bounded from above by
$480/12=40$, which implies that $h_{\ell,3}\leqslant 27.$  As
$\varphi_2(6,27)<10$, bound $(8)$ implies that $D_k^{1/6}\leqslant
D_\ell^{1/12}<10 .$ Now since $N_{c}(21)>10$, we see that the class
number 
of $\ell$ cannot be larger than $21/12<2.$ Hence,  
$h_{\ell} =1=h_{\ell,3}.$  We may now apply bound $(8)$ again to conclude
that
$D_k^{1/6}<\varphi_2(6,1)<8.365.$  
Checking from the table t66.001 of [1], we know that the two smallest discriminants of
totally real sextics are $300125$ as mentioned in 6.1, followed by $371293$.
As $371293^{1/6}>8.47,$ the second case is not possible and we are left with
only one candidate, $D_k=300125.$
As ${\mathfrak p}(6,300125,1)<1.3,$ we conclude from bound $(9)$ that 
$D_\ell/D_k^2=1.$  Hence, if $d =6$,  $(D_k,D_\ell)=(300125,300125^2)$ is the
only possibility. From the tables in [1] we find that there is a
unique totally real number field $k$ of degree $6$ with $D_k =
300125$. Moreover, the class number of this field is $1$. Gunter
Malle, using the procedure described in 8.1 below, has shown that
there does not exist a totally complex quadratic extension $\ell$ of this field
with $D_{\ell} = 300125^2$. Therefore $d$ cannot be $6$.  
\vskip2mm

\ni{\bf 7.6.} For $d=5,$  bound $(10)$ implies that $D_\ell^{1/10}<
f(0.7,5)<26.1.$ It is seen from the table in 7.3 that 
$R_\ell/w_\ell\geqslant 0.2261$.
Hence, $D_\ell^{1/10}<\varphi_1(5,0.2261,0.72)<21.42.$ As 
$N_c(2400)> 21.53$, arguing as in 7.4 we see that the class number
$h_{\ell}$ of $\ell$ is bounded from above by 
$2400/10=240.$ Hence, $h_{\ell,3}\leqslant 81=3^4.$  Now we note that
$\varphi_2(5,81)<10.43$,  but  
$N_c(23)> 10.43.$ So, $h_\ell< 23/10$, and therefore, 
$h_{\ell,3}=1.$  But then 
$D_k^{1/5}\leqslant D_{\ell}^{1/10}<\varphi_2(5,1)<8.3649.$  
As $M_r(5)^5\geqslant 14641$ and ${\mathfrak p}(5,14641,1)<5.2,$ we
conclude from bound $(9)$ that $D_\ell/D_k^2\leqslant 5.$
\vskip2mm

\ni{\bf 7.7.}  Let now $d=4.$  In this case, $k$ is a totally
real 
quartic and $\ell$ is a totally complex
octic containing $k$. Table 4 of [F] gives the lower bound
$R_k\geqslant 41/50$
for the regulator.  Since $\ell$ is a CM field which is a totally complex quadratic
extension of $k,$ we know that $R_\ell=2^{d-1}R_k/Q,$ where $Q=1$ or $2$ is the unit 
index of $k$ (cf.\,[W]).  We will now estimate $w_\ell$, the number
of roots of unity in $\ell.$  
\vskip1mm

We know that the group of roots of unity in $\ell$ 
is a cyclic group of even order, say $m$. Let $\zeta_m$ be a primitive
$m$-th root of unity. As the degree of the cyclotomic field $\bQ(\zeta_m)$
is $\phi(m),$ where $\phi$ is the Euler function, we know that
$\phi(m)$ is a divisor of $2d=8.$ The following table gives the
values of $m$ and $\phi(m)$ for 
$\phi(m)\leqslant 8$.
$$\begin{array}{ccccccccccccc}
m&2&4&6&8&10&12&14&16&18&20&24&30\\
\phi(m)&1&2&2&4&4&4&6&8&6&8&8&8.
\end{array}$$
\vskip1mm

If $\phi(m)=8$,  then $m= 16$, $20$, $24$ or $30$, and 
$\bQ(\zeta_m)$ equals $\ell.$   Note that $\bQ(\zeta_{30})=\bQ(\zeta_{15}).$   The class number of these four 
cyclotomic fields are all known to be $1$ (see [W], pp.\,230 and
432).  So 
in these
four cases, $h_{\ell,3}=1.$  Bound $(8)$  
implies that 
$D_k^{1/4}\leqslant D_\ell^{1/8}<\varphi_2(4,1)<8.3640.$  
As $M_r(4)^4 = 725$ and ${\mathfrak p}(4,725,1)<21.3,$ we
conclude from bound $(9)$ that $D_\ell/D_k^2\leqslant 21.$
\vskip1mm

Assume now that $\phi(m)\ne 8$. Then $m\leqslant 12$.
Hence,  
$w_\ell\leqslant12.$  So we conclude that except for the four
cyclotomic fields dealt with earlier,  
$$R_\ell/w_\ell\geqslant 2^3R_k/12Q\geqslant R_k/3\geqslant 41/150.$$ 

Applying bound $(7)$, we conclude that
$D_\ell^{1/8}< \varphi_1(4,41/150,0.69)<21.75$ by a 
direct computation.  From Table IV of [Ma], we know that
totally complex number fields of degree $\geqslant 4000$ have unconditional
root-discriminant lower bound $21.7825.$
It follows, as before, using the Hilbert Class
field of $\ell$, that the class number $h_{\ell}$ of $\ell$ is at most 
$4000/8=500.$  Hence,  
$h_{\ell,3}\leqslant 3^5=243.$ Bound $(8)$ now gives that 
$D_\ell^{1/8}< \varphi_2(4,243)<11.8.$ But from Table 2 of [O2]
we find that $N_c(32)> 11.9.$ 
So we conclude $h_\ell\leqslant 32/8 =4.$  Hence, 
$h_{\ell,3}\leqslant 3.$  Applying bound $(8)$ again we infer that    
$D_\ell^{1/8}<\varphi_2(4,3)<8.96.$ As $N_c(18)> 9.2$, we conclude
that $h_\ell< 18/8.$  But then $h_{\ell,3}=1$, and the argument in the
preceding paragraph leads to the conclusion that 
$D_k^{1/4}\leqslant D_\ell^{1/8}<\varphi_2(4,1)<8.3640$ 
and $D_\ell/D_k^2\leqslant 21.$
\vskip2mm

\ni{\bf 7.8.}  We consider now the case $d=3.$
Suppose that $D_\ell^{1/6}<21.7.$  Since according to Table IV
of [Ma], $M_c(4000)\geqslant  21.7825,$ we infer, as above, using the 
Hilbert class field of $\ell$, that $h_{\ell}\leqslant
4000/6<667$. Then $h_{\ell,3}\leqslant 243=3^5.$  It follows from bound $(8)$ that
$D_\ell^{1/6}<\varphi_2(3,243)<13.3.$  From Table 2
of [O2] we find that $N_c(44)>13.37.$  Therefore, 
$h_\ell\leqslant 44/6 <8.$  Hence, $h_{\ell,3}\leqslant 3.$ 
Now we observe that 
$\varphi_2(3,3)<9.17$. But as  $N_c(18)\geqslant 9.28$,
$h_\ell<18/6=3$, 
which implies that $h_{\ell,3}=1.$ We then deduce from bound $(8)$ that 
$D_k^{1/3}\leqslant D_{\ell}^{1/6}<\varphi_2(3,1)<8.3591$.  Also since
$D_k\geqslant 49$ (see 6.1), and   
${\mathfrak p}(3,49,1)<52.8,$ we
conclude from bound $(9)$ that $D_\ell/D_k^2\leqslant 52.$
\vskip1mm

We assume now that $D_\ell^{1/6}\geqslant21.7$ (and $d =3$). We will
make use of a lower bound for $R_\ell/w_\ell$ which is better than the 
one provided in 7.3.  Table 4 of [F] gives that $R_k\geqslant 0.524.$
Recall from 7.7 that $R_\ell/w_\ell=2^{d-1}R_k/Qw_\ell\geqslant
2R_k/w_\ell\geqslant 2(0.524)/w_\ell.$ From the
table 
of values of the Euler function given in 7.7, we see that $\phi(m)$ is a
proper divisor of $6$ only for $m= 2$, $4$, $6$.  So we conclude 
that $w_\ell\leqslant 6$ unless $\ell$ is either $\bQ(\zeta_{14})$ or
$\bQ(\zeta_{18}).$  Since both $\bQ(\zeta_{14})=\bQ(\zeta_{7})$ or
$\bQ(\zeta_{18})=\bQ(\zeta_{9})$ are known to have class number $1$
(cf.\,\,[W], pp.\,229 and 412), 
the bounds obtained in the last paragraph apply to these two cases as
well.
 Hence, it remains only to consider the cases where $w_\ell\leqslant
 6$. So we assume now that $w_{\ell}\leqslant 6$. Then 
$R_\ell/w_\ell\geqslant 2(0.524)/6 > 0.17.$
\vskip1mm
   
   Observe that bounds (2), (3) and (6) imply that 
$$D_k^{1/d}>\xi(d,D_\ell,R_\ell/w_\ell,\delta):=\big[{\frac{(R_{\ell}/w_{\ell})\zeta(2d)^{1/2}}{\delta(\delta+1)}}\big]^{1/d}\frac{(2\pi)^{1+\delta}}{16\pi^5\Gamma(1+\delta)\zeta(1+\delta)^2}(D_\ell^{1/2d})^{4-\delta}.$$

As $D_\ell^{1/6}\geqslant21.7$,
it follows from this bound by a direct computation that 
$D_k^{1/3}>\xi(3,21.7^6,0.17,0.65)>16.4.$
 \vskip1mm
 
Recall now a result of Remak, stated as bound (3.15) in [F],
$$R_k\geqslant \big[\frac{\log D_k-d\log d}{\{\gamma_{d-1}d^{1/(d-1)}(d^3-d)/3\}^{1/2}}\big]^{d-1},$$
where $d=3$ and $\gamma_{2}=2/\sqrt3$ as given on page 613 in [F]. 
Since $R_\ell=2^2R_k/Q\geqslant 2R_k,$ we obtain the following
lower bound 
$$R_\ell/w_\ell\geqslant r(D_k,w_\ell):=\frac2{w_\ell}
\big[\frac{\log D_k-3\log 3}{\{2(3^2-1)\}^{1/2}}\big]^{2}.$$
As in the argument in the last paragraph, we assume that $w_{\ell}\leqslant 6.$ 
Then from the preceding bound we get the following:    
$$R_\ell/w_\ell\geqslant r(16.4^3,6)>0.54.$$
We now use  bound (7) to conclude that 
$D_\ell^{1/6}<\varphi_1(3,0.54,0.66)<20.8<21.7,$
contradicting our assumption that $D_\ell^{1/6}\geqslant
21.7.$ 

Therefore, 
$D_k^{1/3}\leqslant D_{\ell}^{1/6}<8.3591$ and $D_\ell/D_k^2\leqslant 52.$
\vskip2mm

\ni{\bf 7.9.}
 Finally we consider the case $d=2.$ 
In this case, we know from 7.3 that $R_\ell/w_\ell\geqslant1/8$ except
in the three cases mentioned there.  So bound 
(7) implies that $D_\ell^{1/4}<\varphi_1(2,1/8,0.52)<
28.96.$  Hence, $D_\ell\leqslant 703387$. This bound holds for the 
three exceptional cases of 7.3 as well.
 Since quartics of such small discriminant
are all known, we know the class number of all such fields explicitly.
In particular, the number fields are listed in t40.001-t40.057 of [1], 
where each file contains $1000$ number fields listed in ascending order
of the absolute discriminants.  There are altogether $5700$ number
fields in the files, the last one has discriminant $713808$. So [1] is
more 
than adequate for our purpose.  Inspecting by hand, or
using PARI/GP and a simple program, we find that the largest class
number of an $\ell$ with  $D_\ell\leqslant 703387$ is $64.$  The corresponding 
number field has discriminant $654400$ with a defining polynomial
$x^4-2x^3+27x^2-16x+314.$
\vskip1mm

Once we know that $h_\ell\leqslant 64,$ we find that
$h_{\ell,3}\leqslant 27.$  We may now apply
bound (8) to conclude that
$D_\ell^{1/4}<\varphi_2(2,27)< 12.57.$ Now since in Table 2 of [O2] we find that $N_c(38)> 12.73,$ we infer that $h_\ell<38/4<10$, which implies that $h_{\ell,3}\leqslant 9.$  
But $\varphi_2(2,9)<10.96$, and $N_c(26)>11.01$. So 
$h_\ell<26/4<7$, and hence $h_{\ell,3}\leqslant 3.$  
It follows from bound $(8)$ that $D_k^{1/2}\leqslant
D_\ell^{1/4}<\varphi_2(2,3)<9.5491.$ From this we conclude
that $D_k\leqslant 91.$ As $D_k\geqslant 5$ (see 6.1) and 
${\mathfrak p}(2,5,3)<104.2,$ bound $(9)$
implies that $D_\ell/D_k^2\leqslant 104.$

\vskip2mm

\ni{\bf 7.10.} The results in 7.6--7.9 are summarized in the following table.

$$\begin{array}{cccc}
{\ \ d \ \ }&D_k^{1/d}\leqslant D_\ell^{1/2d}\leqslant\ \ &{\ \ h_{\ell,3}
\leqslant \ \ } &D_\ell/D_k^2\leqslant\\
5&8.3649&1&5\\
4&8.3640&1&21\\
3&8.3591&1&52\\
2&9.5491&3&104
\end{array}
$$
\vskip5mm

\ni
\begin{center}
{\bf 8. $(k,\ell)$ with  $d=2,\ 3,\ 4,$ and $5$} 
\end{center}

\ni{\bf 8.1.}  To make a list of all pairs $(k,\ell)$ of interest to us, we
will make use of the 
tables of number fields given in [1]. In the following table, in the column 
under $r_d$ (resp.,\,\,$c_d$) we
list the largest integer less than the $d$-th power (resp.,\,\,an
integer slightly larger than the $2d$-th power) of the numbers
appearing in
the second column of the table in 7.10. The column under $x_d$
reproduces the numbers appearing in the last column of the table in
7.10. Therefore, we need only find 
all totally real number fields $k$ of degree $d$, $2\leqslant d\leqslant
5$, and totally complex quadratic extensions $\ell$ of each $k$, 
such that $D_k\leqslant r_d$, $D_{\ell}\leqslant c_d$, and moreover,
$D_{\ell}/D_k^2\leqslant x_d$. Thanks to a detailed computation
carried out at our request by Gunter Malle,
for each $d$, we know the exact number of pairs of $(k,\ell)$ satisfying these 
constraints.  This number is listed in the last column of the
following table. The data is obtained in the following way.  The number fields $k$ with
$D_k$ in the range we are interested in are listed in [1]. Their class
numbers, and a set of generators of their group of units, are also 
given there. For $d=2,$ the quadratic
extensions $\ell$ are 
also listed in [1].  Any quadratic extension of $k$ is of
the 
form $k(\sqrt{\alpha})$, with $\alpha$ in the ring of integers 
${\mathfrak o}_k$ of $k$. For $d>2$, the class number of any totally real 
$k$ of interest turns out to be  $1$; hence, ${\mathfrak o}_k$ is a unique 
factorization 
domain. Now using factorization of small primes and explicit
generators of 
the group of units of $k$, Malle listed all possible $\alpha$ modulo 
squares, and then for each of the $\alpha$, the discriminant of 
$k(\sqrt{\alpha})$ could be 
computed.  Using this procedure, Malle explicitly determined all
totally 
complex quadratic 
extensions $\ell$ with $D_{\ell}$ satisfying the conditions mentioned
above.

$$\begin{array}{ccccc}
d&r_d&c_d&x_d&\#(k,\ell)\\
5&40954&17 \times 10^8&5&0\\
4&4893&24 \times 10^6&21&7\\
3&584&35 \times 10^4&52&4\\
2&91&8320&104&52
\end{array}$$
\ni Thus there are no $(k,\ell)$ with  $d=5$. 
For $2\leqslant d\leqslant 4,$ 
there are $52+4+7=63$ pairs $(k,\ell)$ satisfying the 
constraints on $r_d, c_d$ and $x_d$ imposed by the considerations  
in 7.6--7.9.
\vskip1mm

\ni{\bf 8.2.}  For each of the $63$ potential pairs $(k,\ell)$
mentioned 
above, we know defining polynomials for $k$ and $\ell$, and also the
values of  $D_k$, 
$D_{\ell}$, and $h_{\ell,3}$. It turns out that $h_{\ell,3}=1$ or
$3.$ We are able to further cut down the list of pairs $(k,\ell)$ such that there is a $k$-form of ${\rm SU}(2,1)$, described in terms of the quadratic extension $\ell$ of $k$, which may provide an arithmetic subgroup $\Gamma$ of ${\rm SU}(2,1)$ with $\chi(\Gamma)\leqslant 1$, by making use of
bound $(9)$ for $D_{\ell}/D_k^2$, and the fact that this number is an integer.  We are then left with only $40$
pairs. These are listed below.

\vskip1mm
In the lists below, there are only three pairs $(k,\ell)$ with $d=3$. In the list provided by Malle there was a fourth pair with $(D_k, D_{\ell}, h_{\ell}) = (321, 309123, 1)$. Bound $(9)$ for this pair gives us $D_{\ell}/D_k^2 <2.7$, and therefore, $D_{\ell}\leqslant 2D_k^2$. But $309123> 2\times 321^2$, that is why the fourth pair with $d=3$ does not appear in the lists below.  

$$\begin{array}{cll}
{\ \ \ (k,\ell)\ \ \ }&k&\ell\\
\cC_1&x^2 - x - 1&x^4 - x^3 + x^2 - x + 1\\
\cC_2&x^2 - x - 1&x^4 - x^3 + 2x^2 + x + 1\\
\cC_3&x^2 - x - 1&x^4 + 3x^2 + 1\\
\cC_4&x^2 - x - 1&x^4 - x^3 + 3x^2 - 2x + 4\\
\cC_5&x^2 - x - 1&x^4 - x^3 + 5x^2 + 2x + 4\\
\cC_6&x^2 - x - 1&x^4 - 2x^3 + 6x^2 - 5x + 5\\
\cC_7&x^2 - x - 1&x^4 + 6x^2 + 4\\
\cC_8&x^2 - 2&x^4 + 1\\
\cC_9&x^2 - 2&x^4 + 2x^2 + 4\\
\cC_{10}&x^2 - 2&x^4 - 2x^3 + 5x^2 - 4x + 2\\
\cC_{11}&x^2 - 3&x^4 - x^2 + 1\\
\cC_{12}&x^2 - 3&x^4 + 4x^2 + 1\\
\cC_{13}&x^2 - x - 3&x^4 - x^3 + 4x^2 + 3x + 9\\
\cC_{14}&x^2 - x - 3&x^4 - x^3 + 2x^2 + 4x + 3\\
\cC_{15}&x^2 - x - 4&x^4 - x^3 - 2x + 4\\
\cC_{16}&x^2 - x - 4&x^4 - x^3 + 5x^2 + 4x + 16\\
\cC_{17}&x^2 - x - 5&x^4 - x^3 - x^2 - 2x + 4\\
\cC_{18}&x^2 - 6&x^4 - 2x^2 + 4\\
\cC_{19}&x^2 - 6&x^4 + 9\\
\cC_{20}&x^2 - 7&x^4 - 3x^2 + 4\\
\cC_{21}&x^2 - x - 8&x^4 - x^3 - 2x^2 - 3x + 9\\
\cC_{22}&x^2 - 11&x^4 - 5x^2 + 9\\
\cC_{23}&x^2 - 14&x^4 - 2x^3 + 9x^2 - 8x + 2\\
\cC_{24}&x^2 - x - 14&x^4 - x^3 - 4x^2 - 5x + 25\\
\cC_{25}&x^2 - 15&x^4 - 5x^2 + 25\\
\cC_{26}&x^2 - 15&x^4 - 7x^2 + 16\\
\cC_{27}&x^2 - x - 17&x^4 - x^3 - 5x^2 - 6x + 36\\
\cC_{28}&x^2 - 19&x^4 - 9x^2 + 25\\
\cC_{29}&x^2 - x - 19&x^4 + 9x^2 + 1\\
\cC_{30}&x^2 - 22&x^4 - 2x^3 + 11x^2 - 10x + 3\\
\cC_{31}&x^3 - x^2 - 2x + 1&x^6 - x^5 + x^4 - x^3 + x^2 - x + 1\\
\cC_{32}&x^3 - x^2 - 2x + 1&x^6 - x^5 + 3x^4 + 5x^2 - 2x + 1\\
\cC_{33}&x^3 - 3x - 1&x^6 - x^3 + 1\\
\cC_{34}&x^4 - x^3 - 4x^2 + 4x + 1&x^8 - x^7 + x^5 - x^4 + x^3 - 
x + 1\\
\cC_{35}&x^4 - 5x^2 + 5&x^8 - x^6 + x^4 - x^2 + 1\\
\cC_{36}&x^4 - 4x^2 + 2&x^8 + 1\\
\cC_{37}&x^4 - 4x^2 + 1&x^8 - x^4 + 1\\
\cC_{38}&x^4 - 2x^3 - 7x^2 + 8x + 1&x^8 - 3x^6 + 8x^4 - 3x^2 + 1\\
\cC_{39}&x^4 - 6x^2 - 4x + 2&x^8 - 4x^7 + 14x^6 - 28x^5 + 
43x^4 - 44x^3 + 30x^2 - 12x + 2\\
\cC_{40}&x^4 - 2x^3 - 3x^2 + 4x + 1&x^8 - 4x^7 + 5x^6 + 
2x^5 - 11x^4 + 4x^3 + 20x^2 - 32x + 16.\\
\end{array}
$$

\newpage

The relevant numerical values are given below, where 
$\mu$ is the expression $2^{-2d}\zeta_k(-1)L_{\ell|k}(-2).$
\vskip1mm
$$\begin{array}{cccccc}
{\ \ \ (k,\ell)\ \ \ }&D_k&D_\ell&\zeta_k(-1)\ \ &L_{\ell|k}(-2)\ \ &\mu\\
\cC_1&5&125&1/30&4/5&1/600\\
\cC_2&5&225&1/30&32/9&1/135\\
\cC_3&5&400&1/30&15&1/2^{5}\\
\cC_4&5&1025&1/30&160&1/3\\
\cC_5&5&1225&1/30&1728/7&18/35\\
\cC_6&5&1525&1/30&420&7/8\\
\cC_7&5&1600&1/30&474&79/2^4\cdot5\\
\cC_8&8&256&1/12&3/2&1/2^7\\
\cC_9&8&576&1/12&92/9&23/2^4\cdot3^3\\
\cC_{10}&8&1088&1/12&64&1/3\\
\cC_{11}&12&144&1/6&1/9&1/2^5\cdot3^3\\
\cC_{12}&12&2304&1/6&138&23/2^4\\
\cC_{13}&13&1521&1/6&352/9&11/3^3\\
\cC_{14}&13&2197&1/6&1332/13&111/104\\
\cC_{15}&17&2312&1/3&64&4/3\\
\cC_{16}&17&2601&1/3&536/9&67/54\\
\cC_{17}&21&441&1/3&32/63&2/189\\
\cC_{18}&24&576&1/2&2/3&1/48\\
\cC_{19}&24&2304&1/2&23&23/32\\
\cC_{20}&28&784&2/3&8/7&1/21\\
\cC_{21}&33&1089&1&4/3&1/12\\
\cC_{22}&44&1936&7/6&3&7/32\\
\cC_{23}&56&3136&5/3&48/7&5/7\\
\cC_{24}&57&3249&7/3&44/9&77/108\\
\cC_{25}&60&3600&2&60/9&5/6\\
\cC_{26}&60&3600&2&8&1\\
\cC_{27}&69&4761&2&32/3&4/3\\
\cC_{28}&76&5776&19/6&11&209/96\\
\cC_{29}&77&5929&2&96/7&12/7\\
\cC_{30}&88&7744&23/6&18&69/16\\
\cC_{31}&49&16807&-1/21&-64/7&1/147\\
\cC_{32}&49&64827&-1/21&-2408/9&43/2^3\cdot3^3\\
\cC_{33}&81&19683&-1/9&-104/27
&13/2^3\cdot3^5\\
\cC_{34}&1125&1265625&4/15&128/45&2/3^3\cdot5^2\\
\cC_{35}&2000&4000000&2/3&
12&1/2^5\\
\cC_{36}&2048&16777216&5/6&411&5\cdot137/2^9\\
\cC_{37}&2304&5308416&1&46/3&23/2^7\cdot3\\
\cC_{38}&3600&12960000&8/5&160/3&1/3\\
\cC_{39}&4352&18939904&8/3&96&1\\
\cC_{40}&4752&22581504&8/3&928/9&29/27.\\
\end{array}
$$
\newpage

\ni{\bf 8.3 Remark} The second table above lists the values of $\zeta_k(-1)$ and $L_{\ell|k}(-2)$.  These were obtained with the help of PARI/GP and the functional equations
$$\zeta_k(2) = (-2)^d\pi^{2d}D_k^{-3/2}\zeta_k(-1),\ \ \ L_{\ell|k}(3)= (-2)^d\pi^{3d}(D_k/D_{\ell})^{5/2}L_{\ell|k}(-2).$$
The values have been rechecked using MAGMA.  The latter software gives us precision up to
more than 
40 decimal places.  On the other hand, we know from
a result of Siegel [Si] that both $\zeta_k(-1)$ and $L_{\ell|k}(-2)$ 
are rational numbers.  Furthermore, the denominator of $\zeta_k(-1)$
can be effectively estimated as explained in [Si].  Similar estimates
for $L_{\ell|k}(-2)$ are given in [Ts].  In this way, we know that
the values listed in the above table are exact.  Alternatively, the
values can also be obtained from the formulae in [Si] and [Ts],
but the computations are quite tedious.

\vskip2mm

Using Proposition 2.12, and the value of $\mu$ given in the second table of 8.2,  we conclude the following at once.
\vskip2mm

\ni{\bf 8.4.} {\it The pair $(k,\ell)$, with degree $d= [k: \bQ ]>1$, can only be one of the
  following fifteen: 
$\cC_1$, $\cC_2,$ $\cC_3$, $\cC_4,$ $\cC_8$, $\cC_{10}$, $\cC_{11}$, 
$\cC_{18}$, $\cC_{20}$, $\cC_{21}$, $\cC_{26}$, $\cC_{31}$,
  $\cC_{35}$, $\cC_{38}$ and $\cC_{39}$.}
\vskip2mm

It is convenient to have the
following concrete description provided to us by Tim Steger of the
fifteen 
pairs occurring above. As before, in the sequel, $\zeta_n$ will denote a primitive
$n$-th root of unity. 
$$\begin{array}{ll}
\cC_1 = (\bQ(\sqrt{5}),\bQ(\zeta_5)), 
&\cC_2 =
(\bQ(\sqrt{5}),\bQ(\sqrt{5},\zeta_3)),\\ 
\cC_3 =(\bQ(\sqrt{5}),\bQ(\sqrt{5},\zeta_4)),&
\cC_4 = (\bQ(\sqrt{5}),\bQ(\sqrt{(-13+\sqrt{5})/2)}),\\
\cC_8=(\bQ(\sqrt{2}),\bQ(\zeta_8)),& \cC_{10} =
(\bQ(\sqrt{2}),\bQ(\sqrt{-7+4\sqrt{2}})),\\
\cC_{11} = (\bQ(\sqrt{3}),\bQ(\zeta_{12})),&
\cC_{18}=(\bQ(\sqrt{6}),\bQ(\sqrt{6},\zeta_3)),\\
\cC_{20}=
(\bQ(\sqrt{7}),\bQ(\sqrt{7},\zeta_4)),&
\cC_{21} = (\bQ(\sqrt{33}),\bQ(\sqrt{33},\zeta_3)),\\
\cC_{26}= (\bQ(\sqrt{15}),\bQ(\sqrt{15},\zeta_4)),& \cC_{31}=
(\bQ(\zeta_7+\zeta_7^{-1}), \bQ(\zeta_7)),\\
\cC_{35} =(\bQ(\zeta_{20}+\zeta_{20}^{-1}),\bQ(\zeta_{20})),& 
\cC_{38}= (\bQ(\sqrt{3},\sqrt{5}),\bQ(\sqrt{3},\sqrt{5},\zeta_4)),\\ 
\cC_{39}=(\bQ(\sqrt{5+2\sqrt{2}}), \bQ(\sqrt{5+2\sqrt{2}},\zeta_4)).&
\end{array}
$$
\vskip1mm

{\it The class number of $\ell$ in all the above pairs except $\cC_{26}$ is $1$ and the class number of $\ell$ in $\cC_{26}$ is $2.$}

\vskip2mm
 
\ni{\bf 8.5.} We will now assume that the pair $(k,\ell)$ is one of the fifteen listed above; $\cD$ and the $k$-group $G$ be as in 1.2. Let 
$(P_v)_{v\in V_f}$ be a coherent collection of parahoric subgroups  $P_v$ of
$G(k_v)$. Let $\Lambda = G(k)\cap\prod_{v\in V_f}P_v$, and $\Gamma$ be its normalizer in $G(k_{v_o})$. Let  $\cT$  be the set of nonarchimedean places $v$ of $k$ such that $P_v$ is not maximal, and also all those $v$ which are unramified in $\ell$ and $P_v$ is {\it not} a hyperspecial parahoric subgroup. Let $\cT_0$ be the subset of $\cT$ consisting of places where $G$ is anisotropic. The places in $\cT_0$ split in $\ell$, cf.\,2.2.
\vskip1mm

We first treat the case where $\cD$ is a cubic division algebra. In this case, $\cT_0$ is nonempty. 
\vskip2mm

\ni{\bf 8.6.\:Proposition.} {\it Assume that $\cD$ is a cubic division algebra. If the orbifold Euler-Poincar\'e characteristic $\chi(\Gamma)$ of $\Gamma$ is a reciprocal integer, then the pair $(k,\ell)$ must be one of the following nine: $\cC_2$, $\cC_3$, $\cC_{10}$, $\cC_{18}$, $\cC_{20}$, $\cC_{26}$, $\cC_{31}$, $\cC_{35}$ and $\cC_{39}$. Moreover, $\cT_0$ consists of exactly one place $\mathfrak{v}$, and $\cT = \cT_0$ except in the case where $(k,\ell)$ is either $\cC_{18}$ or $\cC_{20}$. Except for the pairs $\cC_3$, $\cC_{18}$ and $\cC_{35}$, $\mathfrak{v}$ is the unique place of $k$ lying over $2$; for $\cC_3$ and $\cC_{35}$, it is the unique place of $k$ lying over $5$, and for $\cC_{18}$ it is the unique place of $k$ lying over $3$.}

 \vskip1mm

{{\it Description of the possible $\cT$ if the pair $(k,\ell)$ is either $\cC_{18}$ or $\cC_{20}$:    If  $(k,\ell)= \cC_{18} = (\bQ(\sqrt{6}),\bQ(\sqrt{6}, \zeta_3))$, the possibilities are $\cT =\cT_0=\{\mathfrak{v}\}$, and $\cT =\{\mathfrak{v},\mathfrak{v}_2\}$, where $\mathfrak{v}_2$ is the unique place of $k=\bQ(\sqrt{6})$ lying over $2$.  On the other hand, if $(k,\ell) = \cC_{20} = (\bQ(\sqrt{7}), \bQ(\sqrt{7}, \zeta_4))$, let $\mathfrak{v}'_3$ and $\mathfrak{v}''_3$ be the two places of $k =\bQ(\sqrt{7})$  lying over $3$. Then either $\cT = \cT_0 =\{\mathfrak{v}\}$, or $\cT =\{\mathfrak{v}, \mathfrak{v}'_3\}$, or $\cT = \{\mathfrak{v}, \mathfrak{v}''_3\}$.}}    

\vskip2mm
 
\ni{\it Proof.} We recall from \S\S 1 and 2 that $\chi(\Gamma) = 3
\mu(G(k_{v_o})/\Gamma)$, and $\mu(G(k_{v_o})/\Gamma) =
\mu\cdot\prod_{v\in\cT}e'(P_v)/[\Gamma : \Lambda]$, where, as before, 
$\mu = 2^{-2d}\zeta_k(-1)L_{\ell|k}(-2)$. Moreover, $e'(P_v)$ is an
integer for every $v$, and, as we have shown in 2.3, $[\Gamma :
  \Lambda]$, which is a power of $3$, is at most
$3^{1+\#\cT_0}h_{\ell, 3}\prod_{v\in \cT-\cT_0}\#\Xi_{\Theta_v}$. We
note that $h_{\ell,3} =1$ for the $\ell$ occurring in any of the fifteen pairs
$(k,\ell)$ listed in 8.4. From 2.5({\it ii}) we know that for $v\in \cT_0$, $e'(P_v) = (q_v -1)^2(q_v +1)$. Now the proposition 
can be proved by a straightforward case-by-case analysis carried out 
for each of the fifteen pairs $(k,\ell)$, keeping
in 
mind Proposition 2.12, the fact that every $v\in \cT_0$
splits in $\ell$, and making use of the values of $e'(P_v)$ and 
$\#\Xi_{\Theta_v}$ given in
2.5 and 2.2 respectively. We can show that unless (i) $(k,\ell)$ {\it is one of the following nine pairs $\cC_2$, $\cC_3$, $\cC_{10}$, $\cC_{18}$, $\cC_{20}$, $\cC_{26}$, $\cC_{31}$, $\cC_{35}$ and $\cC_{39}$}, (ii) $\cT_0$ and $\cT$ {\it are as in the proposition}, and (iii) $P_v$ {\it is maximal for all $v\in V_f$, except when the pair is $\cC_{18}$}, at least one of the following two assertions will hold:

\vskip1.5mm

\ni{$\bullet$} The numerator of $\mu\cdot \prod_{v\in \cT}e'(P_v)$ is
divisible by a prime other than $3$.

\vskip1.5mm

\ni{$\bullet$} $\mu\cdot\prod_{v\in \cT}e'(P_v)/3^{\#\cT_0}\prod_{v\in
  \cT-\cT_0}\#\Xi_{\Theta_v} >1$.    

\vskip3mm

\ni{\bf 8.7.} Let $k$, $\ell$, and $G$ be as in 1.2 with $\cD =\ell$. We assume here that $d= [k:\bQ ]>1$, 
  $h_{\ell,3} =1$, and $\ell$ contains a root $\zeta$ of unity of
  order $s$. {\it We will now show that then given
  any coherent collection $(P_v^{\mathfrak m})_{v\in V_f}$ of {\it maximal} parahoric subgroups,
  the principal arithmetic subgroup $\Lambda^{\mathfrak m} :=
  G(k)\bigcap \prod_{v\in V_f}P_v^{\mathfrak m}$ contains an element of order $s$. 
In particular, $\Lambda^{\mathfrak m}$ contains an element of order $2$.} 
%Moreover, if for every nonarchimedean place $v$ of $k$ 
%which does not split in $\ell$, $\ell_v :=\ell\otimes_k k_v$ contains a primitive cube root of unity, then $\Lambda^{\mathfrak{m}}$ contains an element of order $3$.   
(In several cases of interest, Tim Steger has shown  us an explicit element of order $2$ in $\Lambda^{\mathfrak m}$.)  
For the proof, let $Q$ be the 
  quaternion division algebra with center $k$, which is unramified at
  every nonarchimedean places of $k$, and which is ramified at
  all real places of $k$ if $d$ is even, and if $d$ is odd, it is
  ramified at all real places $v\ne v_o$. 
%  Let $\ell' = k(\omega) = k[X]/(X^2+X+1)$, where $\omega$ is a primitive cube root of unity. 
  It is obvious that as $\ell$ is a totally complex quadratic extension of $k$ it embeds  in
  $Q$. We will view $\ell$ as a field contained in $Q$ in terms of a fixed embedding, and will view $Q$ as a $\ell$-vector space of dimension $2$ 
  (the action of $\ell$ on
  $Q$ is by multiplication on the left). Then the
  reduced-norm-form 
  on $Q$ gives us an hermitian form $h_0$ on the two-dimensional $\ell$-vector
  space $Q$. Now we choose 
  $a\in k^{\times}$ so that the hermitian form $h_0 \perp \langle
  a\rangle$ is indefinite at $v_o$, and definite at all real places
  $v\ne v_o$. We may (and we do) assume that $h$ is this form, see 1.2. We
  will view $G_0 := {\rm SU}(h_0)$ as a subgroup of $G = {\rm SU}(h)$
  in terms of its natural embedding.  
  \vskip1mm
  
  Let $c_{\zeta}\in G(k)$ be the 
  element which on $Q$ acts by multiplication on the left by $\zeta$,
  and on 
the $1$-dimensional $\ell$-subspace of the above hermitian space on which the hermitian form is 
$\langle a
  \rangle$ it acts by multiplication by $\zeta^{-2}$. It is obvious
  that $c_{\zeta}$ is of order $s$, and it commutes with $G_0$. 
  \vskip1mm
  
%  We assume in this paragraph that $\ell$ does not contain a primitive cube root of unity but for every nonarchimedean place $v$ of $k$ which does not split in $\ell$, $\ell_v :=\ell\otimes_k k_v$ does contain a primitive cube root of unity. Let $c$ be the element of $G_0(k)$ which acts on $Q$ by multiplication by $\omega\, (\in \ell')$ on the right; 
%  $c$ is of order $3$ and the eigenvalues of $c$ considered as a $\ell$-linear transformation of the $2$-dimension $\ell$-vector space $Q$ are $\omega$ and $\omega^2$. The centralizer of $c$ in $Q$ is clearly $\ell'$, and hence the centralizer of $c$ in $G_0$ is the $1$-dimensional torus $S_0$ consisting of elements of norm $1$ in the torus $R_{\ell'/k}({\rm{GL}}_1)$.  The centralizer of $c$ in $G$ is then the product $S$ of $S_0$ with the $1$-dimensional torus $S_1$ consisting of elements of norm $1$ in $R_{\ell/k}({\rm{GL}}_1)$.    
%  
%  \vskip1mm
 
  As $c_{\zeta}$ is a  $k$-rational element, it lies in $P_v^{\mathfrak m}$ 
  for all but finitely many $v\in V_f$. We
  assert that for every $v\in V_f$, $c_{\zeta}$ belongs to a conjugate of
  $P_v^{\mathfrak m}$ under an element of ${\overline{G}}(k_v)$. This is clear if
  $v$ splits in $\ell$ since then the maximal parahoric subgroups of
  $G(k_v)$ form a single conjugacy class under ${\overline G}(k_v)$. On the other hand, if  $v$ 
  does not split in $\ell$, then both $G$ and $G_0$ are of rank $1$ over $k_v$, and as  $c_{\zeta}$ 
  commutes with $G_0$, it fixes pointwise the apartment corresponding to any maximal
  $k_v$-split torus of $G$ contained in $G_0$. From these observations our assertion
  follows. Now Proposition
  5.3 implies that a conjugate of $c_{\zeta}$ under an element of  ${\overline{G}}(k)$
  lies in $\Lambda^{\mathfrak m}$.      
\vskip2mm

We will now prove the following proposition in which  $\Gamma$ is as in 2.1, $\Lambda =\Gamma \cap G(k)$, and $\overline\Gamma$ is the image of $\Gamma$ in ${\overline G}(k_{v_o})$.
\vskip2mm

\ni{\bf 8.8. Proposition.} {\it If $\cD = \ell$, and $\overline\Gamma$ contains a torsion-free subgroup $\Pi$ which is cocompact in $G(k_{v_o})$ and whose Euler-Poincar\'e characteristic is $3$, 
then the pair $(k,\ell)$ can only be one of the following
  five: $\cC_1$, $\cC_8$, $\cC_{11}$, $\cC_{18}$ and
  $\cC_{21}$.}

\vskip1mm

\ni{\it Proof.} It follows from 4.1 that $d>1$, so $(k,\ell)$ can only
be one of the fifteen pairs listed in 8.4. Let $\widetilde{\Pi}$ be the inverse image of $\Pi$ in $G(k_{v_o})$. As observed in 1.3, the orbifold Euler-Poincar\'e characteristic $\chi({\widetilde\Pi})$ of $\widetilde\Pi$ is $1$, 
hence the orbifold Euler-Poincar\'e characteristic $\chi(\Gamma)$ of $\Gamma$ is a reciprocal integer. Moreover, $[\Gamma : \Lambda]$ is a power of
$3$. Let $\Lambda^{\mathfrak m}$ be a maximal principal arithmetic subgroup 
of $G(k)$ containing
$\Lambda$. From the volume formula $(11)$ we see
that $\mu(G(k_{v_o})/\Lambda^{\mathfrak m})$ is an integral multiple
$a\mu$ of $\mu =2^{-2d}\zeta_k(-1)L_{\ell|k}(-2)$. We assume now 
that $\cD = \ell$, and $(k,\ell)$ is one of the following ten pairs: 
$\cC_2$, $\cC_3$, $\cC_4$, $\cC_{10}$, $\cC_{20}$, $\cC_{26}$, $\cC_{31}$, $\cC_{35}$,
$\cC_{38}$ and $\cC_{39}$. These are the pairs appearing in 8.4
excluding the five listed in the proposition. To each of these pairs we associate a prime $p$ as follows. For all these pairs except $\cC_3$ and 
$\cC_{35}$, $p$ is $2$. For $\cC_3$ and $\cC_{35}$, $p$ is $5$. We observe that the denominator of $\mu$, for each of the ten  
pairs, is prime to the corresponding $p$.  
\vskip1mm

   We will first exclude the pair $\cC_3= (\bQ(\sqrt{5}), \bQ(\sqrt{5}, \zeta_4))$.  For $v\in V_f$, let $P_v$ be as in 2.1. Then $\Lambda = G(k)\cap\prod_{v\in V_f}P_v$. Using the volume formula $(11)$, the values of $e'(P_v)$ given in 2.5, and the value of $\mu$ given in the second table in 8.2, it is easy to see that for all $v\in V_f$, $P_v$ is a maximal parahoric subgroup of $G(k_v)$ and it is hyperspecial except when $v$ is the unique place of $k$ lying over 2 (this place ramifies in $\ell$). Hence, $\Lambda$ is a maximal principal arithmetic subgroup of $G(k)$, and $\chi(\Lambda) = 3\mu(G(k_{v_o})/\Lambda) = 3\mu$. We know from 5.4 that $[\Gamma : \Lambda] = 3$ since $h_{\ell, 3} = 1$  and $\cT_0$ is empty. Then $\chi(\Gamma) = \chi(\Lambda)/3 = \mu$. Since $\chi({\widetilde\Pi}) = 1$, the index of $\widetilde\Pi$ in $\Gamma$ is $1/\mu$, which is a power of $2$ in the case presently under consideration. As $\ell$ does not contain a primitive cube root of unity, the center of $G(k)$, and so also of $\Lambda$, is trivial, and therefore, $\Gamma = \Lambda\cdot C(k_{v_o})$, where $C( k_{v_o})$ is the center of $G(k_{v_o})$ which is a cyclic group of order $3$. We conclude from this that the image $\overline\Gamma$ of $\Gamma$ in the adjoint group ${\overline G}(k_{v_o}) = {\rm{PU}}(2,1)$ coincides with the image $\overline\Lambda$ of $\Lambda$, and the index of $\Pi$ in $\overline\Lambda$ is a power of $2$.  Cartwright and Steger [CS2] have shown that $\overline\Lambda$ contains an element of order $5$. Then any subgroup of $\overline\Lambda$ of index a power of $2$, in particular, $\Pi$, contains an element of order $5$. This contradicts the fact that $\Pi$ is torsion-free.  This shows that the pair $(k, \ell)$ cannot be $\cC_3$. 
\vskip1mm   

We will now use the result proved in 8.7 to exclude the remaining nine pairs: $\cC_2$, $\cC_4$, $\cC_{10}$, $\cC_{20}$, 
$\cC_{26}$, $\cC_{31}$, $\cC_{35}$,
$\cC_{38}$ and $\cC_{39}$. This will prove the proposition. As $\ell$ occurring in each of these pairs contains a root of unity 
of order $p$, and  $h_{\ell,3} =1$, it follows from 8.7 that
   $\Lambda^{\mathfrak m}$ contains an element of order
   $p$. Hence, either $\Lambda$ contains an element of order $p$, or
   its index in $\Lambda^{\mathfrak m}$ is a multiple of $p$. This
   implies that either $\Gamma$ contains an element of order $p$, or
   the numerator of $\mu(G(k_{v_o})/\Gamma)=$\,\,$\mu(G(k_{v_o})
/{\Lambda^{\mathfrak m}})[\Lambda^{\mathfrak m}:\Lambda]/
[\Gamma : \Lambda]=$\,\,$a\mu\cdot [\Lambda^{\mathfrak m}:\Lambda]/
[\Gamma : \Lambda]$ is a multiple of $p$. This in turn implies that
   either $\widetilde{\Pi}$ contains an element of order $p$, or the
   numerator of $\chi(\widetilde{\Pi})= 3
   \mu(G(k_{v_o})/{\widetilde{\Pi}})$ is a multiple of $p$. Both these
   alternatives are impossible, the former because any element of
   finite order in $\widetilde{\Pi}$ is of order $3$, whereas $p=2$ or
   $5$, and the latter because $\chi (\widetilde{\Pi}) =1$.

\vskip8mm
\ni
\begin{center}
{\bf 9. Ten additional classes of fake projective planes}
\end{center}

\vskip3mm
              
\ni{\bf 9.1.} In this section,  $(k,\ell)$ will be one of the following nine pairs (see
Proposition 8.6).
$$\begin{array}{ll}
\cC_2 =(\bQ(\sqrt{5}),\bQ(\sqrt{5},\zeta_3)),&\cC_3 =(\bQ(\sqrt{5}),\bQ(\sqrt{5},\zeta_4)),\\
 \cC_{10}= (\bQ(\sqrt{2}),\bQ(\sqrt{-7+4\sqrt{2}})),&\cC_{18} = (\bQ(\sqrt{6}),\bQ(\sqrt{6},\zeta_3)),\\ 
  \cC_{20}=
(\bQ(\sqrt{7}),\bQ(\sqrt{7},\zeta_4)),&\cC_{26}= (\bQ(\sqrt{15}),\bQ(\sqrt{15},\zeta_4)),\\
\cC_{31} =(\bQ(\zeta_7+\zeta_7^{-1}),\bQ(\zeta_7)),& \cC_{35} =(\bQ(\zeta_{20}+\zeta_{20}^{-1}),\bQ(\zeta_{20})),\\
 \cC_{39}= (\bQ(\sqrt{5+2\sqrt{2}}),\bQ(\sqrt{5+2\sqrt{2}}, \zeta_4)).
\end{array}$$

Let $\mathfrak v$ be the unique place of $k$ lying over $p :=2$ if $(k,\ell)\ne
\cC_3,\,\cC_{18}$ and $\cC_{35}$; if $(k,\ell) =\cC_3$ or $\cC_{35}$, let $\mathfrak v$ be the unique place of $k$ lying over $p:=5$; and 
if $(k,\ell)= \cC_{18}$, let ${\mathfrak v}$ be the unique place of
$k$ lying over $p :=3$. Let $q_{\mathfrak v}$ be the cardinality of the
residue field of $k_{\mathfrak v}$. 
\vskip1mm

\ni{\bf 9.2.} Let $\cD$ be a cubic division algebra with center $\ell$ 
whose local
invariants at the two places of $\ell$ lying over $\mathfrak v$ are nonzero and negative of
each other, and whose local invariant at all the other places of $\ell$ is 
zero. There are two such division algebras, they are opposite of each
other. $k_{\mathfrak v}\otimes_k \cD= (k_{\mathfrak v}\otimes_k\ell)\otimes_{\ell}\cD={\mathfrak D}\times{\mathfrak D}^o$, where ${\mathfrak D}$ is a cubic division algebra with center $k_{\mathfrak v}$, 
and ${\mathfrak D}^o$ is its opposite. 

We fix a real place $v_o$ of $k$, and an involution $\sigma$ 
of $\cD$ of the second kind so
that $k = \{ x\in \ell\ |\ \sigma(x) = x \}$, and if $G$ is the simple simply
connected $k$-group with $$G(k)= \{ x\in \cD^{\times}\ | \ x\sigma(x)
= 1 \ \ {\rm and}\ \ {\rm Nrd}(x) =1\},$$ then  $G(k_{v_o})\cong {\rm
  SU}(2,1)$, and $G$ is anisotropic at all real places of $k$
different from $v_o$. Any other such involution of $\cD$, or of its
opposite, similarly determines a $k$-group which is $k$-isomorphic to
$G$ (1.2).  
\vskip 1mm

The set $\cT_0$ of nonarchimedean places of $k$ where $G$ is anisotropic equals $\{{\mathfrak v}\}$. As $\sigma({\mathfrak D})= {\mathfrak D}^o$, it is easily seen 
that  $G(k_{\mathfrak v})$ is the compact group 
${\rm SL}_1({\mathfrak D})$ of elements of reduced norm
$1$ in $\mathfrak D$. 
The first congruence subgroup ${\rm SL}_1^{(1)}({\mathfrak D})$
of ${\rm SL}_1({\mathfrak D})$
is known to be a pro-$p$ group, and 
${\mathfrak C} := {\rm SL}_1({\mathfrak D})/{\rm SL}_1^{(1)}({\mathfrak D})$ is a
cyclic group 
of order $(q_{\mathfrak v}^3-1)/(q_{\mathfrak v}-1)= q_{\mathfrak
  v}^2+q_{\mathfrak v}+1$, see Theorem 7(iii)(2) of [Ri].    
\vskip1mm

Let $(P_v)_{v\in V_f}$ be a coherent collection of maximal parahoric
subgroups $P_v$ of $G(k_v)$, $v\in V_f$, such that $P_v$ is hyperspecial whenever 
$G(k_v)$ contains such a
subgroup. Let 
$\Lambda = G(k)\bigcap\prod_{v\in V_f} P_v$. Let $\Gamma$ be the
normalizer of $\Lambda$ in $G(k_{v_o})$. It follows from 5.4 that 
$[\Gamma :\Lambda] = 9$ since $\# \cT_0 = 1$.  Then $\chi(\Lambda) =
3\mu(G(k_{v_o})/\Lambda) = 3\mu\cdot e'(P_{\mathfrak v})$, and (see 2.5$(ii)$) $e'(P_{\mathfrak v}) = (q_{\mathfrak v}-1)^2 (q_{\mathfrak
  v}+1)$.       

We list $q_{\mathfrak v}$, $\mu$ and $\chi(\Lambda)$ in the
table given below. 
$$\begin{array}{cccccccccc}
{(k,\ell)}&\cC_2&\cC_3&\cC_{10}&\cC_{18}&\cC_{20}&\cC_{26}&\cC_{31}&\cC_{35}&\cC_{39}\\
{q_{\mathfrak v}}&4&5&2&3&2&2&8&5&2\\
\mu&{1/135}&{1/32}&{1/3}&{1/48}&{1/21}&1&{1/147}&{1/32}&{1}\\
{\chi(\Lambda)}&1&9&3&1&{3/7}&9&9&9&9.\\
\end{array}$$

In case the pair $(k,\ell)$ is $\cC_{20} = (\bQ(\sqrt{7}), \bQ(\sqrt{7},\zeta_4))$ we will need the following three subgroups of $\Lambda$ in 9.9. Let $\fv'_3$ and $\fv''_3$ be the two places of $k= \bQ(\sqrt{7})$ lying over $3$. Note that these places do not split in $\ell =\bQ(\sqrt{7}, \zeta_4)$. We fix non-hyperspecial {\it maximal} parahoric subgroups $P'$ and $P''$ of $G(k_{\fv'_3})$ and $G(k_{\fv''_3})$ respectively. As recalled above,  there is a cubic division algebra $\fD$ with center $k_{\fv}$ such that $G(k_{\fv})$ is the compact group $\mathrm{SL}_1(\fD)$ of elements of reduced norm $1$ in $\fD$.  The first congruence subgroup $G(k_{\fv})^+ :={\mathrm{SL}}_1^{(1)}(\fD)$ of $G(k_{\fv})=\mathrm{SL}_1(\fD)$ is the unique maximal  normal pro-$p$ subgroup of $G(k_{\fv})$, and the quotient $\mathfrak{C} =G(k_{\fv})/G(k_{\fv})^+$ is of order $7$. Now  let  $\Lambda^+ = \Lambda\cap G(k_{\fv})^+$, $\Lambda' = G(k)\cap P'\cap\prod_{v\in V_f-\{\fv'_3\}}P_v$ and $\Lambda'' = G(k)\cap P''\cap \prod_{v\in V_f-\{\fv''_3\}}P_v$.  Then $\chi(\Lambda') = 3\mu (G(k_{v_o})/\Lambda') = 3 = 3\mu(G(k_{v_o})/\Lambda'')= \chi(\Lambda'')$. By the strong approximation property, $\Lambda^+$ is a subgroup of index $7$ ($= [G(k_{\fv}):G(k_{\fv})^+]$) of $\Lambda$. Hence, $\chi(\Lambda^+) = 3\mu(G(k_{v_o})/\Lambda^+) = 21\mu(G(k_{v_o})/\Lambda) =3$. 

\vskip1mm

 We will now prove the following lemma.
\vskip1mm

\ni{\bf 9.3. Lemma.} {\it Let $(k,\ell)$ be one of the nine pairs listed in 9.1. Then 

  (1) $G(k)$ is torsion-free except when
  $(k,\ell)$ is either $\cC_2$ or $\cC_{18}$ or $\cC_{20}$. 
  
  (2) If $(k,\ell)=\cC_{2}$ or $\cC_{18}$, then any
  nontrivial element of $G(k)$ of finite order is central and hence is of order
  $3$. 
  
  (3) If $(k,\ell) = \cC_{20}$, then any nontrivial element of $G(k)$ of finite order is of order $7$; $\Lambda^+$, $\Lambda'$ and $\Lambda''$ are torsion-free.}
\vskip2mm

\ni{\it Proof.} Let $x\in G(k)$ ($\subset \cD$) be a nontrivial
element of finite order, say of order $m$. As the reduced norm of $-1$
is $-1$, $-1\notin G(k)$, and so $m$ is {\it odd}.  Let $L$ be the
$\ell$-subalgebra 
of $\cD$ generated by $x$. Then $L$ is a field extension of $\ell$ of degree
$1$ or $3$. If $L =\ell$, then $x$ is clearly central, and hence it is
of order $3$. As $\ell$ does not contain a nontrivial cube root of unity
unless $(k,\ell)$ is $\cC_2$ or $\cC_{18}$, to prove the
lemma, we can assume that $L$ is an extension of $\ell$ of degree $3$. Then
$[L:\bQ ] = 6d$, where $d=2$, $3$ or $4$.   
\vskip1mm

(i)  If $(k,\ell) = \cC_2$ or $\cC_{18}$, then  $d=2$,  $[L: \bQ] = 12$, and $\zeta_3$ is 
in $\ell$. Hence, if $(k,\ell)$ is one of these two pairs, we can assume that $m$ is a multiple of
$3$. Then as $\phi (m)$, where $\phi$ is the Euler function,
must divide $12$, we conclude that $m$ is either $9$ or $21$. We assert that if $(k,\ell) = \cC_2$ or $\cC_{18}$, 
then $m=9$. For if $m=21$, then $L\cong \bQ(\zeta_{21})$, and since $3$ and $7$ are the only
primes which ramify in $\bQ(\zeta_{21})$, whereas $5$ ramifies in $k\subset L$, if $(k,\ell) = \cC_2$, so 
$m$ cannot be $21$ in this case.  Next  we observe that if $(k,\ell) = \cC_{10}$, $\cC_{18}$, or $\cC_{39}$, then as 
$7\nmid D_{\ell}$, $7$ does not ramify in $\ell$, and hence the ramification index of $L$ at $7$ is at most $3$. But the ramification index of $\bQ(\zeta_7)$ at $7$ is $6$. So if  $(k,\ell) = \cC_{10}$, $\cC_{18}$, or $\cC_{39}$, then $L$ cannot contain a nontrivial $7$-th root of unity. We conclude, in particular, that if $(k,\ell) =\cC_{18}$, $m=9$. 

Now let $(k,\ell) = \cC_2$ or $\cC_{18}$. Then, as $\ell$
  contains $\zeta_3$, and $x^3$ is of order $3$, the latter is contained in $\ell$. So any automorphism  
  of $L/\ell$ will fix $x^3$, and hence it will map $x$ to either $x$, or to $x^4$, or to 
  $x^7$. Therefore, ${\rm Nrd}(x) = x^{12} =x^3\ne 1$ and $x$ cannot belong to $G(k)$.
\vskip1mm

(ii) $(k,\ell) = \cC_3$: Then ${\rm SL}_1^{(1)}({\mathfrak D})$ is a pro-$5$ group, and $\mathfrak C$ is a group of order $31$. Since $\phi(31) =30>6d=12$, we conclude that $m$ must be a power of $5$. But  $\ell$, and hence $L$,  contains  $\zeta_4$, so $L$ contains $\zeta_{4m}$. This is impossible since $\phi(4m)$ is not a divisor of $12$.    

\vskip1mm

(iii) $(k,\ell) =\cC_{10}$, {\it or}\,\, $\cC_{31}$, {\it or}\,\, $\cC_{39}$: Then
${\rm SL}_1^{(1)}({\mathfrak D})$ is a pro-$2$ group, and $\mathfrak
C$ is a group of order $7$ if $(k,\ell) = \cC_{10}$ or $\cC_{39}$,
and is of order $73$ if $(k,\ell) =\cC_{31}$. Therefore, if
$(k,\ell)=\cC_{31}$, $m=73$, but this is impossible since $\phi(73)
=72>6d=18$.  On the other hand, if $(k,\ell) =\cC_{10}$ or $\cC_{39}$,
then $m = 7$. But this is impossible since, as we observed in (i), $L$ does not 
contain a nontrivial $7$-th root of unity.  
\vskip1mm

(iv) Let us assume now that $(k,\ell) = \cC_{20}=(\bQ(\sqrt{7}),$ $\bQ(\sqrt{7},\zeta_4))$. In this case, $L$ is of degree $12$ over $\bQ$, and as $\zeta_4\in \ell$, $L$ contains a primitive $4m$-th root of unity. This implies that $\phi(4m)$ divides $12$. From this we conclude that $m$ is either $3$, $7$ or $9$. Now since $G(k_{\fv})^+$  is a normal pro-2 subgroup of index $7$ in $G(k_{\fv})$, it is clear that the order of a nontrivial element of $G(k_{\fv})$ of odd order can only be $7$, and moreover, $G(k_{\fv})^+$ does not contain any nontrivial elements of odd order. We observe now that if $P'^+$ and $P''^+$ are the unique maximal normal pro-3 subgroups of $P'$ and $P''$ respectively, then $[P': P'^+] = 2^5\cdot 3=[P'':P''^+]$, and hence any nontrivial element of odd order of either $P'$ or $P''$ is of order $3$. Assertion (3) of the lemma follows at once from these observations.
\vskip1mm

(v) Let now $(k,\ell) = \cC_{26}=(\bQ(\sqrt{15}), \bQ(\sqrt{15},\zeta_4))$. Then again $L$ is of degree $12$ over $\bQ$, and as $\zeta_4\in\ell$, we conclude, as above, that $m$ is either $3$, $7$ or $9$. As in the case considered above, $G(k_{\fv})^+$  is a normal pro-2 subgroup of index $7$ in $G(k_{\fv})$, therefore the order of any nontrivial element of $G(k_{\fv})$ of odd order can only be $7$. This implies that $\zeta_7\in L$, and hence, $L =\bQ(\zeta_{28})$. Since the only primes which ramify in this field are $2$ and $7$, whereas $3$ ramifies in $k = \bQ(\sqrt{15}) \subset  L$, we conclude that $G(k)$ is torsion-free if $(k,\ell) =\cC_{26}$. 
\vskip1mm

(vi) Let us now consider $(k,\ell) = \cC_{35}=(\bQ(\zeta_{20}+\zeta_{20}^{-1}), \bQ(\zeta_{20}))$. In this case, $L$ is of degree $24$ over $\bQ$, and as  $L$ is an extension of degree $3$ of $\ell$, and $\zeta_{20}$ (and hence $\zeta_5$) lies in the latter, $5$ does not divide $m$. But there does not exist  such an  $m\ne 1$ for which  $\phi(20m)$ divides $24$. This implies that $G(k)$ is torsion-free.   

\vskip2mm
In the rest of this section, $\overline G$ will denote the 
adjoint group of $G$ and $\overline\Lambda$
(resp.,\,$\overline\Gamma$) the image of $\Lambda$ (resp.,\,$\Gamma$)
in ${\overline G}(k_{v_o})$. 
\vskip1mm

\ni{\bf 9.4.} {\bf Classes of fake projective planes arising from $\cC_2$
  and $\cC_{18}$ with $\cT = \cT_0$.} We assume here that $(k,\ell)$ is either $\cC_2$ or
  $\cC_{18}$, and $\cT = \cT_0$ (which is automatically the case if the pair is $\cC_2$, see Proposition 8.6). 
  Then $\ell$ contains a nontrivial cube root of unity, and
  hence the center $C(k)$  of $G(k)$ is a group of order $3$ which is contained in $\Lambda$. The naural homomorphism 
$\Lambda\rightarrow{\overline\Lambda}$ is surjective and its kernel equals 
  $C(k)$. Hence, $\chi({\overline\Lambda}) = 3\chi(\Lambda)=
  3$.  Lemma 9.3 implies that $\overline\Lambda$ is
  torsion-free. According to Theorem 15.3.1 of [Ro], $H^1(\Lambda, \bC)$
  vanishes which implies that so does $H^1({\overline\Lambda},
  \bC)$. By Poincar\'e-duality, $H^3({\overline\Lambda},\bC)$ also
  vanishes. We conclude that if $B$ is the symmetric space of
  $G(k_{v_o})$, then $B/{\overline\Lambda}$ is a fake projective
  plane. Its fundamental group is $\overline\Lambda$. There is a natural faithful action of ${\overline\Gamma}/{\overline\Lambda}$ on $B/{\overline\Lambda}$. As the normalizer of $\overline\Lambda$ in ${\overline G}(k_{v_o})$ is $\overline\Gamma$, the automorphism group of $B/{\overline\Lambda}$ equals ${\overline\Gamma}/{\overline\Lambda}$.

   Clearly, $[{\overline\Gamma} : {\overline\Lambda}]=[\Gamma : \Lambda]
   =9$.  Now let $\Pi$ be a torsion-free subgroup of 
$\overline\Gamma$ of index $9$. Then $\chi(\Pi) =3$, and so if
   $H^1(\Pi,\bC) = 0$ (or, equivalently,
   $\Pi/[\Pi,\Pi]$ is finite), then $B/\Pi$ is a fake projective
   plane, and its fundamental group is $\Pi$. The set of these fake projective planes is the class associated with $\Gamma$. For every fake projective plane belonging to this class, $\cT = \cT_0$. 
\vskip1mm

\ni{\bf 9.5. Remark.} Let $(k,\ell) = \cC_2$, and $\fD$, $\Lambda$ and $\overline\Lambda$ be as in 9.2. Then as ${\rm {SL}}_1({\mathfrak D})/{\rm {SL}}_1^{(1)}({\mathfrak D})$ is a cyclic group of order $21$, ${\rm{SL}}_1({\mathfrak D})$ contains a (unique) normal subgroup $N$ of index $3$ 
containing ${\rm{SL}}_1^{(1)}({\mathfrak D})$.  Let $\Lambda^+ = \Lambda\cap N$. Then since ${\rm{SL}}_1^{(1)}({\mathfrak D})$ is a pro-$2$ group, $\Lambda^+$ is a torsion-free normal subgroup of $\Lambda$ of index $3$. It maps isomorphically onto $\overline\Lambda$. 

\vskip2mm

 \ni{\bf 9.6.}  {\it In this subsection we will deal exclusively with $(k,\ell) =\cC_{18}$ and $\cT = \{\fv, \fv_2\}$, where $\fv_2$ is the unique place of $k= \bQ(\sqrt{6})$ lying over $2$} (see Proposition 8.6). Note that $\ell = \bQ(\sqrt{6},\zeta_3)=\bQ(\sqrt{-2},\sqrt{-3})$, the class number of $\ell$ is $1$, and $\ell_{\fv_{2}} := k_{\fv_2}\otimes_k \ell$ is an unramified field extension of $k_{\fv_{2}}$. We fix an Iwahori subgroup $I$ of $G(k_{\fv_2})$, and a {\it non-hyperspecial} {maximal} parahoric subgroup $P$ (of $G(k_{\fv_2})$) containing $I$. Let $\Lambda_P =G(k)\cap P\cap \prod_{v\in V_f-\{ \fv_2\}} P_v$ and $\Lambda_I= \Lambda_P\cap  I$.  Let $\Gamma_I$ and $\Gamma_P$ be the normalizers of $\Lambda_I$ and $\Lambda_P$ respectively in $G(k_{v_o})$. Then ${\Gamma}_I\subset {\Gamma}_P$. Let $\overline{\Lambda}_I$, $\overline{\Gamma}_I$, $\overline{\Lambda}_P$ and $\overline{\Gamma}_P$ be the images of $\Lambda_I$, $\Gamma_I$, $\Lambda_P$ and $\Gamma_P$ respectively  in $\overline{G}(k_{v_o})$. Note that $\overline{\Gamma}_P$ is
  contained in $\overline{G}(k)$, see, for example, [BP, Proposition 1.2].
\vskip1mm

 It follows from the result in 5.4 that $$[\overline{\Gamma}_I:\overline{\Lambda}_I] =[\Gamma_I : \Lambda_I] = 9 =[\Gamma_P:\Lambda_P] = [\overline{\Gamma}_P:\overline{\Lambda}_P].$$ For the pair $(k, \ell) = \cC_{18}$, using the value $\chi(\Lambda) =1$ given in 9.2, and the values $e'(I) = 9$, and $e'(P) = 3$ obtained from 2.5({\it iii}), we find that $\chi(\Lambda_I) = 9$ and $\chi(\Lambda_P) = 3$, and hence, $\chi(\overline{\Gamma}_I) = 3$ and $\chi(\overline{\Gamma}_P) = 1$. Furthermore, it follows from Theorem 15.3.1 of [Ro] that $H^1(\overline{\Gamma}_I,\bC)$, and for any subgroup $\Pi$ of $\overline{\Gamma}_P$ containing $\overline{\Lambda}_P$, $H^1(\Pi,\bC)$ vanish. 
We conclude from these observations that $\overline{\Gamma}_I$ is the fundamental group of a fake projective plane if and only if it is torsion-free, and a subgroup $\Pi$ of $\overline{\Gamma}_P$ containing $\overline{\Lambda}_P$ is the fundamental group of a fake projective plane if and only if it is torsion-free and is of index $3$ in $\overline{\Gamma}_P$.  We will now prove the following proposition.
\vskip2mm

\ni{\bf Proposition.} (i) $\overline{\Gamma}_I$ {\it is torsion-free and hence it is the fundamental group of a fake projective plane.} 
\vskip1mm

(ii) {\it There are three torsion-free subgroups of $\overline{\Gamma}_P$ containing 
$\overline{\Lambda}_P$ which are fundamental groups of fake projective planes.}
    
\vskip1mm

\ni{\it Proof.} Let $\cG$ be the connected reductive $k$-subgroup of ${\rm GL}_{1,\cD}$, which contains $G$ as a normal subgroup, such that $$\cG(k) = \{ z\in \cD^{\times}\, |\, z\sigma(z)\in k^{\times}\}.$$ Then the center $\cC$ of $\cG$ is $k$-isomorphic to $R_{\ell/k}({\rm GL}_1)$. The adjoint action of $\cG$ on the Lie algebra of $G$ induces a $k$-isomorphism $\cG/\cC\to {\overline G}$. As $H^1(k, \cC) =\{0\}$, the natural homomorphism $\cG(k)\to {\overline G}(k)$ is surjective.
\vskip1mm

Let $C$ be the center of $G$, and $\varphi: \,G\rightarrow \overline{G}$ be the natural isogeny. 
Let $\delta:\, \overline{G}(k)\rightarrow H^1(k,C)\subset \ell^{\times}/{\ell^{\times}}^3$ be the 
coboundary homomorphism. Its kernel is $\varphi(G(k))$. 
Given $\overline{g}\in \overline{G}(k)$, let $g$ be any element of $\cG(k)$ which maps onto $\overline{g}$. Then $\delta(\overline{g}) ={\rm Nrd}(g)$ modulo ${\ell^{\times}}^3$. 
\vskip1mm

Since  $\overline{\Lambda}_I$ is torsion-free (cf.\,Lemma 9.3), and $[\overline{\Gamma}_I:\overline{\Lambda}_I] = 9$, if $\overline{\Gamma}_I$ contains an element of finite order, then it contains an element $\overline{g}$ of order $3$. We fix an element $g\in \cG(k)$ which maps onto $\overline{g}$. Then $a :=g\sigma(g)\in k^{\times}$, and $\lambda : = g^3$ lies in $\ell^{\times}$. The reduced norm of $g$ is clearly $\lambda$; the norm of $\lambda$ over $k$ is $a^3 \in {k^{\times}}^3$.   Hence, the image $\delta(\overline{g})$ of $\overline{g}$ in $H^1(k,C)\,(\subset \ell^{\times}/{\ell^{\times}}^3)$ is the class of $\lambda$ in $\ell^{\times}/{\ell^{\times}}^3$.   Since $\overline{g}$ stabilizes the collection $(P_v)_{v\in V_f-\{\fv_2\}}$, as in the proof of Proposition 5.8 (cf.\:also 5.4), we conclude that $w(\lambda)\in 3\bZ$ for any normalized valuation of $\ell$ which does not lie over $2 $ or $3$. But as $\fv_2$ does not split in $\ell$, and the norm of  $\lambda$ lies
  in ${k^{\times}}^3$, it is automatic that for the normalized valuation $w$ of $\ell$ lying over $2$, $w(\lambda)\in 3\bZ$.   Therefore, $\lambda\in \ell^{\bullet}_{\{ 3\}}$, where the latter denotes the subgroup of $\ell^{\times}$ consisting of $z$ such that $N_{\ell/k}(z)\in {k^{\times}}^3$, and for all normalized valuations $w$ of $\ell$, except for the two lying over $3$, $w(z)\in 3\bZ$. Now let  $\alpha= (1+\sqrt{-2})/(1-\sqrt{-2})$. It is not difficult to see that $\ell^{\bullet}_{\{ 3\}} = \bigcup_{0\leqslant m,\, n<3}\  \alpha^m\zeta_3^n {\ell^{\times}}^3$.

\vskip1mm

Let $L$ be the field extension of $\ell$ in $\cD$ generated by $g$. Let $T$ be the centralizer of $g$ in $G$. Then $T$ is a maximal $k$-torus of $G$; its group of $k$-rational points is $L^{\times}\cap G(k)$. It can be shown that if $\lambda =g^3\in \alpha^m\zeta_3^n{\ell^{\times}}^3$, with $0\leqslant m, n <3$, then $\ell_{\fv_{2}}\otimes_{\ell}L$ is the direct product of three  copies of $\ell_{\fv_{2}}$, each stable under $\sigma$ if $n =0$, and it is an unramifield field extension of $\ell_{\fv_{2}}$ of degree $3$ if $n\ne 0$. We conclude from this that the $k$-torus $T$ is anisotropic over $k_{\fv_{2}}$.

\vskip1mm

According to the main theorem of [PY], the subset of points fixed by $g$ in the Bruhat-Tits building of $G/k_{\fv_{2}}$ is the building of $T/k_{\fv_{2}}$. But as $T$ is anisotropic over $k_{\fv_{2}}$, the building of $T/k_{\fv_{2}}$ consists of a single point. Since the two maximal parahoric subgroups of $G(k_{\fv_{2}})$ containing $I$ are nonisomorphic, if $g$ normalizes $I$, then it fixes the edge corresponding to $I$ in the Bruhat-Tits building of $G/k_{\fv_{2}}$. But as $g$ fixes just a single point in this building, we conclude that $g$ (and hence $\overline{g}$) cannot normalize $I$. This proves that $\overline{\Gamma}_I$ is torsion-free, and we have proved assertion (i) of the proposition.

\vskip1mm

We will now prove assertion (ii) of the proposition. It can be seen, using Proposition 2.9 of [BP], cf.\:5.4, that,  under the homomorphism induced by $\delta$, $\overline{\Gamma}_P/\overline{\Lambda}_P$ is isomorphic to  the subgroup $\ell^{\bullet}_{\{ 3\}}/{\ell^{\times}}^3$ of $\ell^{\times}/{\ell^{\times}}^3$. As has been noted above, $\ell^{\bullet}_{\{ 3\}} = \bigcup_{0\leqslant m, \, n<3}\  \alpha^m\zeta_3^n {\ell^{\times}}^3$, and hence, $\ell^{\bullet}_{\{ 3\}}/{\ell^{\times}}^3$ is isomorphic to $\bZ/3\bZ\times \bZ/3\bZ$. There are three subgroups of $\ell^{\bullet}_{\{ 3\}} /{\ell^{\times}}^3$ of index $3$ generated by an element of the form $\alpha^m\zeta^n_3$ with $n\ne 0$. Let $\Pi$ be the inverse image in $\overline{\Gamma}_P$ of any of these three subgroups. Then, as we will show presently, $\Pi$ is torsion-free and so it is the fundamental group of a fake projective plane.

\vskip1mm

Let us assume that $\Pi$ contains a nontrivial element $\overline{g}$ of finite order. Since $\overline{\Lambda}_P$ is torsion-free (cf.\:Lemma 9.3), and $[\Pi :\overline{\Lambda}_P] = 3$, the order of $\overline{g}$ is $3$.  As in the  proof of assertion (i), we fix $g\in \cG(k)$ which maps onto $\overline{g}$, and let $\lambda = g^3$. Then $\lambda$ is the reduced norm of $g$ and it lies in $\ell^{\bullet}_{\{3\}}$. The image $\delta(\overline{g})$ of $\overline{g}$ in $\ell^{\bullet}_{\{3\}}/{\ell^{\times}}^3$ is the class of $\lambda$ modulo ${\ell^{\times}}^3$. 
Since $\Pi$ is the inverse image in $\overline{\Gamma}_P$ of the subgroup generated by $\alpha^m\zeta^n_3$ for some $m,n<3$, with $n\ne 0$, and $\lambda$ is not a cube in $\ell$, $\lambda\in  (\alpha^m\zeta_3^n){\ell^{\times}}^3\cup  (\alpha^m\zeta_3^n)^2{\ell^{\times}}^3 $. Let $L$ be the field extension of $\ell$ in $\cD$ generated by $g$, and let $T$ be the centralizer of $g$ in $G$. Then $T(k) = L^{\times}\cap G(k)$. As observed in the proof of assertion (i), $T$ is a maximal $k$-torus of $G$ which is anisotropic over $k_{\fv_2}$, and its splitting field over $k_{\fv_2}$ is clearly $\ell_{\fv_2}\otimes_{\ell} L$ which is an unramified field extension of $\ell_{\fv_2}$ of degree $3$.  This implies that the unique point in the Bruhat-Tits building of $G/k_{\fv_2}$ fixed by $g$ is hyperspecial. But since $P$ is a 
non-hyperspecial maximal parahoric subgroup of $G(k_{\fv_2})$, it cannot be normalized by $g$. This implies that $\overline{g}$ does not lie in $\overline{\Gamma}_P$, and we have arrived at a contradiction.

\vskip2mm

\ni{\bf 9.7. Remark.}  The above proposition implies that the pair $\cC_{18}$  gives two classes of fake projective planes
with $\cT = \{ \fv, \fv_2\}$:\:the class consisting of a unique fake projective plane with the fundamental group isomorphic to $\overline{\Gamma}_I$, and the class consisting of the fake projective planes whose fundamental group is embeddable in $\overline{\Gamma}_P$, but not in $\overline{\Gamma}_I$. Cartwright and Steger [CS1] have shown that the latter class consists of just three fake projective planes up to isometry (hence, six up to biholomorphism), the ones with the fundamental group as in (ii) of the above proposition.
 
 \vskip2mm

\ni{\bf 9.8.} {\bf The classes of fake projective planes arising from the pair  
$\cC_{10}$.} We now assume that $(k,\ell) = \cC_{10}$. Then $\Lambda$
is torsion-free (9.3). Hence, ${\overline\Lambda}\cong{\Lambda}$, and
therefore, $\chi({\overline\Lambda}) =\chi(\Lambda) = 3$. Theorem
15.3.1 of [Ro] once again implies that $H^1(\Lambda,\bC)$, and so also 
$H^1({\overline\Lambda},\bC)$, vanishes. From this we conclude, as
above, that if $B$ is the symmetric space of $G(k_{v_o})$, then
$B/\Lambda$ is a fake projective plane. Its fundamental group is 
$\Lambda\cong{\overline\Lambda}$. There is a natural faithful action of ${\overline\Gamma}/{\overline\Lambda}$ on $B/{\overline\Lambda}$. As the normalizer of $\overline\Lambda$ in ${\overline G}(k_{v_o})$ is $\overline\Gamma$, the automorphism group of $B/{\overline\Lambda}$ equals ${\overline\Gamma}/{\overline\Lambda}$.

Since $[\Gamma:\Lambda] =9$, $[{\overline\Gamma}:{\overline\Lambda}] =3$, any torsion-free subgroup $\Pi$ of
$\overline\Gamma$ of index $3$ with vanishing $H^1(\Pi, \bC)$ is
the fundamental group of a fake projective plane, 
namely, that of $B/\Pi$. The set of these fake projective planes is the class associated with $\Gamma$.

\vskip2mm
\ni{\bf 9.9. Three classes of fake projective planes arising from $(k,\ell)=\cC_{20}$.}  Let $\Lambda$, 
$\Lambda^+$, $\Lambda'$ and $\Lambda''$ be as in 9.2. Let $\Gamma$, $\Gamma'$ and $\Gamma''$ be the normalizers of $\Lambda$, $\Lambda'$ and $\Lambda''$ in $G(k_{v_o})$, and $\overline{\Gamma}$, ${\overline{\Gamma}}'$ and ${\overline{\Gamma}}''$ be their images in $\overline{G}(k_{v_o})$. Let 
${\overline{\Lambda}}^+$, ${\overline{\Lambda}}'$ and ${\overline{\Lambda}}''$ be the images of 
$\Lambda^+$, $\Lambda'$ and $\Lambda''$ in $\overline{G}(k_{v_o})$. By Lemma 9.3, these groups are torsion-free.

Theorem 15.3.1 of [Ro] implies that the first cohomology (with coefficients $\bC$) of 
${\overline{\Lambda}}^+$, ${\overline{\Lambda}}'$ and ${\overline{\Lambda}}''$ vanish. As the Euler-Poincar\'e characteristic of each of these three groups is $3$, we conclude that these groups are the fundamental groups of fake projective planes $B/{\overline{\Lambda}}^+$, $B/{\overline{\Lambda}}'$ and $B/{\overline{\Lambda}}''$ respectively. The automorphism groups of these fake projective planes are respectively $\overline{\Gamma}/\overline{\Lambda}^+$, $\overline{\Gamma}'/\overline{\Lambda}' $, and $\overline{\Gamma}''/\overline{\Lambda}''$, which are of order $21$, $3$ and $3$.  Any subgroup $\Pi$ of $\overline{\Gamma}$ (resp., ${\overline{\Gamma}}'$ or ${\overline{\Gamma}}''$) of index $21$ (resp., $3$),  with vanishing $H^1(\Pi, \bC)$, is the fundamental group of a fake projective plane, namely, that of $B/\Pi$. We thus obtain three distinct classes of fake projective planes from $\cC_{20}$.

\vskip2mm

\ni{\bf 9.10.}  The constructions in 9.4, 9.6, 9.8 and 9.9  give us ten distinct classes
of fake projective planes. To see this, note that the construction is
independent of the choice of a real place of $k$ since in 9.4, 9.6, 9.8 and
9.9, $k$ is a quadratic extension of $\bQ$ and the nontrivial Galois
automorphism of $k/\bQ$ interchanges the two real places of $k$. On
the other hand, if $v$ is a nonarchimedean place of $k$ which is
unramified 
in $\ell$,
the parahoric $P_v$ involved in the construction of $\Lambda$ is 
hyperspecial, and the hyperspecial parahoric
subgroups of $G(k_v)$ are conjugate to each other under ${\overline
  G}(k_v)$, see [Ti2], 2.5. 
But if $v$ is a nonarchimedean
place of $k$ which ramifies in $\ell$, there are two possible choices of a
maximal parahoric subgroup $P_v$ of $G(k_v)$ up to conjugation. Hence, it
follows from Proposition 5.3 that each of the pairs $\cC_2$ and
$\cC_{10}$ gives two distinct classes of fake projective planes, and
the pairs  $\cC_{18}$ and $\cC_{20}$ give three each since in case $(k,\ell) = \cC_2$ or
$\cC_{10}$, there is (just) one nonarchimedean place of $k$ which
ramifies 
in $\ell$, and
if $(k,\ell) = \cC_{18}$ or $\cC_{20}$, every nonarchimedean place of $k$ is
unramified in $\ell$ since $D_{\ell} = D_k^2$.    
\vskip2mm

\ni{\bf 9.11.} We will now show that {\it the remaining five pairs $\cC_3$, $\cC_{26}$, $\cC_{31}$, $\cC_{35}$ and
$\cC_{39}$ do not give rise to any fake projective planes}. None of the fields $\ell$ occurring in these five pairs contains a nontrivial 
cube root of unity; the class number of $\ell$ is $1$ except for $\ell$ in $\cC_{26}$ which has class number  $2$.  
Let $(k,\ell)$ be one of the five pairs.  We first 
recall (9.3 and 9.2) that $\Lambda$ is a torsion-free subgroup and its
Euler-Poincar\'e characteristic is $9$. Therefore,
$\chi({\overline\Lambda})=9$. As $[\Gamma: \Lambda] =9$,
$[{\overline\Gamma}:{\overline\Lambda}] = 3$. Hence, the 
orbifold Euler-Poincar\'e characteristic
$\chi({\overline\Gamma})$ of $\overline\Gamma$ equals $3$. 
So no proper subgroup of
$\overline\Gamma$ can be the fundamental group of a fake projective
plane. We will prove presently that $\overline\Gamma$ contains an element of order $3$. This will  imply that it cannot be
the fundamental group of a fake projective plane either.           
\vskip1mm

Before embarking on the long proof of the above assertion, we give a brief outline. We construct a cubic extension $L$ 
of $\ell$ generated by an element $x$ and an involution $\tau$ of $L$ which  
restricted the subfield $\ell$ coincides with $\sigma\vert_{\ell}$ such that $x^3\in \ell$, $x\tau(x) =1$, and $(L,\tau)$ embeds 
as an $\ell$-algebra with involution into the cubic division algebra $\cD$ given with the involution $\sigma$. We identify $L$ with a maximal subfield of $\cD$ in terms of such an 
embedding and will use $\sigma$ in place of $\tau$ in the next paragraph.   

\vskip1mm

For $h\in (L^{\sigma})^{\times}$, we will define a $k$-subgroup $G_h$ of ${\rm GL}_{1,\cD}$ using the involution $\sigma_h :={\rm Int}\,h\cdot\sigma$ of $\cD$.  This subgroup is normalized by the above element $x$ so it provides a $k$-rational element of order $3$ in the adjoint group ${\overline G}_h$ of $G_h$. We will use Chebotarev's density theorem and local and global class field theory to find $h\in (L^{\sigma})^{\times}$ so that 
(1) $G_h$ is $k$-isomorphic to  $G$,  and we fix a $k$-isomorphism $\psi: G\rightarrow G_h$; (2) the $k$-rational element $g$ of ${\overline G}_h$ normalizes a coherent collection $(P'_v)_{v\in V_f}$ of parahoric subgroups $P'_v$ of $G_h(k_v)$ such that $P'_v$ is conjugate to $\psi(P_v)$ under ${\overline G}_h(k_v)$ for all $v\in V_f$, where $(P_v)_{v\in V_f}$ is the coherent collection of maximal parahoric subgroups $P_v$ of $G(k_v)$ as in 9.2. Now Proposition 5.3 implies  that a conjugate of $\psi^{-1}(g)$ in ${\overline{G}}(k)$ normalizes $(P_v)_{v\in V_f}$, and hence it normalizes $\Lambda$, and therefore lies in $\overline\Gamma$.      

\vskip1mm
We will now give a detailed proof of the above assertion. Let $\mathfrak v$ be as in 9.1, and let 
${\mathfrak v}'$ and ${\mathfrak v}''$ be the two places of $\ell$ lying over ${\mathfrak v}$. 
Recall that the cubic division algebra $\cD$ ramifies only at ${\mathfrak v}'$ and ${\mathfrak v}''$. 
Hence, $\mathfrak v$ is the only nonarchimedean place of $k$ where $G$ is anisotropic, at all the other 
nonarchimedean places of $k$ it is quasi-split. 
Let $v'$ and $v''$ be the normalized valuations of $\ell$ corresponding to ${\mathfrak v}'$ and ${\mathfrak v}''$ respectively. 

\vskip1mm
 
To find an element of ${\overline G}(k)$ of order $3$
which normalizes $\Lambda$ (and hence lies in $\overline\Gamma$) we
proceed as follows. Since the class number of $\ell$ is either $1$ or $2$, there exists an element 
$a\in \ell^{\times}$ such that $v'(a)=1$ or $2$, and for all the other normalized valuations $v$ of $\ell$, $v(a) =0$. 
Let $\lambda = a/\sigma(a)$. Then $v'(\lambda) = 1$ or $2$, $v''(\lambda ) =-v'(\lambda)$; for all normalized  valuations $v\ne v', v''$, of $\ell$, $v(\lambda) =0$,  
and $N_{\ell/k}(\lambda) =1$. We will denote the field $\ell[X]/(X^3-\lambda)$ by $L$ in the sequel, 
and $x$ will denote the unique cube root of $\lambda$ in $L$. The field $L$ admits an 
involution $\tau$ 
(i.\,e., an automorphism of order 2) whose restriction to the 
subfield $\ell$ coincides with $\sigma |_{\ell}$ and $\tau (x)=x^{-1}$. 
\vskip1mm

 We assert that there is an embedding $\iota$ of $L$ in $\cD$ such that, in
 terms of this embedding, $\sigma|_L = \tau$. Since 
$k_{\mathfrak v}\otimes_k L= (k_{\mathfrak v}\otimes_k\ell)\otimes_{\ell}L= k_{\mathfrak v}\otimes_k\ell\otimes_k L^{\tau}$ is clearly a direct product 
 of two cubic field extensions of $k_{\mathfrak v}$, $L$ does embed in
 $\cD$. Now to see that there is an embedding such that $\sigma|_L =
 \tau$, we can apply Proposition A.2 of [PrR]. The existence of local
 embeddings respecting the involutions $\sigma$ and $\tau$ need to be checked only at the real places of $k$, since
 at all the nonarchimedean places of $k$, except for $\mathfrak{v}$, $G$ is quasi-split (see
 p.\,340 of [PlR]). We will now show that for every real place $v$ of
 $k$, there is an embedding $\iota_v$ of $k_v\otimes_k L$ in $k_v\otimes_k \cD$ 
 such that $\tau =\iota_v^{-1}\sigma\iota_v\vert_L$. This will imply
 that there is an embedding $\iota$ of $L$ in $\cD$ with the desired
 property. 
\vskip1mm
 
Let $y= x+\tau(x) =x+ x^{-1}$. Then $L^{\tau}= k(y)$. As $y^3 =
x^3+x^{-3}+3(x+x^{-1})= \lambda + \sigma(\lambda)+3y$, $y^3-3y-b=0$,
where $b = \lambda+\sigma(\lambda)\in k$. The discriminant of
the 
cubic polynomial $Y^3-3Y-b$ is
$27(4-b^2)= 27\{ 4\lambda\sigma(\lambda)-(\lambda+\sigma(\lambda))^2\}
=-27(\lambda-\sigma(\lambda))^2$. Since $\ell$ is totally complex, for any real place $v$ of $k$, 
$k_v\otimes_{k} \ell= \bC$, and $\lambda-\sigma (\lambda)$ is purely imaginary. So the discriminant $-27(\lambda-\sigma(\lambda))^2$ is positive in $k_v =\bR$. Therefore, for any real place $v$ of $k$, all the roots of $Y^3-3Y-b$ are in $k_v$. This implies the following fact which will be used later in this proof.
 \vskip1mm
 
 {\it The smallest Galois extension of $k$ containing $L^{\tau}$ is totally
real, and so it is linearly disjoint from the totally complex
quadratic extension $\ell$ of $k$.} 
\vskip1mm

\ni Moreover, $k_v\otimes_k L= (k_v\otimes_k L^{\tau})\otimes_k \ell$ is a direct product of three copies of $\bC$, each of
which is stable under $\tau$. Hence there is an
embedding $\iota_v$ of $k_v\otimes_k L$ in $k_v\otimes_k \cD$, and so also an embedding $\iota$ of $L$ in $\cD$, 
with the desired property. We use $\iota$ to identify $L$ with a maximal subfield of $\cD$ and from now 
on denote the involution $\tau$ of $L$ by $\sigma$. As $\sigma(x-x^{-1}) = -(x-x^{-1})$, $(x-x^{-1})^2\in L^{\sigma}$ and $L = L^{\sigma}(x-x^{-1})$. 
\vskip1mm

We will denote the center ${\rm GL}_{1,\ell}$ of the reductive group ${\rm GL}_{1,\cD}$ by $\mathscr C$.  
Given a field extension $K$ of $k$, and $h\in (K\otimes_k L^{\sigma})^{\times}$, we denote by $\sigma_h$ the involution of $A\otimes_k\cD$, for any commutative $K$-algebra $A$,  defined as
follows $$\sigma_h (z) = {\rm Int}\,h(\sigma(z))=h\sigma(z)h^{-1} \ {\rm {for}} \ z\in A\otimes_k\cD.$$ 
Then $\sigma_1 =\sigma$ and $\sigma_h\vert_{A\otimes_k L} = \sigma\vert_{A\otimes_k L}$. Let $G_h$ (resp.,\,\,$\cG_h$) be 
the absolutely simple simply connected 
(resp.,\,\,reductive) $K$-subgroup of ${\rm GL}_{1,K\otimes_k\cD}$ such that for any commutative $K$-algebra $A$, 
$$G_h(A) =\{ z\in {\rm GL}_{1, K\otimes_k\cD}(A) = (A\otimes_k \cD)^{\times} \  |\  z \sigma_h(z)=1 \ {\rm {and}}\
{\rm {Nrd}}(z) =1\}$$ $$\cG_h(A) =\{ z\in{\rm GL}_{1,K\otimes_k\cD}(A)= (A\otimes_k\cD)^{\times}\  |
\  z\sigma_h(z)\in A^{\times}\}.$$ 

\ni $G_h$ is a normal subgroup of $\cG_h$, and the
center of $\cG_h$ equals the center $\mathscr{C}_K:={\rm GL}_{1,K\otimes_k\ell}$ of  ${\rm GL}_{1,K\otimes_k\cD}$. The conjugation action of $\cG_h$ on $G_h$  
induces a $K$-isomorphism of $\cG_h /{\mathscr C}_K$ onto the adjoint group 
${\overline G}_h$ of $G_h$. Therefore, ${\overline G}_h$ has a natural identification with a $K$-subgroup of the adjoint group 
${\rm GL}_{1,K\otimes_k\cD}/{\mathscr{C}}_K$ of ${\rm GL}_{1,K\otimes_k\cD}$. 
\vskip1mm

As above, let $x$ be the unique cube root of $\lambda$ in $L$. We will view $x$ as  the element $1\otimes x$ of $(K\otimes_kL)^{\times}\,(\subset (K\otimes_k\cD)^{\times})$. 
Since $x\sigma_h(x)=1$, $x$ is an element of $\cG_h(K)$. The image of $x$ in $({\rm GL}_{1,K\otimes_k\cD}/{\mathscr C}_K)(K)$ as well as in  
$\overline{G}_h(K)$ will be denoted by $g$. Since $x^3 =\lambda\in\ell$, $g$ is an
element of order $3$.   By the Bruhat-Tits fixed point theorem, for every nonarchimedean place $v$ of $k$ and $h\in (L^{\sigma})^{\times}$, any element of the automorphism group of $G_h$ of finite order, so in particular, $g$, normalizes a parahoric subgroup of $G_h(k_v)$.  As $x$ generates the maximal subfield $L$ of $\cD$ over its center $\ell$, and the centralizer 
of $L$ in $\cD$ is $L$,  the centralizer $\cT$ of $g$ in ${\rm GL}_{1,K\otimes_k\cD}$ is the maximal $K$-torus 
${\rm GL}_{1, K\otimes_kL}$ whose group of $A$-rational points for any commutative $K$-algebra $A$ is the subgroup $(A\otimes_k L)^{\times}$ of  $(A\otimes_k \cD)^{\times}$. 
Let $T_h$  be the centralizer of $g$ in $G_h$; then $T_h= G_h\cap \cT$ and it is a maximal $K$-torus of $G_h$. The group of $K$-rational points of $T_h$ is 
$\{z\in (K\otimes_k L)^{\times}\ |\ z\sigma(z)=1\  {\rm and}\ N_{L/\ell}(z)=1\} $.  As a torus in ${\rm GL}_{1,K\otimes_k \cD}$, $T_h$ does not depend on the choice of $h\in (K\otimes_kL^{\sigma})^{\times}$. 
So we will denote it by $T$ below. 
We shall denote $k$-groups $G_1$, ${\overline G}_1$ and ${\cG}_1$ (defined as above using $h =1$) simply by $G$, $\overline G$ and ${\cG}$ respectively.
\vskip1mm

Let $(P_v)$ be a coherent collection of parahoric subgroups $P_v$ of $G(k_v)$. Since $g$ is a $k$-rational automorphism of 
$G$ of finite order, it normalizes an arithmetic subgroup of $G(k)$ 
(for example, $\cap_i g^i\cdot\Lambda$ is an arithmetic subgroup of $G(k)$ normalized by $g$). 
By strong approximation property, the closure $U$ of such an arithmetic subgroup in the group $G(A_f)$ 
of finite ad\`eles is a compact-open subgroup normalized by $g$. So the projection $U_v$ of $U$ in $G(k_v)$ is a compact-open subgroup 
normalized by $g$. We know from 2.1 that for all but finitely many $v$, $U_v =P_v$ and $P_v$ is hyperspecial in $G(k_v)$. Thus, for all but finitely many $v$, $g$ normalizes $P_v$.

\vskip1mm

{\it In the next three paragraphs we  will show that if $v$ ramifies in $\ell$, then $g$ normalizes some conjugate of each parahoric subgroup of $G(k_v)$, and  if either 
$v\ne \mathfrak{v}$ splits in $\ell$ or it lies over $3$, then $g$ normalizes a hyperspecial parahoric subgroup of $G(k_v)$.}
\vskip1mm

Let us assume that $v$ ramifies in $\ell$, we will show in this case that  $g$ normalizes some conjugate of each parahoric subgroup of $G(k_v)$. 
For  $(k,\ell) = \cC_{26}$, $\cC_{35}$ and $\cC_{39}$, since $D_{\ell} = D_k^2$, every nonarchimedean
place of $k$ is unramified in $\ell$. If $(k,\ell) = \cC_3$, the only place of $k$ which ramifies in $\ell$ is the place 
$v$ lying over $2$, the residue field of $k_v$ has $4$ elements,  so in the Bruhat-Tits building of 
$G(k_v)$, $5$ edges emanate from every vertex. If $(k,\ell) = \cC_{31}$,
the only place $v$ of $k$ which ramifies
in $\ell$ is the place over $7$. The residue
field of $k_v$ has $7$ elements, so in the Bruhat-Tits building of
$G(k_v)$, $8$ edges emanate from every vertex. We infer that if $(k,\ell)$ is either $\cC_3$ 
or $\cC_{31}$, $g$ must fix 
at least two edges. This implies that $g$ normalizes some conjugate of each  
parahoric subgroup of $G(k_v)$ (and, in view of the main theorem of [PY], this also 
implies that the tori $T$, for the pairs $\cC_3$ and $\cC_{31}$, are isotropic over $k_v$).  
\vskip1mm

We next observe that the reduced norm of 
$x$ ($x$ considered as an element of $\cD$) is $\lambda$, and the image of $g$ in 
$H^1(k, C)\subset
\ell^{\times}/{\ell^{\times}}^3$, where $C$ is the center of $G$, is
the class of $\lambda$ in $\ell^{\times}/{\ell^{\times}}^3$. Now let $v\ne{\mathfrak v}$ be a 
nonarchimedean place of $k$ which splits in $\ell$. Then
$G(k_v)\cong {\rm SL}_3(k_v)$, and hence every maximal parahoric
subgroup of $G(k_v)$ is hyperspecial.  
As $\lambda$ is a unit in both the
embeddings of $\ell$ in $k_v$, $g$ does normalize a maximal parahoric
subgroup of $G(k_v)$, see [BP], 2.7 and 2.3(i). 
\vskip1mm

Let $v$ be a nonarchimedean place of $k$ which does not split in $\ell$, and $\ell_v := k_v\otimes_k \ell$ 
is an unramified field extension of $k_v$.  If 
$3$ does not divide $q_v+1$ (for example, if $v$ lies over $3$), 
then $g$ must normalize a hyperspecial
parahoric subgroup of $G(k_v)$. For if $g$
normalizes an Iwahori subgroup of $G(k_v)$, then it also normalizes the two nonisomorphic maximal parahoric 
subgroups of $G(k_v)$ containing this Iwahori subgroup, one of them is hyperspecial. Let us assume now that  
$g$ normalizes  a  non-hyperspecial {\it maximal} parahoric subgroup 
of $G(k_v)$. The number
of edges in the Bruhat-Tits building of
$G(k_v)$ emanating from the vertex corresponding to this parahoric 
subgroup is $q_v+1$. As $g$ is a
$k$-automorphism of $G$ of order $3$, and $3$ does not divide $q_v+1$, at least one of
these edges is fixed by $g$ and $g$ normalizes the hyperspecial
parahoric subgroup of $G(k_v)$ corresponding to the other vertex of every edge fixed by $g$.  
\vskip1mm

{\it The assertions in the next  two paragraphs hold for an arbitrary $h\in (L^{\sigma}){^\times}$ $($in particular, for $h =1)$. 
We will choose  $h$ in the fourth and the seventh paragraph below.}
 \vskip2mm
 
If $v$ is a nonarchimedean place of $k$ which does not lie over $3$, then according to the main theorem of [PY], the 
set of points fixed by $g$ in the Bruhat-Tits building of $G_h(k_v)$ is the 
Bruhat-Tits building of $T(k_v)$. Therefore, if $T$ is anisotropic at $v$ (i.\,e., $T(k_v)$ is compact), 
then as the building of  $T(k_v)$ consists of a single point, 
$g$ fixes a unique point in the building of $G_h(k_v)$. This implies that if $T$ is anisotropic at $v$ (and $v$ does not lie over $3$), 
then $g$ normalizes a unique parahoric subgroup of $G_h(k_v)$; in case $T$ splits over an unramified extension of $k_v$ 
then this parahoric subgroup is the unique parahoric subgroup of $G_h(k_v)$ containing $T(k_v)$ ([Ti2], 3.6.1). 
\vskip1mm

We assume now that $v$ is a nonarchimedean place of $k$ which does not split in $\ell$, does not lie over $3$, and  
$\ell_v := k_v\otimes_k \ell$ is an unramified field extension of $k_v$.  Then $\ell_v$ contains all the cube roots of unity, 
and {either} (i) $k_v\otimes_k L= \ell_v\otimes_{\ell}L$ is an unramified field extension of  $k_v$ in which case $k_v\otimes_k L^{\sigma}$ 
is also an unramified field extension of $k_v$, {or} (ii) $\ell_v\otimes_{\ell}L$ is a direct product of three copies of $\ell_v$ 
in which case $k_v\otimes_k L^{\sigma}$ is either the direct product of $k_v$ and $\ell_v$, or it is the direct 
product of three copies of $k_v$.  In case $k_v\otimes_k L$ is
a field, the torus $T$ is anisotropic over $k_v$ and its
splitting field is the unramified extension $k_v\otimes_k L$ of $k_v$ of degree 6. 
This implies at once that the unique parahoric subgroup
of $G_h(k_v)$ containing $T(k_v)$ is hyperspecial. This parahoric subgroup is
normalized by $g$. On the other hand, if $k_v\otimes_k L^{\sigma} = k_v\times \ell_v$, then $T$ 
is isotropic over $k_v$. The apartment in the Bruhat-Tits building of $G_h(k_v)$ corresponding to
this torus is fixed pointwise by $g$. Hence $g$ normalizes infinitely many  
hyperspecial parahoric subgroups of $G_h(k_v)$. 
\vskip1mm

Let $V$ be the set of all places $v$ of $k$ such that (i) $v$ does not lie over $3$, (ii) $\ell_v := k_v\otimes_k \ell$ is an unramified field extension of $k_v$, and (iii)  $k_v\otimes_kL^{\sigma}$ is the 
direct product of three copies of $k_v$. It may not always be the case that for every $v\in V$, $g$ normalizes a hyperspecial parahoric subgroup of $G(k_v)$. 
This is the reason why we look for and find an $h\in (L^{\sigma})^{\times}$ below such that for the $k$-group $G_h$ the following conditions hold: ({\bf 1}) $G_h(k_{v_o})$ is isomorphic to ${\rm SU}(2,1)$, 
and for all real places $v\ne v_o$ of $k$, $G_h(k_v)$ is
  isomorphic to the compact group ${\rm SU}(3)$. This condition will
  clearly hold if for every real place $v$ of $k$, $h$ is a square in $k_v\otimes_k L^{\sigma}$, 
  or, equivalently, in every embedding of $L^{\sigma}$ in $\bR$, $h$ is
  positive. {\it It will imply that the group $G_h$ defined here
  in terms of the involution $\sigma_h$ of $\cD$ is $k$-isomorphic to the
  group $G$ introduced in} 9.2\,(see 1.2).  ({\bf 2}) For every nonarchimedean place $v$ of
  $k$ such that $G_h(k_v)$ contains a hyperspecial parahoric subgroup
  (this is the case if, and only if, $v\ne {\mathfrak v}$ and $v$ does not ramify in $\ell$), $g$ normalizes one.  
\vskip1mm 

 For $v\in V$,  $k_v\otimes_k L= (k_v\otimes_k L^{\sigma})\otimes_k\ell$ is the  direct product of three copies of $\ell_v$ each of which is stable under $\sigma$. 
 This implies that for all $v\in V$, $T(k_v)$ is compact, i.e., $T$ is anisotropic over $k_v$ and it splits over the unramified extension $\ell_v$. We have shown above, 
 that for all $v\notin V$, $g$ normalizes a hyperspecial parahoric subgrop of $G(k_v)$ if $v\ne \mathfrak{v}$ and it does not ramify in $\ell$; if $v$ does ramify in $\ell$ 
 then $g$ normalizes some conjugate of any given parahoric subgroup of $G(k_v)$. Let $S$ be the finite set of nonarchimedean places $v\ne \mathfrak{v}$ of $k$ such that $v$ does not ramify in $\ell$ and  
$g$ {\it does not} normalize any hyperspecial parahoric subgroup of $G(k_v)$. Then $S\subset V$. 
If $S$ is empty, then we take $h =1$. 

\vskip1mm
We will assume now that $S$ is not empty.  
For each $v\in S$, we fix an isomorphism $\phi_v$ of $k_v\otimes_{k}\cD$ with the matrix algebra $M_3(\ell_v)$ 
which maps $k_v\otimes_kL= \ell_v\otimes_{\ell}L$ onto the subalgebra of diagonal matrices. The idempotents in $\ell_v\otimes_{\ell}L$ are all contained in $k_v\otimes_{k}L^{\sigma}$ and generate it. The idempotents of the algebra of diagonal matrices in $M_3(\ell_v)$ have entries $0$ or $1$ (hence these idempotents lie in $M_3(k_v)$). Therefore, $\phi_v(k_v\otimes_kL^{\sigma})$ is the algebra of diagonal matrices with entries in $k_v$. Now as the involution $\sigma$ is identity on $k_v\otimes_kL^{\sigma}$ we conclude that the involution of $M_3(\ell_v)$ induced from the involution $\sigma$ under $\phi_v$ is of the form ${\rm Int}\, d_v\cdot\tau_v$, where $d_v$ is a diagonal matrix in $M_3(k_v)$ and $\tau_v$ is the standard involution of the second kind  of $M_3(\ell_v)$ over $k_v$.  Let $a_v = \phi_v^{-1}(d_v)\in k_v\otimes_kL^{\sigma}$.    
 
\vskip1mm
As the smallest
Galois extension of $k$ containing $L^{\sigma}$ is linearly disjoint
from $\ell$ over $k$ (see p.\,47), using Chebotarev's density theorem we see that there are infinitely 
many nonarchimedean places $w$ of $k$ such
that $k_w\otimes_k L\,(= k_w\otimes_k
L^{\sigma}\otimes_k \ell)$ is an unramified field extension of $k_w$ of
degree $6$. As we saw above, for any such $w$ not lying over 3, and any $h\in (L^{\sigma})^{\times}$, 
$g$ normalizes a hyperspecial parahoric subgroup of $G_h(k_w)$. 
We now fix one such $w$ which does not lie over 3. We shall denote the unique extension of $w$ to $L$, as well as to $L^{\sigma}$,  
by $w$. From the definition of $V$ it is clear that $w\notin V$. 
The field $L_w = k_w\otimes_k L$ is a quadratic extension of the subfield
$L^{\sigma}_w =k_w\otimes_k L^{\sigma}$. Hence, by local class field theory, 
$N_{\ell/k}(L_w^{\times})$ is a subgroup of index $2$ of 
$(L_w^{\sigma})^{\times}$. Also, since $L$ is a quadratic extension of $L^{\sigma}$, by global
class field theory  
$N_{L/L^{\sigma}}(I_L)\cdot (L^{\sigma})^{\times}$ is a subgroup of index $2$ of $I_{L^{\sigma}}$,
where $I_L$ and $I_{L^{\sigma}}$ are the id\`ele groups of $L$ and
$L^{\sigma}$ respectively, and $N_{L/L^{\sigma}}:I_L\rightarrow
I_{L^{\sigma}}$ is the norm map.
\vskip1mm

Recall that $L = L^{\sigma}(x-x^{-1})$ and $(x-x^{-1})^2$ lies in $L^{\sigma}$. An id\`ele $c\in I_{L^{\sigma}}$ whose $v$-component $c_v$ is $1$ for $v\ne w$, and the $w$-component $c_w$ is not a norm of  any element of $L_w^{\times}$, cannot lie in the subgroup $N_{L/L^{\sigma}}(I_L)\cdot (L^{\sigma})^{\times}$ since the product $\prod_v (x-x^{-1}, c_v)_v$ of Hilbert symbols for the quadratic extension $L/L^{\sigma}$ equals $(x-x^{-1}, c_w )_w\ne 1$. Therefore, $I_{L^{\sigma}}$ is a disjoint union of $N_{L/L^{\sigma}}(I_L)\cdot (L^{\sigma})^{\times}$ and $c\cdot N_{L/L^{\sigma}}(I_L)\cdot(L^{\sigma})^{\times}$. Let $I'_L$ (resp., $I'_{L^{\sigma}}$) denote the restricted direct
product of  $(k_v\otimes_k L)^{\times}$ (resp., $(k_v\otimes_k L^{\sigma})^{\times}$), for all places $v\ne w$ of $k$.   Now considering the natural projection of $I_{L^{\sigma}}$ onto $I'_{L^{\sigma}}$,  we conclude that
$N_{L/L^{\sigma}}(I'_L)\cdot(L^{\sigma})^{\times} = I'_{L^{\sigma}}$. 
From this we see that there exists $z =(z_v)\in I'_L$, $z_v\in (k_v\otimes_k L)^{\times}$, and $h\in (L^{\sigma})^{\times}$,  
such that for $v\notin S\cup \{w\}$, the $v$-component 
$h_v :=z_v\sigma(z_v)h$ of $N_{L/L^{\sigma}}(z)h$ is 1 in $(k_v\otimes_kL^{\sigma})^{\times}$, and for $v$ in the finite set $S$, the $v$-component  
$h_v=z_v\sigma(z_v)h$ of $N_{L/L^{\sigma}}(z)h$ equals $a_v^{-1}$.  It is obvious that  $h$ is positive in every embedding of $L^{\sigma}$ in $\bR$. We choose this $h$ to define $G_h$. 
\vskip1mm

For $v\ne w$, let $G_{h_v}$ be the special unitary subgroup of the $k_v$-group ${\rm GL}_{1,k_v\otimes_k\cD}$ defined using the involution ${\rm Int}(h_v)\cdot\sigma$. It is easily seen that the conjugation action of $z_v^{-1}$ on $(k_v\otimes_k\cD)^{\times}$ gives an isomorphism $\psi_v: G_{h_v}\rightarrow G_{h}$ defined over $k_v$. As $g$ acts trivially on $k_v\otimes_kL$, $\psi_v$ commutes with the action of $g$ on $G_{h_v}$ and $G_h$. We conclude that $g$ {\it normalizes a parahoric subgroup of $G_{h_v}(k_v)$ if and only if it normalizes its image in $G_h(k_v)$ under} $\psi_v$. For $v\notin S\cup\{w\}$, since $h_v=z_v\sigma(z_v)h =1$, $G_{h_v}$ is the group $G$ defined using the involution $\sigma$. We know from the above discussion that unless $v\in S$, $g$  normalizes some  conjugate of a given parahoric subgroup of $G(k_v)$, which is assumed to be hyperspecial if $v$ does not ramify in $\ell$ and is different from $\mathfrak{v}$, so the same assertion holds for $G_h(k_v)$ for all $v\notin S\cup\{w\}$. We also know that $g$ normalizes a hyperspecial parahoric subgroup of $G_h(k_w)$ (for any $h$). 
\vskip1mm

On the other hand, for $v\in S$,  as $z_v\sigma(z_v)ha_v=1$, the involution of $M_3(\ell_v)$ induced from the involution ${\rm Int}(h_v)\cdot\sigma$ of $k_v\otimes_k\cD$ via $\phi_v$ is just the standard involution of the second kind of $M_3(\ell_v)$ over $k_v$, hence the corresponding hermitian form on $\ell_v^3$ is the standard hermitian form given by the identity matrix and $\phi_v$ maps $G_{h_v}$ isomorphically onto the special unitary group, to be denoted ${\rm SU}_3$, of this hermitian form. This isomorphism carries the centralizer $T$ of $g$ in $G_{h_v}$ 
onto the diagonal torus of ${\rm SU}_3$.  The diagonal torus of ${\rm SU}_3$  is anisotropic over $k_v$, splits over the unramified extension $\ell_v$, and the group of its $k_v$-rational points is the group of diagonal matrices of determinant $1$ whose diagonal entries are elements of $\ell_v$ of norm $1$ over $k_v$. The diagonal subgroup of 
${\rm SU}_3(k_v)$ is contained in a unique parahoric subgroup $P_v$ of ${\rm SU}_3(k_v)$ ([Ti2], 3.6.1); $P_v$ consists of matrices in ${\rm SU}_3(k_v)$ with entries in the ring of integers of $\ell_v$; $P_v$ is clearly hyperspecial. This implies that for every $v\in S$,  there is a  unique parahoric subgroup of $G_{h_v}(k_v)$ containing $T(k_v)$, and this parahoric subgroup is hyperspecial. This hyperspecial parahoric subgroup  is normalized by  $g$  as $T(k_v)$ is the centralizer of $g$ in $G_{h_v}(k_v)$. This then implies that $g$ normalizes a hyperspecial parahoric subgroup of $G_h(k_v)$ also. 
\vskip1mm

As observed above, $G_h$ is $k$-isomorphic{\footnote{To display a concrete isomorphism between $G$ and $G_h$, we observe using Landherr's theorem (see Theorem 6.27 and the remark following its proof in [PlR]) that $h = d\sigma(d)$ for a $d\in \cD^{\times}$, and conjugation action of $d$ on $\cD^{\times}$ provides a $k$-isomorphism: $G\rightarrow G_h$.}} to the group $G$ introduced in 9.2. So we will assume that 
for the five pairs of fields under consideration,  the involution $\sigma$ in 9.2 has been replaced with $\sigma_h$ for the $h$ chosen above. This amounts to replacing $G$ there with $G_h$. So we will now use the notation introduced in 9.2. In particular $(P_v)_{v\in V_f}$ is a coherent collection of maximal parahoric subgroups $P_v$ of $G(k_v)$, $v\in V_f$, such that $P_v$ is hyperspecial if $v\ne \mathfrak{v}$ and it does not ramify in $\ell$.  Since any two hyperspecial parahoric subgroups
of $G(k_v)$ are conjugate to each other under an element of
${\overline G}(k_v)$,  from the observations above we conclude
that  $g$ normalizes a coherent collection
$(P'_v)_{v\in V_f}$ of maximal parahoric subgroups such that for every
$v\in V_f$, $P'_v$ is conjugate to $P_v$  under an
element of ${\overline G}(k_v)$. Proposition 5.3 implies that a
conjugate of $g$ (in ${\overline G}(k)$) normalizes $(P_v)_{v\in
  V_f}$, and hence it normalizes $\Lambda$, and therefore lies in 
$\overline\Gamma$. This proves that $\overline\Gamma$ contains an
element of order $3$. 
\vskip2mm

Combining the results of 8.6, 9.4, and 9.6--9.11 we obtain the following.
\vskip2mm

\ni {\bf 9.12. Theorem.} {\it There exist exactly ten distinct classes of fake projective planes with the underlying totally real number field $k$ of degree $>1$, a totally complex quadratic extension $\ell$ of $k$, and a cubic division algebra $\cD$ with center $\ell$. The pair 
$(k,\ell)=\cC_2=(\bQ(\sqrt{5}),\bQ(\sqrt{5},\zeta_3))$ gives two of these ten, the pair $\cC_{10}=(\bQ(\sqrt{2}), \bQ(\sqrt{-7+4\sqrt{2} } ))$ also gives two, the pairs $\cC_{18}=(\bQ(\sqrt{6}),\bQ(\sqrt{6},\zeta_3))$ and $\cC_{20}=(\bQ(\sqrt{7}), \bQ(\sqrt{7},\zeta_4))$ give three each.}
\vskip2mm

Now combining the above theorem with the results of  8.8, the first paragraph of 1.5,  and 5.14 we obtain the following:
\vskip2mm

\ni{\bf 9.13. Theorem.} {\it  There are exactly twenty-eight non-empty classes of fake projective planes. The  underlying pairs $(k,\ell)$ and $\cT$, and the number of classes associated to each of them is given in the following tables:}
%$$\begin{array}{|l|l|c|}
%\hline
%\ \ \ (k,\ell) &{\cT}&\mbox{classes}\\
%\hline
%\hline
%
%
%\cC_2: (\bQ(\sqrt{5}),\bQ(\sqrt{5},\zeta_3))&\{\mathfrak{v}\}&2\\
%\hline
%\hline
%\cC_{10}: (\bQ(\sqrt{2}),\bQ(\sqrt{-7+4\sqrt{2}}))&\{\mathfrak{v}\}&2\\
%\hline
%\hline
%\cC_{18}: (\bQ(\sqrt{6}),\bQ(\sqrt{6},\zeta_3))&\{\mathfrak{v}\}&1\\
%\cline{2-3}
%&\{\mathfrak{v},\mathfrak{v}_2\}&2\\
%\hline
%\hline
%\cC_{20}: (\bQ(\sqrt{7}),\bQ(\sqrt{7},\zeta_4))&\{\mathfrak{v}\}&1\\
%\cline{2-3}
%&\{\mathfrak{v},\mathfrak{v}^\prime_3\}&1\\
%\cline{2-3}
%&\{\mathfrak{v},\mathfrak{v}^{\prime\prime}_3\}&1\\
%\hline
%\end{array}
%$$

$$\begin{array}{|l|l|c|}
\hline
\ \ \ (k,\ell) &{\cT}&\mbox{classes}\\
\hline
\hline
{\ (\bQ,\bQ(\sqrt{-1})) \ \ }&\{5\}&2\\ \cline{2-3}
&\{2,5\}&1\\
\hline
\hline
{ \ (\bQ,\bQ(\sqrt{-2})) \ \ }&\{3\}&2\\\cline{2-3}
&\{2,3\}&1\\
\hline
\hline
{ \ (\bQ,\bQ(\sqrt{-7})) \ \ }&\{2\}&2\\\cline{2-3}
&\{2,3\}&2\\\cline{2-3}
&\{2,5\}&2\\
\hline
\hline
{\ (\bQ,\bQ(\sqrt{-15})) \ \ }&\{2\}&4\\
\hline
\hline
{\ (\bQ,\bQ(\sqrt{-23})) \ \ }&\{2\}&2\\
\hline

\hline
\hline 

\cC_2: (\bQ(\sqrt{5}),\bQ(\sqrt{5},\zeta_3))&\{\mathfrak{v}\}&2\\
\hline
\hline
\cC_{10}: (\bQ(\sqrt{2}),\bQ(\sqrt{-7+4\sqrt{2}}))&\{\mathfrak{v}\}&2\\
\hline
\hline
\cC_{18}: (\bQ(\sqrt{6}),\bQ(\sqrt{6},\zeta_3))&\{\mathfrak{v}\}&1\\
\cline{2-3}
&\{\mathfrak{v},\mathfrak{v}_2\}&2\\
\hline
\hline
\cC_{20}: (\bQ(\sqrt{7}),\bQ(\sqrt{7},\zeta_4))&\{\mathfrak{v}\}&1\\
\cline{2-3}
&\{\mathfrak{v},\mathfrak{v}^\prime_3\}&1\\
\cline{2-3}
&\{\mathfrak{v},\mathfrak{v}^{\prime\prime}_3\}&1\\
\hline
\end{array}
$$

\vskip1cm
\ni
\begin{center}
{\bf 10. Some geometric properties of the fake projective planes}
\end{center}
\vskip4mm

In the following, $P$ will denote any fake projective plane,  
and $\Pi$ will 
denote its fundamental group. Let the pair $(k,\ell)$, the $k$-form $G$ of 
${\rm SU}(2,1)$, and the real place $v_o$ of $k$, be the ones
associated 
to $\Pi$. Let $\overline G$ be the adjoint group of $G$, $C$ the center
of $G$, and $\varphi
:\, G\rightarrow {\overline G}$ be the natural isogeny. 
Then $\Pi$ is a torsion-free cocompact arithmetic subgroup of 
${\overline G}(k_{v_o})$ ($\cong {\rm PU}(2,1)$). Let $\widetilde\Pi$ be the inverse image of 
$\Pi$ in $G(k_{v_o})$. Let $\cD$, $\Lambda$ and $\Gamma\,(\supseteq \widetilde{\Pi})$ be as in 1.3. Then 
$\Lambda =\Gamma \cap G(k)$, and $\Gamma$ is the normalizer of
$\Lambda$ in $G(k_{v_o})$. In view of the result mentioned in the first paragraph of 1.5, $\cD$ cannot be $\ell$, so it is a cubic division algebra with center $\ell$. 
\vskip2mm

\ni{\bf 10.1. Theorem.} {\it $H_1(P, \bZ) =H_1(\Pi,\bZ)={\Pi}/[{\Pi}, {\Pi}]$ is nontrivial. Therefore,  a smooth complex surface with the same integral homology groups as $\bP_{\bC}^2$ is biholomorphic to $\bP_{\bC}^2.$}

\vskip2mm

\ni{\it Proof.} There is a nonarchimedean place
$\mathfrak v$ of $k$ such that $k_{\mathfrak v}\otimes_k \cD = \mathfrak{D}\times\mathfrak{D}^o$, where $\mathfrak{D}$ is a cubic division algebra with center $k_{\mathfrak v}$, and $\mathfrak{D}^o$ is its opposite (cf.\,5.7 and 9.2). Then the group $G(k_{\mathfrak v})$ is the
compact 
group ${\rm SL}_1(\mathfrak D )$ of elements of reduced norm $1$ in $\mathfrak D$. The image $\overline\Gamma$ of $\Gamma$, and 
hence the image $\Pi$ of $\widetilde\Pi$, 
in ${\overline G}(k_{v_o})$ is contained in ${\overline G}(k)$, see 
Proposition 1.2 of [BP]. We will view $\Pi\subset {\overline G}(k)$ as a subgroup of 
${\overline G}(k_{\mathfrak v})$. We observe that ${\overline G}(k_{\mathfrak v})$ ($\cong {\mathfrak D}^{\times}/
{k_{\mathfrak v}^{\times}}$) is a pro-solvable
group, i.\,e., if we define the decreasing sequence $\{ \cG_i\}$ of subgroups of $\cG : ={\overline G}(k_{\mathfrak v})$ inductively as follows: $\cG_0 = \cG$, and $\cG_i =[\cG_{i-1},\cG_{i-1}]$, then $\bigcap \cG_i$ is trivial, to see this use [Ri], Theorem 7(i). From this it is obvious that for any subgroup $\cH$ of $\cG$, $[\cH,\cH]$ is a proper subgroup of $\cH$. We conclude, in particular, that ${\Pi}/[{\Pi}, {\Pi}]$ is nontrivial.          
\vskip2mm

\ni{\bf 10.2. Remark.} We can use the structure of ${\rm
  SL}_1({\mathfrak D})$ to provide an explicit lower bound for the order of $H_1(P,\bZ)$.
\vskip1mm

In the following proposition, $P$ is any fake projective plane whose underlying pair of number fields is neither $\cC_2 = (\bQ(\sqrt{5}), \bQ(\sqrt{5},\zeta_3))$ nor $\cC_{18} = (\bQ(\sqrt{6}), \bQ(\sqrt{6},\zeta_3))$ (these are the only pairs which give rise to fake projective planes and in which $\ell$ contains $\zeta_3$).
  
\vskip2mm

\ni{\bf 10.3. Proposition.} {\it The short exact sequence 
$$\{1\}\rightarrow C(k_{v_o})\rightarrow {\widetilde\Pi}
\rightarrow {\Pi}\rightarrow \{1\}$$ splits.}
\vskip2mm
\ni{\it Proof.} We know from 5.4 that $[\Gamma :\Lambda] =9$. As observed in the proof of the preceding theorem, the
image 
$\overline\Gamma$ of $\Gamma$, so the image $\Pi$ of $\widetilde\Pi$, 
in ${\overline G}(k_{v_o})$ is contained in ${\overline G}(k)$. Hence, $\Gamma\subset 
G({\overline{k}})$, where $\overline{k}$ is an algebraic closure of 
$k$. Now let $x$ be an element of $\Gamma$. As $\varphi (x)$ lies in 
${\overline G}(k)$, 
for every $\gamma \in {\rm Gal}({\overline {k}}/{k})$, 
$\varphi(\gamma(x))=\varphi(x)$, and hence ${\gamma(x)}x^{-1}$ lies in
$C
({\overline k})$. Therefore, $(\gamma(x)x^{-1})^3 = \gamma(x)^3 x^{-3} =1$, 
i.\,e., ${\gamma(x)}^3 = x^3$, which implies that $x^3\in \Gamma\cap G(k)
= \Lambda$.
\vskip1mm

Let $\overline\Lambda$ be the image of 
$\Lambda$ in ${\overline G}(k_{v_o})$. Then $\overline\Lambda$ 
is a normal subgroup 
of $\overline\Gamma$ of index $3$ (we have excluded the fake
projective planes arising from the pairs $\cC_2$ and $\cC_{18}$ to ensure this). Now we observe that ${\widetilde\Pi}
\cap\Lambda$ is torsion-free. This is obvious from Lemmas 5.6 and 9.3
if $\ell \ne\bQ(\sqrt{-7})$,  
since then $G(k)$, and hence $\Lambda$, is torsion-free. On the other hand, 
if $\ell = \bQ(\sqrt{-7})$, then any nontrivial element of finite order of $\Lambda$, and so 
of ${\widetilde\Pi}\cap \Lambda$, is 
of order $7$ (Lemma 5.6), but as $\Pi$ is torsion-free, the order of 
such an element must be $3$. We conclude that ${\widetilde\Pi}\cap\Lambda$ is always 
torsion-free. Therefore, it maps isomorphically onto ${\Pi}
\cap{\overline\Lambda}$. In particular, if $\Pi\subset {\overline\Lambda}$, 
then 
the subgroup ${\widetilde\Pi}\cap\Lambda$ maps isomorphically onto $\Pi$ and 
we are done. 
\vskip1mm

Let us assume now that $\Pi$ is not contained in $\overline\Lambda$. 
Then $\Pi$ projects onto ${\overline\Gamma}/{\overline\Lambda}$, which implies that ${\Pi}\cap{\overline\Lambda}$ is a normal subgroup of 
$\Pi$ of index $3$. We pick an element $g$ of 
${\Pi -\overline\Lambda}$ and let $\tilde g$ be an element of ${\widetilde\Pi}$ which 
maps onto $g$. Then ${\tilde g}^3\in {\widetilde \Pi}\cap G(k)= {\widetilde
\Pi}\cap\Lambda$, and $\bigcup_{0\leqslant i\leqslant 2}\, {\tilde g}^i({\widetilde\Pi}\cap\Lambda)$ is a subgroup of $\widetilde\Pi$ which maps isomorphically onto $\Pi$. This proves the proposition. 

\vskip2mm 

\ni{\bf 10.4.}  We note here that whenever the assertion of Proposition
10.3 holds, we get the geometric
result that {\it the canonical line
bundle $K_P$ of $P$ is three times a  holomorphic line
bundle}. To see this, we will use the following embedding of the open unit ball $B$ as an ${\rm SU}(2,1)$-orbit in $\bP^2_{\bC}$ given in Koll\'ar [Ko], 8.1. We think of ${\rm SU}(2,1)$ as the subgroup of ${\rm SL}_3(\bC)$ which keeps the hermitian form $h(x_0,x_1,x_2) = -|x_0|^2 +|x_1|^2 +|x_2|^2$ on $\bC^3$ invariant. We use the homogeneous coordinates $(x_0:x_1:x_2)$ on $\bP^2_{\bC}$. The affine plane described by $x_0\ne 0$ admits affine coordinates $z_1 = x_1/x_0$ and $z_2 = x_2/x_0$, and the open unit ball $B=\{ (z_1, z_2)\,|\,|z_1|^2+|z_2|^2<1\,\}$ in this plane is an ${\rm SU}(2,1)$-orbit. We identify $B$ with the universal cover $\tP$ of $P$. In the subgroup (of the Picard group) consisting of 
${\rm SU}(2,1)$-equivariant line bundles on $\bP^2_{\bC}$, the canonical line bundle $K_{\bP_\bC^2}$ of $\bP_\bC^2$ equals $-3H$ for
the hyperplane line bundle $H$ on $\bP_\bC^2$ ([Ko], Lemma 8.3).  Proposition 10.3 implies that 
$\Pi$ can be embedded in ${\rm SU}(2,1)$ as a discrete subgroup, and hence, 
$K_{\bP_\bC^2}|_{\tP}$ and $-H|_{\tP}$ descend to holomorphic line bundles
$K$ and $L$ on the fake projective plane $P$.  As $K=3L$ and 
$K$ is just the canonical line bundle $K_P$
of $P$, the assertion follows. 
\vskip2mm

\ni{\bf 10.5. Remark.} It follows from Theorem 3(iii) of Bombieri [B] that three times the canonical line bundle $K_P$ of $P$ is very ample, and it provides an embedding of $P$ in $\bP_{\bC}^{27}$ as a smooth surface of degree $81$.
\vskip1mm

 From the facts that (i) the second Betti number of $P$ is $1,$ 
(ii) $K_P$ is ample (since $P$ is ball-quotient), and (iii) $c_1^2=9$, we conclude as in subsection 1.1 of Chapter V of [BHPV], that there is an ample line bundle $L$ on $P$ such that $K_P = 3L \,\ {\mbox {modulo\  torsion}}$. From Theorem 1 of Reider [Re], $K+4L$ is very ample.  Kodaira
Vanishing Theorem implies that $h^i(P,K+4L)=0$ for $i>0.$
It follows from Riemann-Roch, using the Noether formula for surfaces, that 
$$h^0(P,K+4L)=\frac12 c_1(K+4L)(c_1(4L)) +\frac1{12}(c_1^2(K)+c_2(P))=15.$$ 
Let $\Phi:P\rightarrow \bP_\bC^{14}$
be the projective embedding associated to $K+4L.$  
The degree of the image is given by 
$$\deg_\Phi(P)=\int_{\Phi(P)}c_1^2(H_{\bP_\bC^{14}})
=\int_Pc_1^2(\Phi^*H_{\bP_\bC^{14}}
)=c_1^2(K+4L)=c_1^2(7L)=49.$$  Hence, holomorphic sections of $K+4L$ give an embedding 
of $P$ as a smooth surface of degree $49$ in $\bP_{\bC}^{14}$. 
\vskip2mm

\newpage

\ni
\begin{center}
{\bf Appendix: Table of class numbers}  
\end{center}
\vskip2mm

The following table lists  $(D_\ell, h_\ell, n_{\ell,3})$ for all complex quadratic extensions $\ell$ of $\bQ$ with $D_{\ell}\leqslant 79$. 

$$\begin{array}{ccccc}
(3,1,1)& (4, 1,1)& (7, 1,1)& (8, 1,1)& (11,1,1)\\
 (15,2,1)& (19, 1,1)& (20, 2,1)& (23, 3,3)&(24, 2,1)\\
 (31, 3,3)& (35, 2,1)& (39, 4,1)& (40, 2,1)& (43, 1,1)\\
 (47, 5,1)& (51, 2,1)& (52, 2,1)& (55, 4,1)& (56, 4,1)\\
 (59, 3,3)& (67, 1,1)& (68, 4,1)& (71, 7,1)& (79, 5,1).\\

\end{array}$$

\bs
\ni{\bf Acknowledgments.} We thank Shigeaki Tsuyumine for 
helpful correspondence on the values of $L$-functions. He computed 
for us the precise values of $L_{{\bQ(\sqrt{-a})}|\bQ}(-2)$ given in 3.5 and several more. We thank Pierre Deligne, Tim Dokchitser,
Igor Dolgachev, J\"urgen Kl\"uners, J\'anos Koll\'ar, 
Ron Livn\'e, Gunter Malle, Andrew Odlyzko, Chris Peters, Andrei Rapinchuk, Peter Sarnak, Tim Steger and Domingo Toledo
for their interest in this work and for many useful 
conversations, comments and correspondence. We thank Donald Cartwright for pointing out the omission of three pairs $\cC_{20}$, $\cC_{26}$, and $\cC_{35}$ in the original version of Proposition 8.6, and thank Tim Steger for pointing out that  our original assertion that $\cT =\cT_0$ is not correct for the pair $\cC_{18}$. We thank Anoop Prasad for checking our 
computations, and thank the referee whose comments have led to improvements in the exposition.

The first-named author would like to 
acknowledge support from the Institute for Advanced Study and the National Science Foundation (grant DMS-0400640).  The second-named author would also like
to acknowledge partial support from the National Science Foundation (grant 
DMS-0104089), and the hospitality of the Institute of Mathematical Research of the University of Hong Kong.

\bs
\centerline{\bf References}
\vskip3mm

\ni[BHPV] Barth, W.\,P., Hulek, K., Peters, Chris A.\,M., Van de Ven, A.:  Compact complex surfaces, Springer-Verlag, Berlin (2004).
\vskip1.5mm

\ni[B] Bombieri, E., {\it Canonical models of surfaces of general type.} Publ.
Math. IHES No. {\bf 42}(1972), 171-220.
\vskip1.5mm

\ni[BP] Borel,\,A., Prasad,\,G., {\it Finiteness theorems for discrete subgroups of bounded covolume in semisimple groups.} Publ.\,Math.\,IHES No.\,{\bf 69}(1989), 119--171.
\vskip1.5mm

\ni[BS] Borevich, Z.\,I., Shafarevich, I.\,R., {\it Number theory.} Academic Press, New York (1966).
\vskip1.5mm

\ni[CS1] Cartwright,\:D.\,I.,\:Steger,\:Tim: {\it Enumeration of the 50 fake projective planes.} C.\,R.\,Acad. Sc.\,Paris, Ser.\,I\, {\bf 348}, 11-13 (2010).
\vskip1.5mm

\ni[CS2] Cartwright,\:D.\,I.,\:Steger,\:Tim: {\it Enumerating the fake projective planes: eliminating the matrix algebra cases}, preprint.
\vskip1.5mm

\ni[Fr] Friedman, E., {\it Analytic formulas for the regulator of a number field.}
Invent.\,Math. {\bf 98}(1989), 599--622.
\vskip1.5mm

\ni[Ho] Holzapfel, R-P., {\it Ball and surface arithmetics.} Vieweg \& Sohn, Wiesbaden (1998).
\vskip1.5mm

\ni[IK] Ishida, M.-N., Kato, F., {\it The strong rigidity theorem for non-archimedean uniformization.} Tohoku Math.\,J.\,{\bf 50}(1998), 537--555.
\vskip1.5mm

\ni[Ke] Keum,\:J., {\it A fake projective plane with an order $7$ automorphism.} Topology {\bf 45}(2006), 919-927.
\vskip1.5mm

\ni[KK] Kharlamov,\:V., Kulikov,\:V.: {\it On real structres on rigid surfaces.}  Izv.\,Math. {\bf 66}, 133-150 (2002).
\vskip1.5mm

\ni [Kl] Klingler, B., {\it Sur la rigidit\'e de certains groupes fonndamentaux,
l'arithm\'eticit\'e des r\'eseaux hyperboliques complexes, et les
`faux plans projectifs'.} Invent.\,Math.\,{\bf 153} (2003), 105--143.
\vskip1.5mm
 
\ni [Ko] Koll\'ar, J., {\it Shafarevich maps and automorphic forms.} Princeton University 
Press, Princeton (1995).
\vskip1.5mm

\ni [Ma] Martinet, J.,  {\it Petits discriminants des corps de nombres.} Number theory days, 1980 (Exeter, 1980), 151--193, London Math. Soc. Lecture Note Ser., 56, Cambridge Univ. Press, Cambridge-New York, 1982.

\vskip1.5mm
\ni [Mo] Mostow, G. D., {\it Strong rigidity of locally symmetric spaces.} Annals of 
Math.\, Studies {\bf 78}, Princeton U.\,Press, Princeton (1973).  
\vskip1.5mm

\ni
[Mu] Mumford, D., {\it An algebraic surface with $K$ ample, $K^2=9$, $p_g=q=0.$}
 Amer. J.
Math.\,{\bf 101}(1979), 233--244.
\vskip1.5mm

\ni[N] Narkiewicz, W., {\it Elementary and analytic theory of algebraic numbers,}
third edition. Springer-Verlag, New York (2000).
\vskip1.5mm

\ni[O1] Odlyzko, A.\,M., {\it Some analytic estimates of class numbers and 
discriminants.} Invent. Math. {\bf 29}(1975), 275--286.
\vskip1.5mm

\ni[O2] Odlyzko, A.\,M., {\it Discriminant bounds,} unpublished, available
from:

 http://www.dtc.umn.edu/$\sim$odlyzko/unpublished/index.html.
\vskip1.5mm

\ni[PlR] Platonov, V.\,P., Rapinchuk, A.\,S., {\it Algebraic groups and Number 
theory.} Academic Press, New York (1994).
\vskip1.5mm

\ni[P] Prasad, G., {\it Volumes of $S$-arithmetic quotients of semi-simple groups.} Publ.\,Math. IHES No.\,{\bf 69}(1989), 91--117.
\vskip1.5mm

\ni[PrR] Prasad, G., Rapinchuk, A.\,S., {\it Computation of the metaplectic
kernel.} Publ. Math. IHES No.\,{\bf 84}(1996), 91--187. 
\vskip1.5mm

\ni[PY] Prasad, G., Yu, J.-K., {\it On finite group actions on reductive groups and buildings.} Invent.\,Math.\,{\bf 147}(2002), 545--560.
\vskip1.5mm

\ni[Re] Reider, I., {\it Vector bundles of rank $2$ and linear systems on algebraic
surfaces.} Ann.\,of Math.\,{\bf 127}(1988), 309-316.
\vskip1.5mm

\ni[R\'e] R\'emy,\:B., Covolume de groupes $S$-arithm\'etiques et faux plans projectifs, d'apr\'es Mumford, Prasad, Klingler, Yeung, Prasad-Yeung. S\'eminaire Bourbaki, expos\'e {\bf 984}, November 2007.
\vskip1.5mm
 
\ni[Ri] Riehm, C., {\it The norm $1$ group of $\mathfrak p$-adic division algebra.} Amer.\,J.\,Math.\,{\bf 92}(1970), 499--523. 
\vskip1.5mm
 
\ni[Ro] Rogawski, J.\,D., {\it Automorphic representations of unitary groups in three variables.} Annals of Math.\,Studies\,{\bf 123}, Princeton U.\,Press, 
Princeton (1990).
\vskip1.5mm
 
 \ni[Sc] Scharlau,W., {\it Quadratic and hermitian forms.}  Springer, Berlin (1985).
 \vskip1.5mm
 
\ni[Se1] Serre, J-P., {\it Cohomologie des groupes discrets,} in Annals of 
Math.\,Studies {\bf 70}. Princeton U.\,Press, Princeton (1971).
\vskip1.5mm
 
\ni[Se2] Serre, J-P., {\it A course in arithmetic.} Springer-Verlag, New York (1973). 
\vskip1.5mm
 
\ni[Se3] Serre, J-P., {\it Galois cohomology.} Springer-Verlag, New York (1997).
\vskip1.5mm

\ni[Si] Siegel, C.\,L., {\it Berechnung von Zetafunktionen an ganzzahligen 
Stellen.} Nachr. Akad.\,Wiss. G\"ottingen 1969, 87--102.
\vskip1.5mm

%\ni[Si2] Siegel, C.\,L., \"Uber die Fourierschen Koeffizienten
%von Modulformen.  Nachr. Akad.\,Wiss. G\"ottingen 1970, 15--56.
%\vskip1.5mm
  
\ni[Sl] Slavutskii, I.\,Sh., {\it On the Zimmert estimate for the regulator of an 
algebraic field.} English translation of Mat.\,Zametki in Math.\,Notes {\bf 51}(1992), 531--532.  
\vskip1.5mm

\ni
[Ti1] Tits, J., {\it Classification of algebraic semisimple groups. Algebraic Groups and Discontinuous
Subgroups.} Proc.\,A.M.S.\,Symp.\,Pure Math. {\bf 9}(1966) pp. 33--62. 
\vskip1.5mm

\ni
[Ti2] Tits, J., {\it Reductive groups over local fields.} Proc.\,A.M.S.\,Symp.\,Pure 
Math. {\bf 33}(1979), Part I, 29--69.
\vskip1.5mm

\ni
[Ts] Tsuyumine, S., {\it On values of $L$-functions of totally real algebraic number fields at integers,} Acta Arith. 76 (1996), no. 4, 359--392. 

%\ni [U] Ueno, K., Classification theory of algebraic varieties and 
%compact complex spaces, Lect. Notes in Math. 439, Springer-Verlag 1975.
\vskip1.5mm

\ni
[W] Washington, L. C., {\it Introduction to cyclotomic fields,} Second edition. Graduate Texts in Mathematics, 83. Springer-Verlag, New York, 1997.
\vskip1.5mm

\ni
[Y] Yeung, S.-K., {\it Integrality and arithmeticity of co-compact lattices corresponding
to certain complex two-ball quotients of 
Picard number one.} Asian J. Math.\,{\bf 8} (2004), 107--130, Erratum,
Asian J. of Math.,\,{\bf 13} (2009), 283-286.
\vskip1.5mm

\ni[Z] Zimmert, R., {\it Ideale kleiner Norm in Idealklassen und eine 
Regulatorabsch\"atzung.}  Invent.\, Math. {\bf 62}(1981), 367--380.
\vskip1.5mm

\ni
[1] The Bordeaux Database, Tables obtainable from:

    ftp://megrez.math.u-bordeaux.fr/pub/numberfields/.

\vskip5mm
\ni{\sc  
University of Michigan, Ann Arbor, MI 48109},

\ni e-mail: gprasad@umich.edu

\vskip.5cm
\ni{\sc Purdue University, West Lafayette, IN 47907}

\ni{email: yeung@math.purdue.edu}
\end{document}